%% file: amsart2.tex
\documentclass[10pt,reqno]{amsart}
\usepackage[utf8]{inputenc} 
\usepackage[T1]{fontenc} 
\usepackage[french,english]{babel} 
\usepackage{lmodern} 
\usepackage{enumerate} 
\usepackage{geometry} 
\usepackage{fancyhdr} 
\usepackage{verbatim} 
\usepackage{amsmath,amssymb, amsthm} 
\usepackage{dsfont} 
\usepackage{graphicx} 
\usepackage{tikz} 
\usepackage{float} 
\usetikzlibrary{matrix,arrows}
\usepackage{url}
\usepackage{xypic} 
\usepackage[all]{xy}
\usepackage{stmaryrd}

\usepackage{hyperref} 
\hypersetup{pdfstartview=XYZ} 

\usepackage[final]{pdfpages} 

\title{Unitary dual of p-adic $ U(5) $}

\author{Claudia Schoemann}
\date{October 2014}

\usepackage[leqno]{amsmath}     
\usepackage{amssymb}            
\usepackage{mathrsfs}           
\usepackage{dsfont}             
\usepackage{amsxtra}            

\newtheoremstyle{th}{8pt}{8pt}{\itshape}{}{\bfseries}{ ---}{5pt}{\thmname{#1}\thmnumber{ #2} \thmnote{\bfseries (#3).}}
\theoremstyle{th}

\theoremstyle{break}
\newtheorem{Theorem}{Theorem}[section]

\newtheorem{Lemma}[Theorem]{Lemma}

\newtheorem{Proposition}[Theorem]{Proposition}
\newenvironment{Proof}{\begin{proof}}{\end{proof}}

\newtheoremstyle{def}{8pt}{8pt}{}{}{\bfseries}{}{5pt}{\thmname{#1}\thmnumber{ #2} \thmnote{\bfseries (#3).}}
\theoremstyle{def}
\newtheorem{Definition}[Theorem]{Definition}

\newtheorem{Remark}[Theorem]{Remark}

\theoremstyle{nonumberplain}

\newcommand{\field}[1]{\mathds{#1}}
\newcommand{\N}{\field{N}}              

\newcommand{\Q}{\field{Q}}              
\newcommand{\C}{\field{C}}              

\newcommand{\trans}[1]{{^t\!{#1}}}
\newcommand{\Langdual}[1]{{^L\!{#1}}}

\DeclareMathOperator{\GL}{GL}

\DeclareMathOperator{\Hom}{Hom}

\DeclareMathOperator{\G}{G}

\DeclareMathOperator{\id}{id}

\DeclareMathOperator{\Ind}{Ind}

\DeclareMathOperator{\Gal}{Gal}

\DeclareMathOperator{\Lg}{Lg}

\DeclareMathOperator{\cusp}{cusp}

\DeclareMathOperator{\St}{St}

\DeclareMathOperator{\supp}{supp}

\DeclareMathOperator{\Res}{Res}

\usepackage[scaled=.90]{helvet}

\begin{document}

\include{Abstract_Article2}
\maketitle

\tableofcontents
\input{Article2}

\bibliographystyle{amsalpha}
\include{amsart2.bbl}

\bibliography{biblio}
\nocite{*}
\end{document}

%% file: Abstract_Article2.tex
\begin{abstract}
We study the parabolically induced complex representations of the unitary group in 5 variables, $ U(5), $ defined over a p-adic field.

Let $ F $ be a p-adic field. Let $ E : F $ be a field extension of degree two. Let $ \Gal(E : F ) = \{ \id, \sigma \}. $
We write $ \sigma(x) = \overline{x} \; \forall x \in E. $ Let $ E^* := E 
\setminus \{ 0 \} $ and let $ E^1 := \{x \in E \mid x \overline{x} = 1 \}. $

$ U(5) $ has three proper standard Levi subgroups, the minimal Levi subgroup $ M_0 \cong E^* \times E^* \times
E^1 $ and the two maximal Levi subgroups $ M_1 \cong \GL(2, E) \times E^1 $ and $ M_2 \cong E^* \times U(3). $

We consider representations induced from the minimal Levi subgroup $ M_0, $ representations
induced from non-cuspidal, not
fully-induced representations of the two maximal Levi subgroups $ M_1 $ and $ M_2, $
and representations induced from cuspidal representations
of $ M_1.$

We describe - except several particular cases - the unitary dual in terms of Langlands-quotients.
\end{abstract}

%% file: Article2.tex
\section{Introduction}

We study the parabolically induced complex representations of the unitary group in 5 variables - $ U(5) $ - defined over a 
non-archimedean local field of characteristic $0$. This is $ \Q_p $ or a finite extension of $ \Q_p, $ where
$ p $ is a prime number. We speak of a 'p-adic field'.

\medskip

Let $ F $ be a p-adic field. Let $ E:F $ be a field extension of degree two. Let $ \Gal(E:F) = \{ \id, \sigma \} $ be the
Galois group. We write $ \sigma(x) = \overline{x} \; \forall x \in E. $ Let $ E^* := E \setminus \{ 0 \} $
and let $ E^1 := \{ x \in E \mid  x \overline{x} = 1 \}. $ Let $ \mid \; \mid_E $ denote the p-adic norm on 
$ E. $

\medskip

$ U(5) $ has three proper parabolic subgroups. Let $ P_0 $ denote the minimal standard parabolic
subgroup and $ P_1 $ and $ P_2 $ the two maximal standard parabolic subgroups. One has the Levi decomposition
$ P_i = M_i N_i, \; i \in \{ 0, 1, 2 \}, $ where $ M_i $ denote the standard Levi subgroups and $ N_i $ are the corresponding
unipotent subgroups of $ U(5). $

\medskip

$ M_0 \cong E^* \times E^* \times E^1 $ is the minimal Levi subgroup, $ M_1 \cong \GL(2,E) \times E^1 $ and
$ M_2 \cong E^* \times U(3) $ are the two maximal Levi subgroups.

\medskip

We consider representations of the Levi subgroups, extend them trivially to the unipotent subgroups to obtain 
representations of the parabolic groups. We now perform normalised parabolic induction to obtain representations of
$ U(5). $

\medskip

In this article we give a classification of the irreducible unitary representations of $ U(5) $ in terms of Langlands 
quotients.
At first we consider the irreducible subquotients obtained by induction from representations of $ M_0 $ and from
non-cuspidal,
not fully-induced representations of $ M_1 $ and $ M_2 $ (the subquotients determined in \cite{CS}). We then
consider the irreducible subquotients
of representations induced from cuspidal representations of $ M_1 \cong \GL(2,E) \times E^1. $

\bigskip

\section{Definitions}

Let $ \G $ be a connected reductive algebraic group, defined over a non-archimedean local field of characteristic $ 0. $
Let $ V $ be a vector space, defined over the complex numbers. Let $ \pi $ be a representation of
 $ \G $ on $ V. $ We denote it by $ (\pi,V) $ and sometimes by $ \pi $ or $ V. $
\smallskip

Let $ \mathcal{K} $ denote the set of open compact subgroups of $ \G. $

\begin{Definition}
 A representation $ (\pi, V) $ of $ \G $ is said to be \textbf{smooth} if every $ v \in V $ is fixed by a neighbourhood of the
unity in $ \G. $ This is equivalent to saying that $ \exists K \in \mathcal{K} $ such that $ \pi(k)v = v \; \forall k \in 
K. $ 
\end{Definition}

\begin{Definition}
 A representation $ ( \pi, V) $ of $ \G $ is said to be \textbf{admissible} if it is smooth and for every $ K \in \mathcal{K}
$ the space
$ V^K $ of fixed vectors under $ K $ is finite dimensional.
\end{Definition}

From now on let $ (\pi, V) $ be an admissible representation of $ \G. $

\begin{Definition} Let $ (\pi, V) $ be a smooth representation of $ \G. $ Let $ V^* = \Hom(V, \C) $ be the space of
linear forms on $ V. $ One defines a representation $ (\pi^*, V^*)$ of $ \G: $ if $ v^* \in V^* $ and $ g \in G $ then 
$ \pi^*(g)(v^*) $ is defined by $ v \mapsto v^*(\pi(g^{-1})(v)). $ Let $ \overset{\sim}{V} $ be the subspace of smooth
vectors in $ V^* ,$ i.e. the subset 
$ \overset{\sim}{V} \subset 
V^* $ of elements $ \overset{\sim}{v} $ 
such that the stabiliser of $ \overset{\sim}{v} $ is an open subgroup of $ \G. $ One shows that $ \overset{\sim}{V} $
 is invariant under $ \G $ for the action of $ \pi^*. $ So $ \pi^* $ induces a representation
$ \overset{\sim}{\pi} $ on $ \overset{\sim}{V}; \; (\overset{\sim}{\pi}, \overset{\sim}{V}) $ is
called the \textbf{dual representation} of $ (\pi, V). $  

\end{Definition}

\begin{Definition} A representation $ \pi  $ is called \textbf{hermitian} if 
 $  \pi \cong  \overline{\overset{\sim}{\pi}}. $
\end{Definition}

\begin{Definition} A representation $ (\pi,V) $ of a group $ \G $ is called a \textbf{unitary representation} if and only if on the
vector space
 $ V $ there exists a positive definite hermitean form $ \langle \; , \; \rangle: V \times V \rightarrow \C $ that is
 invariant under the action of $ \G: $
\end{Definition}

$$ \langle \pi(g) v, \pi(g) w \rangle = \langle v,w \rangle \; \forall g \in \G, \forall v,w \in V.  $$

\begin{Definition} 
A \textbf{matrix coefficient} $ f_{v,\overset{\sim}{v}} $ of a representation $ (\pi,V) $ is a (locally constant) function: $ f_{v,
\overset{\sim}{v}}: \G \rightarrow \C: 
g \mapsto \overset{\sim}{v}(\pi(g)(v)), $ where $ v \in V $ and $ \overset{\sim}{v} \in \overset{\sim}{V}. $
\end{Definition}

\begin{Definition} An irreducible representation $ \pi $ of $ \G $ is called \textbf{cuspidal} if $ \pi $ has a non-zero
matrix
coefficient 
$ f: \G \rightarrow \C $ that is
compactly supported modulo the center of G.
\end{Definition} 

\begin{Definition} An irreducible representation $ \pi $ of $ \G $  is called \textbf{square-integrable} if $ \pi $ is unitary and if
$ \pi $ has a non-zero matrix coefficient
 $ f: \G \rightarrow \C $ that is square-integrable
modulo the center $ Z $ of $ \G: \; f \in L^2(\G / Z). $ It follows that all matrix 
coefficients of $ \pi $ are square-integrable.
\end{Definition}

\begin{Definition}  An irreducible representation $ \pi $ of a group $ \G $ is \textbf{tempered} if it is unitary and if
$ \pi $ has a
non-zero matrix coefficient
$ f: \G \rightarrow \C $ that verifies $ f \in L^{2 + \epsilon}(\G / Z) \; \forall \epsilon > 0. $
\end{Definition}

\begin{Remark} Square-integrable representations are tempered.
\end{Remark}

\bigskip

Let $ F $ be a non-archimedean local field of characteristic $ 0. $ i.e. $ \Q_p $ or a finite extension of $ \Q_p, $ where
$ p $ is a prime number.

Let $ E:F $ be a field extension of degree two, hence a Galois extension. Let $ \Gal(E:F) = \{ \id, \sigma \}. $
The action of the non-trivial element $ \sigma $ of the Galois group is called the conjugation of elements in $ E $ corresponding
to the extension $ E: F. $
We write $ \sigma(x) = \overline{x} \; \forall \;x \in E $ and extend $ \overset{-}{\;} $ to
matrices with entries in $ E. $

\medskip

Let $ \Phi \in \GL(n,E) $ be a hermitian matrix, i.e. $ \overline{\Phi}^t = \Phi, $ let $ U_{\Phi} $ be the unitary group
definded by $ \Phi: $

$$ U_{\Phi} = \{ g \in \GL(n,E): g \Phi \overline{g}^t = \Phi \}. $$

\bigskip

Let $ \Phi_n = (\Phi_{ij}), $ where $ \Phi_{ij} = (-1)^{i-1} \delta_{i,n+1-j} $ and $ \delta_{ab} $ is the Cronecker delta.

Let $ \zeta \in E^* $ be an element of trace $ 0, $ i.e. $ tr(\zeta) = \zeta + \overline{\zeta} = 0. $

\smallskip
If  n is odd, then
$ \Phi_n = \left(
\begin{smallmatrix}
&&&&&&1\\
&&&&&\cdotp&\\
&&&&\cdotp&&\\
&&&\cdotp&&&\\
&&1&&&&\\
&-1&&&&&\\
1&&&&&&
\end{smallmatrix}
\right) $ is hermitian. If n is even, 
$ \zeta \Phi_n = \left(
\begin{smallmatrix}
 &&&&&&\zeta\\
&&&&&\cdotp&\\
&&&&\cdotp&&\\
&&&\cdotp&&&\\
&& - \zeta &&&&\\
&\zeta&&&&&\\
- \zeta&&&&&&
\end{smallmatrix}
\right) $ is hermitian.

Denote by \textbf{$ U(n) $} the \textbf{unitary group} corresponding to $ \Phi_n $ if $ n $ is odd or to $ \zeta \Phi_n $ if $ n $ is even,
 respectively.
It is quasi split.

\bigskip

Let $ n $ be an even positive integer. We will call Levi subgroup of $ U(n) $ a subgroup of block diagonal matrices

$ L: = \{ \left(
\begin{smallmatrix}
 A_1&&&&&&-&0\\
&A_2&&&&&&\mid\\
&&\ddots&&&&&\\
&&&A_k&&&&\\
&&&&B&&&\\
&&&&&\trans{\overline{A}}_k^{-1}&&\\
\mid &&&&&&\ddots&\\
0&-&&&&&&\trans{\overline{A}}_1^{-1}
\end{smallmatrix}
\right), \; A_i \in \GL_{n_i}(E) \; \text{ for } \;  1 \leq i \leq k \; \text {and } B \in U(m) \}, $

where $ m, n_1, \ldots, n_k $ are strictly positive integers such that $ m + 2 \underset{i = 1}{\overset{k}{\sum}} n_i = n. $
(If $ k = 0, $ then there are no $ n_i $ and $ L = U(n). $)

It is canonically isomorphic to the product
$ \GL(n_1,E) \times \ldots \times \GL(n_k, E) \times U(m).$
We chose the corresponding parabolic subgroup $ P $ such that it contains $ L $ and the subgroup of upper triangular matrices
in $ U(n). $ We call a parabolic subgroup $ P $ that contains the subgroup of upper triangular matrices standard.
Let $ N $ be the unipotent subgroup with identity matrices for the block diagonal matrices of $ L, $ arbitray entries in 
$ E $ above and $ 0's $ below. Then one has the Levi decomposition $ P = L N. $ 

\smallskip

Let $ \pi_i, \; i = 1, \ldots k, $ be smooth admissible
representations of $ \GL(n_i,E) $ and $ \sigma $ a smooth admissible representation of $ U(m). $ Let 
$ \pi_1 \otimes \ldots \otimes \pi_k \otimes \sigma $ denote the representation of 
$ L = \GL(n_1,E) \times \ldots \times \GL(n_k,E)
\times U(m) $ and denote by
$ \pi: = \Ind_P^{U(n)} (\pi_1 \otimes \ldots \pi_k \otimes \sigma) = \pi_1 \times \ldots \pi_k \rtimes \sigma $ the
normalized parabolically induced representation, where $ P $ is the corresponding standard parabolic subgroup containing
$ L $.    

\begin{Definition} Let $ \pi $ be an irreducible representation of $ \GL(n,E). $ Then there exist irreducible
cuspidal representations $ \rho_1, \rho_2, \ldots, \rho_k $ of general linear groups that are, up to isomorphism,
uniquely defined by $ \pi, $ such that $ \pi  $ is
isomorphic to a subquotient of $ \rho_1 \times \ldots \times \rho_k. $ The multiset of equivalence classess $ (\rho_1,
\ldots, \rho_k) $ is called the \textbf{cuspidal support} of $ \pi. $ It is denoted by $ \supp(\pi).$
\end{Definition}

\begin{Definition}
Let $ n \in \N $ and let $ \tau $ be an irreducible representation of $ U(n). $ Then there exist irreducible cuspidal
representations $ \rho_1, \ldots, \rho_k $ of general linear groups and an irreducible cuspidal representation $ \sigma $
of some $ U(m) $ that are, up to isomorphism and replacement of $ \rho_i $ by $ \rho_i^{-1}(\overset{-}{\;}) $ for some
$ i \in \{ 1, \ldots, k \},  $ uniquely defined by $ \tau, $ s. t. $ \tau $ is isomorphic to a subquotient of
$ \rho_1 \times \ldots \times \rho_k \rtimes \sigma. $ The representation $ \sigma $ is called the \textbf{partial cuspidal
support} of $ \tau $ and is denoted by $ \tau_{\cusp}. $
\end{Definition}

\begin{Definition}
Let $ \pi $ be a smooth representation of finite length of $ G. $ Then  $\hat{\pi} $ denotes the \textbf{Aubert dual} of 
$ \pi, $ as
defined in \cite{MR1285969}.
\end{Definition}

\bigskip

\section{Irreducible unitary representations of $ U(5), $ in terms of Langlands-quotients}

\subsection{Representations with cuspidal support in $ M_0, $ fully-induced}

\bigskip

Let $ \chi_{\omega_{E/F}} \in X_{\omega_{E/F}}. $ Let $ \pi_{1, \chi_{\omega_{E/F}}} $ be the unique irreducible
 square-integrable
subquotient of
$ \mid \; \mid^{1/2} \chi_{\omega_{E/F}} \rtimes \lambda'. $ Let $ \chi_{1_{F^*}} \in X_{1_{F^*}}. $ Recall that
$ \chi_{1_{F^*}} = \sigma_{1, \chi_{1_{F^*}}} \oplus \sigma_{2, \chi_{1_{F^*}}}, $ where $ \sigma_{1, \chi_{1_{F^*}}} $ and
$ \sigma_{2, \chi_{1_{F^*}}} $ are tempered \cite{Ky}.

\begin{Proposition}
Let $ 0 < \alpha_2 \leq \alpha_1, \alpha > 0. $ Let $ \chi_1, \chi_2 $ and $ \chi $ be unitary characters of $ E^*.$ The following
list exhausts all irreducible hermitian representations of $ U(5) $ with cuspidal support in $ M_0: $

\smallskip

0. tempered representations of $ U(5), $

1. $ \Lg(\mid \; \mid^{\alpha_1} \chi_1 ; \mid \; \mid^{\alpha_2} \chi_2 ; \lambda') $ where $ \chi_1, \chi_2 \in X_{N_{E/F}(E^*)} $ or
$ \alpha_1 = \alpha_2 $ and $ \chi_1(x) = \chi_2^{-1}(\overline{x}) \; \forall x \in E^*, $

2. $ \Lg(\mid \; \mid^{\alpha} \chi_1 ; \chi_2 \rtimes \lambda') $ where $ \chi_1 \in X_{N_{E/F}(E^*)} $ and $ \chi_2 \notin X_{1_{F^*}}, $

3. $ \Lg(\mid \; \mid^{\alpha} \chi  \St_{\GL_2} ; \lambda') $ where $ \chi \in X_{N_{E/F}(E^*)}, $

4.  $ \Lg(\mid \; \mid^{\alpha} \chi ; \lambda'(\det) \St_{U(3)}), \Lg( \mid \; \mid^{\alpha} \chi, \pi_{1, \chi_{\omega_{E/F}}}),
\Lg( \mid \; \mid^{\alpha} \chi; \sigma_{1, \chi_{1_{F^*}}}) $ and $ \Lg( \mid \; \mid^{\alpha} \chi ;
\sigma_{2, \chi_{1_{F^*}}}), $ where $ \chi \in X_{N_{E/F}(E^*)}. $
\end{Proposition}

Outline of the proof: 0. Tempered representations are unitary, hence hermitian.

1.-4. Let $ \lambda_i, \; i = 0,1,2 $ be representations of the Levi subgroups $ M_0, M_1 $ and $ M_2. $ By \cite{Ca},
$ \overline{\overset{\sim}{\Ind_{M_i}^{U(5)}(\lambda_i)}} \cong \Ind_{M_i}^{U(5)}(\lambda_i), \; i = 0,1,2, $ is equivalent to
$ \exists w \in W $
such that $ \overline{\overset{\sim}{\lambda_i}} = w \lambda_i, \; i = 0,1,2. $ This holds also for the Langlands quotients,
and the proof is immediate.

\bigskip

\begin{center} Irreducible subquotients of $ \chi_1 \times \chi_2 \rtimes \lambda' $ \end{center}

\medskip

Let $ \chi_1, \chi_2 $ be unitary characters of $ E^*. $

\medskip

All irreducible subquotients of $ \chi_1 \times \chi_2 \rtimes \lambda' $ are tempered, hence unitary.
 
\bigskip

\begin{center}
$  \Lg(\mid \; \mid^{\alpha_1} \chi_1 ; \mid \; \mid^{\alpha_2} \chi_2 ; \lambda'), \; 0 < \alpha_2 \leq 
\alpha_1,
 \; \chi_1 \notin X_{N_{E/F}(E^*)} $ or $ \chi_2 \notin X_{N_{E/F}(E^*)} $ \end{center}

\medskip

\begin{Theorem}
Let $ \chi_1, \chi_2 $ be unitary characters of $ E^* $ such that $ \chi_1 \notin X_{N_{E/F}(E^*)} $ or 
$ \chi_2 \notin X_{N_{E/F}(E^*)}. $

\smallskip

1. Let $ 0 < \alpha_2 \leq \alpha_1. $ Let $ \alpha_1 \neq \alpha_2 $ or $ \exists x \in E^* $ s. t. $ \chi_1(x) \neq
\chi_2^{-1}(\overline{x}). \; \Lg(\mid \; \mid^{\alpha_1}
\chi_1 ; \mid \; \mid^{\alpha_2} \chi_2 ; \lambda') $ is non-unitary.

2. Let $ 0 < \alpha_1 = \alpha_2 $ and $ \chi_1(x) = \chi_2^{-1}(\overline{x}) \; \forall x \in E^*. \;
\Lg(\mid \; \mid^{\alpha_1}
\chi_1 ; \mid \; \mid^{\alpha_2} \chi_2 ; \lambda') $ is unitary for $ 0 < \alpha_1 = \alpha_2 \leq 1/2. $

$  \Lg(\mid \; \mid^{\alpha_1}
\chi_1 ; \mid \; \mid^{\alpha_2} \chi_2 ; \lambda') $ is non-unitary for $ \alpha > 1/2. $

\end{Theorem}

\begin{Proof}
1. Let $ \alpha_1 \neq \alpha_2 $ or $ \exists x \in E^* $ such that $ \chi_1(x) \neq \chi_2^{-1}(\overline{x}). $
The representations $ \mid \; \mid^{\alpha_1} \chi_1 \times \mid \; \mid^{\alpha_2} \chi_2 \rtimes \lambda' $
 are not hermitian,
by
 \cite{Ca} 3.1.2, $ \Lg(\mid \; \mid^{\alpha_1} \chi_1 ; \mid \; \mid^{\alpha_2} \chi_2 ; \lambda') $ is not hermitian, hence not
unitary.

\medskip

2. Let $ \alpha_1 = \alpha_2 $ and $ \chi_1(x) = \chi_2^{-1}(\overline{x}) \; \forall x \in E^*. $ Representations $
 \mid \; \mid^{\alpha_1} \chi_1 \times \mid \; \mid^{\alpha_2} \chi_2 \rtimes \lambda'$ are hermitian.
Let $ \alpha_1 = \alpha_2 < 1/2 $ and $ \chi_1(x) = \chi_2^{-1}(\overline{x}) \; \forall x \in E^*. $ Representations $ 
\mid \; \mid^{\alpha_1} \chi_1 \times \mid \; \mid^{\alpha_2} \chi_2 \rtimes \lambda'$ are
irreducible by \cite[Th. 4.3]{CS} and equal to their Langlands quotients $ \Lg(\mid \; \mid^{\alpha_1} \chi_1 ; \mid \; \mid^{\alpha_2}
\chi_2 ; \lambda'). \; \chi_1 \times \chi_2 \rtimes \lambda' $ is irreducible by \cite[Th. 4.1]{CS} and unitary.
For $ \alpha_1 = \alpha_2 < 1/2, $ representations $ \Lg(\mid \; \mid^{\alpha_1} \chi_1 ; \mid \; \mid^{\alpha_2}
\chi_2 ; \lambda') $ form a continuous 1-parameter family of irreducible hermitian representations that we realize on the same
vector space V (for a detailed version of this argument in a similar case see the proof of Theorem \ref{Lg11}). By Lemma
\cite[Lemma 3.6]{CS} $ \Lg(\mid \; \mid^{\alpha_1} \chi_1 ; \mid \; \mid^{\alpha_2}
\chi_2 ; \lambda') $ is unitary for $ \alpha_1 = \alpha_2 < 1/2. $ For $ \alpha_1 = \alpha_2 = 1/2 $ and $ \chi_1(x) =
\chi_2^{-1}(\overline{x}) \; \forall x \in E^*, \;
\mid \; \mid^{1/2} \chi_1 \times \mid \; \mid^{1/2} \chi_2 \rtimes \lambda'$ is reducible by \cite[Th. 4.3]{CS}. By
\cite{MR0324429} $ \Lg( \mid \; \mid^{1/2} \chi_1 ; \mid \; \mid^{1/2} \chi_2 ; \lambda')$ is unitary.

Let $ \alpha_1 = \alpha_2 > 1/2 $ and $ \chi_1(x) = \chi_2^{-1}(\overline{x}) \; \forall x \in E^*. $ Representations
$ \mid \; \mid^{\alpha_1} \chi_1 \times \mid \; \mid^{\alpha_2} \chi_2 \rtimes \lambda' $ are irreducible by 
\cite[Th. 4.3]{CS} and equal to their Langlands-quotients $ \Lg(\mid \; \mid^{\alpha_1} \chi_1 ; \mid \; \mid^{\alpha_2} \chi_2 ; \lambda'). $
By \cite[Lemmas 3.6 and 3.8]{CS} $ \Lg(\mid \; \mid^{\alpha_1} \chi_1 ; \mid \; \mid^{\alpha_2} 
\chi_2 ; \lambda') $ is non-unitary for $ \alpha_1 = \alpha_2 > 1/2. $
\end{Proof}

\bigskip

\begin{center} $ \Lg(\mid \; \mid^{\alpha} \chi_1 ;  \chi_2 \rtimes \lambda'), \; \alpha > 0, \; 
\chi_1 \notin X_{N_{E/F}(E^*)} $ \textbf{or} $ \chi_2 \notin X_{N_{E/F}(E^*)}$ \end{center}

\medskip

Let $ \chi_{1_{F^*}} \in X_{1_{F^*}}. $ Recall that $ \chi_{1_{F^*}} \rtimes \lambda' = \sigma_{1, \chi_{1_{F^*}}} \oplus
 \sigma_{2, \chi_{1_{F^*}}}, $ where $ \sigma_{1, \chi_{1_{F^*}}} $
and $ \sigma_{2, \chi_{1_{F^*}}} $ are tempered \cite{Ky}.

\medskip

Recall that $ X_{N_{E/F}(E^*)} = 1 \cup X_{\omega_{E/F}} \cup X_{1_{F^*}}. $

\bigskip

\begin{Theorem}

1. Let $ \alpha > 0. $ Let $ \chi_1 \notin X_{N_{E/F}(E^*)} $ and $ \chi_2 \notin X_{1_{F^*}}. \; 
\Lg(\mid \; \mid^{\alpha} \chi_1 ; \chi_2 \rtimes 
\lambda') $ is non-unitary.

2. Let $ \alpha > 0. $ Let $ \chi_1 \notin X_{N_{E/F}(E^*)}. \; 
\Lg(\mid \; \mid^{\alpha} \chi_1 ; \sigma_{1, \chi_{1_{F^*}}}) $ and 
$\Lg(\mid \; \mid^{\alpha} \chi_1 ; \sigma_{2, \chi_{1_{F^*}}}) $ are non-unitary.

3. Let $ \alpha > 0. $ Let $ \chi_1 \in X_{N_{E/F}(E^*)} $ and $ \chi_2 \notin X_{N_{E/F}(E^*)}. $
 
3.1 Let $ \chi_1 = 1. $ Let $ 0 < \alpha \leq 1. \;  \Lg( \mid \; \mid^{\alpha}
1 ; \chi_2 \rtimes \lambda') $ is unitary.

Let $ \alpha > 1. \; \Lg( \mid \; \mid^{\alpha}
1 ; \chi_2 \rtimes \lambda') $ is non-unitary.

3.2 Let $ \chi_1 \in X_{\omega_{E/F}}. $ Let $  0 < \alpha \leq 1/2. \;  \Lg( \mid \; \mid^{\alpha}
\chi_{\omega_{E/F}} ; \chi_2 \rtimes \lambda') $ is unitary.

Let $ \alpha > 1/2. \; \Lg( \mid \; \mid^{\alpha}
\chi_{\omega_{E/F}} ; \chi_2 \rtimes \lambda') $ is non-unitary.

3.3 Let $ \chi_1 \in X_{1_{F^*}}. $ Let $ \alpha > 0. \; \Lg( \mid \; \mid^{\alpha}
\chi_{1_{F^*}} ; \chi_2 \rtimes \lambda') $ is non-unitary.
\end{Theorem}

\begin{Proof}
1. For $ \alpha > 0, $ representations $ \mid \; \mid^{\alpha} \chi_1 \times \chi_2 \rtimes 
\lambda' $ are not hermitian. By \cite{Ca} 3.1.2, $ \Lg(\mid \; \mid^{\alpha} \chi_1 ; \chi_2 \rtimes 
\lambda') $ is not hermitian, hence not unitary.

\medskip

2. For $ \alpha > 0, \; \mid \; \mid^{\alpha} \chi_1 \rtimes \sigma_{1, \chi_{1_{F^*}}} $ and $\mid \; \mid^{\alpha} \chi_1 \rtimes
 \sigma_{2, \chi_{1_{F^*}}} $ are not hermitian. By \cite{Ca} $ \Lg(\mid \; \mid^{\alpha} \chi_1 ; \sigma_{1, \chi_{1_{F^*}}}) $ and 
$\Lg(\mid \; \mid^{\alpha} \chi_1 ; \sigma_{2, \chi_{1_{F^*}}}) $ are not hermitian and hence non-unitary.

\medskip

3.1 $ 1 \times \chi_2 \rtimes \lambda' $ is irreducible by \cite[Th. 4.1]{CS} and unitary. For $ 0 < \alpha < 1, $
 representations
$ \mid \; \mid^{\alpha} 1 \times \chi_2 \rtimes \lambda' $ are irreducible by \cite[Th. 4.5]{CS} and equal to their
Langlands-quotient $ \Lg(\mid \; \mid^{\alpha} 1 ; \chi_2 \rtimes \lambda'). $ By \cite[Lemma 3.6]{CS} these 
Langlands-quotients
are unitary. For $ \alpha = 1, \; \mid \; \mid 1 \times \chi_2 \rtimes \lambda' $ reduces for the first time, by 
\cite[Th. 4.5]{CS}. By 
\cite{MR0324429} $ \Lg(\mid \; \mid 1 ; \chi_2 \rtimes \lambda') $ is unitary. For $ \alpha > 1, $ representations
$  \mid \; \mid^{\alpha} 1 \times \chi_2 \rtimes \lambda' $ are irreducible by \cite[Th. 4.5]{CS} and equal to 
their
Langlands-quotient $  \Lg(\mid \; \mid^{\alpha} 1 ; \chi_2 \rtimes \lambda') $ . By \cite[Lemmas 3.6 and 3.8]{CS}
these Langlands-quotients are non-unitary.

\medskip

3.2 $ \chi_{\omega_{E/F}} \times \chi_2 \rtimes \lambda' $ is irreducible by \cite[Th. 4.1]{CS} and unitary. 
For $ 0 < \alpha < 1/2,
$ representations
$ \mid \; \mid^{\alpha} \chi_{\omega_{E/F}} \times \chi_2 \rtimes \lambda' $ are irreducible by \cite[Th. 4.5]{CS}
 and equal 
to their
Langlands-quotient $ \Lg(\mid \; \mid^{\alpha} \chi_{\omega_{E/F}} ; \chi_2 \rtimes \lambda'). $ By \cite[Lemma 3.6]{CS}
these 
Langlands-quotients
are unitary. For $ \alpha = 1/2, \; \mid \; \mid^{1/2} \chi_{\omega_{E/F}} \times \chi_2 \rtimes \lambda' $ reduces for the first
time, by \cite[Th. 4.5]{CS}. By 
\cite{MR0324429} $ \Lg(\mid \; \mid^{1/2} \chi_{\omega_{E/F}} ; \chi_2 \rtimes \lambda') $ is unitary. For 
$ \alpha > 1/2, $
 representations
$  \mid \; \mid^{\alpha} \chi_{\omega_{E/F}} \times \chi_2 \rtimes \lambda' $ are irreducible by \cite[Th. 4.5]{CS} and equal
to their
Langlands-quotient $ \Lg(\mid \; \mid^{\alpha} \chi_{\omega_{E/F}} ; \chi_2 \rtimes \lambda') $. By 
\cite[Lemmas 3.6 and 3.8]{CS} these Langlands-quotients are non-unitary.

\medskip

3.3 For $ \alpha > 0, $ representations $ \mid \; \mid^{\alpha} \chi_{1_{F^*}} \times \chi_2 \rtimes \lambda' $ are 
irreducible by \cite[Th. 4.5]{CS} and equal 
to their Langlands-quotient $ \Lg(\mid \; \mid^{\alpha} \chi_{1_{F^*}} ; \chi_2 \rtimes \lambda'). $ By \cite[Lemma 3.6]{CS}
and \cite[Lemma 3.8]{CS} these 
Langlands-quotients are non-unitary.
\end{Proof}

We now take $ \chi_1, \chi_2 \in X_{N_{E/F}(E^*)}. $

\bigskip

\begin{center} $ \Lg( \mid \; \mid^{\alpha_1} 1 ; \mid \; \mid^{\alpha_2} 1 ; \lambda'), \; 0 < \alpha_2 \leq
\alpha_1, \; \Lg( \mid \; \mid^{\alpha} 1 ; 1 \rtimes \lambda'), \; \alpha > 0 $ \end{center}

\medskip

\begin{Theorem}
\label{Lg11}

1. Let $ 0 < \alpha_2 \leq \alpha_1 \leq 1 $ and $ \alpha_1 + \alpha_2 \leq 1.  \; \Lg( \mid \; \mid^{\alpha_1} 1 ; \mid \; \mid^{\alpha_2} 1; \lambda') $ is unitary.

2. Let $ 0 < \alpha_2 \leq \alpha_1, \; \alpha_1 + \alpha_2 > 1, \; (\alpha_1, \alpha_2) \neq (2,1). $ If $ \alpha_1 = 1, \; $
then let $ \alpha_2 \notin (0,1]. $ Then $ 
\Lg( \mid \; \mid^{\alpha_1} 1 ; \mid \; \mid^{\alpha_2} 1; \lambda') $ is non-unitary.

2.1 $ \Lg( \mid \; \mid^2 1 ; \mid \; \mid 1; \lambda') = 1_{U(5)} $ is unitary.

\medskip

3. For $ 0 < \alpha \leq 1, \;  \Lg( \mid \; \mid^{\alpha} 1; 1 \rtimes \lambda') $ is unitary.

4. For $ \alpha > 1, \Lg(  \mid \; \mid^{\alpha} 1; 1 \rtimes \lambda') $ is non-unitary.

\end{Theorem}

\begin{Proof}
1. $ 1 \times 1 \rtimes \lambda' $ is irreducible by \cite[Th. 4.1]{CS} and unitary. Let $ 0 < \alpha_2 \leq \alpha_1. $ For $ 
\alpha_1 + \alpha_2 < 1, \; $ representations $ \mid \; \mid^{\alpha_1} 1 \times \mid \; \mid^{\alpha_2} 1 \rtimes \lambda' $
are irreducible by \cite[Th. 4.3]{CS}, hence equal to $ \Lg( \mid \; \mid^{\alpha_1} 1 ; \mid \; \mid^{\alpha_2} 1;
\lambda'). $

We will now construct a continuous one-parameter family of hermitian representations.
Let $ 0 < \alpha_2 \leq \alpha_1 $ such that $ \alpha_1 + \alpha_2 \leq 1. $ Let $ \pi_{ \alpha_1, \alpha_2} $
denote the two-parameter
family of hermitian representations $ \mid \; \mid^{\alpha_1} 1 \times \mid \; \mid^{ \alpha_2} 1 \rtimes \lambda'. $ Let
$ V_{\alpha_1, \alpha_2} $ be the vector space of $ \pi_{\alpha_1, \alpha_2}. $

\medskip
Recall the Levi decomposition $ P_0 = M_0N_0, $ where $ P_0 $ is the minimal parabolic subgroup of $ U(5), \; M_0 =
\{ \left(
\begin{smallmatrix}
 x&&&&0\\
&y&&&&\\
&&k&&\\
&&&\overline{y}^{-1}\\
0&&&&\overline{x}^{-1}
\end{smallmatrix}
\right), \; x,y \in E^*, \; k \in E^1 \} $ is the minimal Levi subgroup and $ N_0 $ the unipotent radical of $ P_0. $

Let $ (\pi, V) $ be the extension to $ P_0 $ of the representation $ 1 \otimes 1 \otimes \lambda'$ 
of $ M_0. $
Let $ \delta_{P_0} $ denote the modulus character of $ P_0. $

$ V_{0,0} = \{ f: \G \rightarrow V : f  \text{ \; is smooth and} \; f(pg) = \delta_{P_0}(p) \pi(p) f(g) \forall g \in \G \}. $

$ V_{\alpha_1, \alpha_2} = \{ h: \G \rightarrow V: h \text{ \; is smooth and} \; h(pg) = \delta_{P_0}(p) \mid x \mid^{\alpha_1} \mid y 
\mid^{\alpha_2}
 \pi(p) h(g) \forall g \in \G \}, $

\medskip
where $ p \in P_0, \; p =
\left(
\begin{smallmatrix}
 x&&&-&*\\
&y&&& \mid \\
&&k&&\\
\mid&&&\overline{y}^{-1}&\\
0&-&&&\overline{x}^{-1}
\end{smallmatrix}
\right), \; x,y \in E^*, k \in E^1, * \in E. $
Let $ \mathcal{O} $ denote the ring of integers of $ E. $
$ \mid \; \mid: E^* \rightarrow F^* $
is non-ramified, hence $ (x,y) \mapsto \mid x \mid^{\alpha_1} \mid y \mid^{\alpha_2}, \; x,y \in E^*, $ is trivial on 
$ E^1 \times E^1 \cong
\mathcal{O}^* \times \mathcal{O}^*. $ Let $ K := U(\mathcal{O}), $ it is a maximal compact subgroup of $ \G. $ We have
$ \G = KP_0. $
Let $ f \in V_{0,0}. $ There exists a unique extension of $ f_{\mid K} : \G \rightarrow V $ to a function 
$ h \in V_{\alpha_1, \alpha_2}, $ so
$ f_{\mid K} = h_{\mid K}. $ This induces an isomorphism $ T_{\alpha_1, \alpha_2}: V_{0,0} \overset{\sim}{\rightarrow} 
V_{\alpha_1, \alpha_2}. $ Via the composition with $ T_{\alpha_1, \alpha_2} $ we consider all representations
 $ \pi_{\alpha_1, \alpha_2} $ in $ V_{0,0}. $

\medskip

Let $ w \in W $ be the longest element of the Weyl group. Let

$ A(w, \lambda): \mid \; \mid^{\alpha_1} 1 \times \mid \; \mid^{\alpha_2} 1 \rtimes \lambda' \rightarrow
\mid \; \mid^{- \alpha_1} 1 \times \mid \; \mid^{- \alpha_2} 1 \rtimes \lambda' $ be the standard long intertwining
operator.

\medskip

On $ V_{0,0} $ we define a set of non-degenerate hermitian
forms $ \langle \; , \; \rangle_{ \alpha_1, \alpha_2} $ by

$$ \langle f, h \rangle_{\alpha_1, \alpha_2} = \underset{U(\mathcal{O})}{\int} A(w, \lambda) 
f(k) \overline{h(k)} dk, \; f, h \in V_{0, 0}, $$

such that $  \langle \; , \; \rangle_{ \alpha_1, \alpha_2} $ is invariant by $ T_{\alpha_1, \alpha_2}^{-1}
 \pi_{\alpha_1, \alpha_2} T_{\alpha_1, \alpha_2}. $

\medskip
 
Fix $ \alpha_1 $ and $ \alpha_2 $ such that $ \alpha_1 + \alpha_2 = 1. $ Let $ \pi_t =
 \pi_{t \alpha_1, t \alpha_2}, \; t \in [0,1] $ denote a continuous one-parameter
family of hermitian representations.
Let $ V_t $ be the vector space of $ \pi_t. $
Via the isomorphism $ T_t: V_0 \overset{\sim}{\rightarrow} V_t, $ we consider all representations $ \pi_t $ in $ V_0, $ as before.

\medskip

Choose a polynomial $ p(t) $ with real
coefficents, such that
$ A(t) = p(t) A(w, \lambda) $ is holomorphic and non-zero for $ t \in [0,1]. $
So for the one-parameter family of representations $ \pi_t $ one obtains, on the same space $ V_0 $,
 a set of non-degenerate
hermitian forms
$ \langle \; , \; \rangle_t $ given by

$$ \langle f, h \rangle_t = \underset{U(\mathcal{O})}{\int} A(t) f(k) \overline{h(k)}
dk, \; f,h \in V_0, $$

such that $ \langle \; , \; \rangle_t $ is invariant under $ T^{-1}_t \pi_t T_t. $

$ \langle \; , \; \rangle_0 $ is positive definite, hence by \cite[Lemma 3.6]{CS} $ \langle \; , \; \rangle_t $ is positive 
definite
until $ \mid \; \mid^{t \alpha_1} 1 \times \mid \; \mid^{t \alpha_2} \rtimes \lambda' $ reduces for the first time,
for $ t = 1. $ By \cite{MR0324429}, for $ t=1,$
 the irreducible subquotients of $ 
\mid \; \mid^{\alpha_1} 1 \times \mid \; \mid^{ \alpha_2} 1
\rtimes \lambda' $ are unitary.
Hence for $ 0 < \alpha_2 \leq \alpha_1 \leq 1, \; \alpha_1 + \alpha_2 \leq
 1, $
the Langlands quotients $ \Lg( \mid \; \mid^{\alpha_1} 1 ; \mid \; \mid^{\alpha_2} 1; \lambda') $ are unitary.

\medskip

2. + 2.1
$ \mid \; \mid^2 1 \times \mid \; \mid 1 
\rtimes \lambda' $ is reducible. $ \lambda'(\det) \St_{U(3)} $ is
the unique irreducible square-integrable subquotient of $ \mid \; \mid 1 \rtimes \lambda' $ \cite{Ky}. By 
\cite[p.915]{MR656064} the subquotient $ \mid \; \mid^{3/2} 1_{\GL_2} \rtimes \lambda' $ is
reducible and has the subquotients $ 1_{U(5)} = \Lg( \mid \; \mid^2 1 ; \mid \; \mid 1 ; \lambda') $ and
 $ \Lg(\mid \; \mid^2 1 ; \lambda'(\det) \St_{U(3)}). \; 1_{U(5)} = \Lg( \mid \; \mid^2 1; \mid \; \mid 1 ; \lambda') $
 is unitary, $ \Lg(\mid \; \mid^2 1 ;
\lambda'(\det) \St_{U(3)}) $ is non-unitary (\cite{MR656064}).
 For $ 1/2 < \alpha < 3/2, $ representations $ \mid \; \mid^{\alpha} 1_{\GL_2} \rtimes \lambda' $ are irreducible by
\cite[Th. 4.7]{CS}, they form a continuous one-parameter 
family of irreducible hermitian 
representations on the space $ V_{\alpha}. $ Like before we identify the vectorspaces
 $ V_{\alpha} $
for $ 1/2 < \alpha < 3/2. $ For $ \alpha = 3/2, $ the irreducible subquotient $ \Lg(\mid \; \mid^2 1 ; \lambda'(\det) \St_{U(3)})
 $ of $  \mid \; \mid^{3/2} 1_{\GL_2} \rtimes \lambda' $ is not unitary \cite{MR656064}. Hence, by \cite{MR0324429} and by 
 \cite[Lemma 3.6]{CS},
$ \mid \; \mid^{\alpha} 1_{\GL_2} \rtimes \lambda' = \Lg(\mid \; \mid^{\alpha_1} 1; \mid \; \mid^{\alpha_2} 1 ; \lambda')
 $  is non-unitary for $ 1/2 < \alpha < 3/2, $ i.e. for $ 1 < \alpha_1 < 2, \; \alpha_1 - \alpha_2 = 1 .$

\smallskip

By \cite[p.915]{MR656064} the subquotient $ \mid \; \mid^2 1 \rtimes \lambda'(\det) 1_{U(3)} $ of the representation
$ \mid \; \mid^2 1 \times \mid \; \mid 1 \rtimes \lambda'$ is
reducible. It has the subquotients $ 1_{U(5)} = \Lg( \mid \; \mid^2 1 ; \mid \; \mid 1 ; \lambda') $ and
 $ \Lg(\mid \; \mid^{3/2} \St_{\GL_2} ; \lambda'). \; 1_{U(5)} $ is unitary, $ \Lg(\mid \; \mid^{3/2} \St_{\GL_2} ;
\lambda') $ is non-unitary \cite{MR656064}.

Let $ 1 < \alpha < 2. $ Representations $ \mid \; \mid^{\alpha} 1 \rtimes \lambda'(\det) 1_{U(3)} $ are irreducible by
\cite[Th. 4.7]{CS}, they form a continuous one-parameter
 family of irreducible hermitian representations
on the space $ V_{\alpha}. $ Similar as before we identify
$ V_\alpha $ for $ 1 < \alpha < 2. $ The irreducible subquotient $ \Lg(\mid \; \mid^{3/2} \St_{\GL_2} ;
\lambda') $ of  $ \mid \; \mid^2 1 \rtimes \lambda'(\det) 1_{U(3)} $ is non-unitary \cite{MR656064}. Hence by \cite{MR0324429} and
 by \cite[Lemma 3.6]{CS}
representations $ \mid \; \mid^{\alpha} 1 \rtimes \lambda'(\det) 1_{U(3)} = \Lg(\mid \; \mid^{\alpha_1} 1 ; \mid \; 
\mid^1
1; \lambda') $ are non-unitary for $  1 = \alpha_2 < \alpha_1 = \alpha < 2. $

\smallskip

Let $ 1 < \alpha_1 < 2, \; 0 < \alpha_2 < 1, \; \alpha_1 - \alpha_2 < 1. \; \mid \; \mid^{\alpha_1} 1 \times \mid \; 
\mid^{\alpha_2} 1 \rtimes \lambda' $ is irreducible by \cite[Th. 4.3]{CS} and equal to its own Langlands quotient 
$ \Lg(\mid \; \mid^{\alpha_1} 1 ; \mid \; \mid^{\alpha_2} 1 ; \lambda') $.
Fix $ 0 < \alpha_2 < 1 $ and  let $ 1 < \alpha_1 \leq \alpha_2 + 1. $ Let $ \pi_{\alpha_1} $ denote the continuous 
one-parameter
family of hermitian representations $ \mid \; \mid^{\alpha_1} 1 \times \mid \; \mid^{\alpha_2} 1 \rtimes \lambda'$ on
the same vector space $ V. $
For $ \alpha_1 = \alpha_2 + 1 $ irreducible subquotients of the representations $ \mid \; \mid^{\alpha_1} 1 \times
\mid \; \mid^{\alpha_2} 1 \rtimes \lambda' $ are non-unitary, as seen
in the previous paragraph. 
By \cite{MR0324429} and by \cite[Lemma 3.6]{CS} the Langlands quotients 
$ \Lg(\mid \; \mid^{\alpha_1} 1 ; \mid \; \mid^{\alpha_2} 1 ; \lambda') $ are non-unitary. (II in Figure ~\ref{figure1} on
page ~\pageref{figure1})

\medskip

Let $ \alpha_1 > 2, \;  \alpha_1 - \alpha_2 = 1. \;
 \Lg(\mid \; \mid^{\alpha_1} 1 ; \mid \; \mid^{\alpha_2} 1 ; \lambda') $ is non-unitary by Lemmas \ref{cont} and \ref{comp}.
For $ \alpha_1 > 2, \; \alpha_2 = 1, \; \Lg(\mid \; \mid^{\alpha_1} 1 ; \mid \; \mid^{\alpha_2} 1 ; \lambda') $ is
 non-unitary
by the same argument.

The same holds for $ \Lg(\mid \; \mid^{\alpha_1} 1 ; \mid \; \mid^{\alpha_2} 1 ; \lambda'), $
where $ \alpha_1 >1, 0 < \alpha_2 < 1, \; \alpha_1 - \alpha_2 > 1, $ for
$ 1 < \alpha_2 \leq \alpha_1, \; \alpha_1 - \alpha_2 < 1 $ and for $ \alpha_1 > 2,\; \alpha_1 - \alpha_2 > 1 $ 
(III, IV, V in Figure ~\ref{figure4.1} on page ~\pageref{figure4.1}).

\medskip

Let $ \alpha_1 = 1, \; \alpha_2 \in (0,1]. $ We have no proof that $ \Lg(\mid \; \mid 1 ; \mid \; \mid^{\alpha_2} 1 ; \lambda') $
is non-unitary.

\medskip

3. $ 1 \times 1 \times \lambda' $ is irreducible by Theorem \cite[Th. 4.1]{CS} and unitary. For $ 0 < \alpha < 1, \; 
\mid \; \mid^{\alpha} 1 \times 1 \rtimes
\lambda' $ is irreducible by \cite[Th. 4.5]{CS} and equal to its own Langlands quotient $ \Lg(\mid \; \mid^{\alpha} 1 ; 1 \rtimes
\lambda'). $ By Lemma \cite[Lemma 3.6]{CS} these representations are unitary. For $ \alpha = 1, \; \mid \; \mid 1 \times 1 \rtimes
\lambda' $ reduces for the first time, hence $ \Lg(\mid \; \mid 1 ; 1 \rtimes \lambda') $ is unitary \cite{MR0324429}.

\medskip

4. For $ \alpha > 1, \; \mid \; \mid^{\alpha} 1 \times 1 \rtimes \lambda' = \Lg(\mid \; \mid^{\alpha} 1 ; 1 \rtimes
\lambda') $ is irreducible by \cite[Th. 4.5]{CS} and by \cite[Lemmas 3.6 and 3.8]{CS} non-unitary.
\end{Proof}

\begin{figure}

\begin{center}
\begin{tikzpicture}
 \draw[->] (0,0) -- (3,0);
\draw (3,0) node[below right] {$\alpha_1$};
\draw[->] (0,0) -- (0,3);
\draw (0,3) node[above left] {$ \alpha_2 $};

\draw [fill=blue!50] (0,0) -- (1,0) -- (0,1) --cycle;
\draw [blue] (1,0) -- (0,1);

\draw [red] (1,1) -- (3,1);
\draw [red] (0,1) -- (2,3);
\draw [red] (1,1) -- (1,3);
\draw [red] (1,0) -- (3,2);
\draw 	    (1,0) -- (1,1);
\draw	    (0,1) -- (1,1);

\draw (1,0) node[below] {$ 1 $};
\draw (0,1) node[left] {$ 1 $};
\draw (2,0) node[below] {$ 2 $};
\draw (0,2) node[left] {$ 2 $};
\draw (1,0) node[blue] {$ \bullet $};
\draw (1,1) node {$ \bullet $};
\draw (2,1) node {$ \bullet $};
\draw (1,0) node[above] {$ \cdot \cdot \cdot \cdot $};
\draw (1.9,1) node[blue, above] { $ \cdot \cdot $};
\draw (2.1, 1) node[above] { $ \cdot \cdot $};
\draw (0.55, 0.4) node[above] { $ \cdot \cdot$};
\draw (1.55, 0.4) node[above] { $ \cdot \cdot$};
\draw (0.55, 0.9) node[above] { $ \cdot \cdot $};

\draw (0.75, 0.6) node {I};
\draw (1.25, 0.6) node {II};
\draw (2.2, 0.5) node {III};
\draw (1.75, 1.6) node {IV};
\draw (2.6, 1.25) node {V};

\draw (1.5, -0.5) node[below] { Figure 1 };
\draw (3, 2.2) node[red, right] { $ - $};
\draw (3.4, 2.2) node[right] {reducible;  irreducible subquotients non-unitary};
\draw (3, 1.8) node[blue, right] { $ - $};
\draw (3.4, 1.8) node[right] {reducible; irreducible subquotients unitary};
\draw (3.1,1.4) node[right] { $ \cdot \; \; $ irreducible non-unitary subquotient };
\draw (3.1, 1) node[blue, right] { $ \cdot $};
\draw (3.4, 1) node[right] { irreducible unitary subquotient}; 
\end{tikzpicture}
\end{center}

\caption[$ \Lg( \mid \; \mid^{\alpha_1} 1 ; \mid \; \mid^{\alpha_2} 1 ; \lambda'), \; 0 < \alpha_2 \leq
\alpha_1, \; \Lg( \mid \; \mid^{\alpha} 1 ; 1 \rtimes \lambda'), \; \alpha > 0 $]{Let $ \alpha_1, \alpha_2 \geq 0. $
Figure 1 shows lines and points of reducibility of the representation 
$ \mid \; \mid^{\alpha_1} 1 \times \mid \; \mid^{\alpha_2} 1 \rtimes \lambda' $ and the unitary dual. Let $ 0 < \alpha_2
\leq \alpha_1. \; \Lg( \mid \; \mid^{\alpha_1} 1 ; \mid \; \mid^{\alpha_2} 1 ; \lambda') $ is unitary for $ 0 < \alpha_2
\leq \alpha_1 < 1, \alpha_1 + \alpha_2 \leq 1 $ and for $ \alpha_1 = 2, \; \alpha_2 = 1. $ Except for 
$ \alpha_1 = 1, \; 0 < \alpha_2 
\leq 1, $ it is non-unitary.
$ \Lg( \mid \; \mid^{\alpha} 1 ;  1 \rtimes \lambda') $ is unitary for $ 0 < \alpha \leq 1 $ and non-unitary for
$ \alpha > 1.$}

\label{figure1}

\end{figure}

\bigskip

\begin{center} $ \Lg( \mid \; \mid^{\alpha_1} \chi_{\omega_{E/F}} ; \mid \; \mid^{\alpha_2} \chi_{\omega_{E/F}} ; \lambda'),
\; 0 < \alpha_2 \leq \alpha_1,  \; \Lg( \mid \; \mid^{\alpha}
\chi_{\omega_{E/F}} ; \chi_{\omega_{E/F}} \rtimes \lambda'), \; \alpha > 0, \; \; \chi_{\omega_{E/F}} \in
 X_{\omega_{E/F}} $ \end{center}

\bigskip

Let $ \chi_{\omega_{E/F}} \in X_{\omega_{E/F}}. $

\begin{Theorem}
\label{lgalpha12chiomega}
1. Let $ 0 < \alpha_2 \leq \alpha_1 \leq 1/2. \;
\Lg( \mid \; \mid^{\alpha_1} \chi_{\omega_{E/F}} ; \mid \; \mid^{\alpha_2} \chi_{\omega_{E/F}} ; \lambda') $
is unitary.

2.  Let $ \alpha_1 > 1/2, \; \alpha_2 \leq \alpha_1, \; (\alpha_1, \alpha_2) \neq (3/2, 1/2). $ If $ 0 < \alpha_2
< 1/2, $ then let $ \alpha_1 \notin (1/2, 1 - \alpha_2].
\; 
\Lg( \mid \; \mid^{\alpha_1} \chi_{\omega_{E/F}} ; \mid \; \mid^{\alpha_2} \chi_{\omega_{E/F}} ; \lambda') $ is
non-unitary.

\medskip

3. Let $ 0 < \alpha \leq 1/2. \; 
\Lg( \mid \; \mid^{\alpha} \chi_{\omega_{E/F}} ; \chi_{\omega_{E/F}} \rtimes \lambda') $ is unitary.

4. Let $ \alpha > 1. \; 
\Lg( \mid \; \mid^{\alpha} \chi_{\omega_{E/F}} ;  \chi_{\omega_{E/F}} \rtimes \lambda') $ is non-unitary.
\end{Theorem}

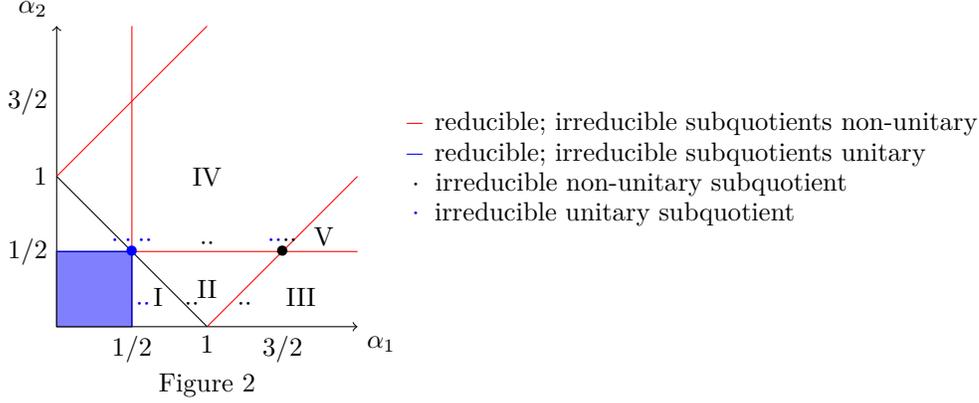
\begin{figure}[!h]

\begin{center}
\begin{tikzpicture}
 \draw[->] (0,0) -- (4,0);
\draw (4,0) node[below right] {$\alpha_1$};
\draw[->] (0,0) -- (0,4);
\draw (0,4) node[above left] {$ \alpha_2 $};

\draw [fill=blue!50] (0,0) -- (1,0) -- (1,1) -- (0,1) -- cycle;

\draw [red] (1 ,1) -- (4,1);
\draw [red] (0,2) -- (2,4);
\draw [red] (1,1) -- (1,4);
\draw [red] (2,0) -- (4,2);
\draw (2,0) -- (0,2);
\draw [blue] (1, 0) -- (1, 1);
\draw [blue] (0, 1) -- (1, 1);

\draw (1,0) node[below] {$ 1/2 $};
\draw (0, 1) node[left] {$ 1/2 $};
\draw (2,0) node[below] {$ 1 $};
\draw (0,2) node[left] {$ 1 $};
\draw (3,0) node[below] { $ 3/2 $};
\draw (0, 3) node[left] { $ 3/2 $};
\draw (1, 1) node[blue] {$ \bullet $};
\draw (3, 1) node {$ \bullet $};
\draw (2.9,0.95) node[blue, above] {$ \cdot \cdot$};
\draw (3.1,0.95) node[above]{$ \cdot \cdot$};
\draw (1,0.95) node[blue, above] { $ \cdot \cdot \cdot \cdot $};
\draw (1.8, 0.1) node[above] { $ \cdot \cdot $};
\draw (1.15, 0.1) node[blue, above] { $ \cdot \cdot $};
\draw (2.5, 0.1) node[above] { $ \cdot \cdot $};
\draw (2, 0.9) node[above]  {$ \cdot \cdot $};

\draw (2, -0.5) node[below] { Figure 2 };
\draw (4.5, 2.7) node[red, right] { $ - $};
\draw (4.9, 2.7) node[right] {reducible; irreducible subquotients non-unitary};
\draw (4.5, 2.3) node[blue, right] { $ - $};
\draw (4.9, 2.3) node[right] {reducible; irreducible subquotients unitary};
\draw (4.6,1.9) node[right] { $ \cdot \; $ irreducible non-unitary subquotient };
\draw (4.6, 1.5) node[blue, right] { $ \cdot $};
\draw (4.9, 1.5) node[right] { irreducible unitary subquotient}; 

\draw (1.35, 0.4) node {I};
\draw (2, 0.5) node {II};
\draw( 3.25, 0.4) node {III};
\draw (2,2) node {IV};
\draw (3.55, 1.2) node {V};

\end{tikzpicture}
\end{center}

\caption[$ \Lg( \mid \; \mid^{\alpha_1} \chi_{\omega_{E/F}} ; \mid \; \mid^{\alpha_2} \chi_{\omega_{E/F}} ; \lambda'),
\; 0 < \alpha_2 \leq \alpha_1,  \; \Lg( \mid \; \mid^{\alpha}
\chi_{\omega_{E/F}} ; \chi_{\omega_{E/F}} \rtimes \lambda'), \; \alpha > 0, \; \; \chi_{\omega_{E/F}} \in
 X_{\omega_{E/F}} $]{ Let $ \alpha_1, \alpha_2 \geq 0, $ let 
$ \chi_{\omega_{E/F}} \in X_{\omega_{E/F}}. $ Figure 2 shows lines and points of reducibility of the representation 
$ \mid \; \mid^{\alpha_1} \chi_{\omega_{E/F}} \times \mid \; \mid^{\alpha_2} \chi_{\omega_{E/F}} \rtimes \lambda' $ and 
the unitary dual.
Let $ 0 < \alpha_2
\leq \alpha_1. \; \Lg( \mid \; \mid^{\alpha_1} \chi_{\omega_{E/F}} ; \mid \; \mid^{\alpha_2} \chi_{\omega_{E/F}} ; 
\lambda') $ is unitary for $ 0 < \alpha_2
\leq \alpha_1 \leq 1/2. $ Except for 
$ 1/2 < \alpha_1 < 1, \; 0 < \alpha_2 \leq 1 - \alpha_1 $ and for $ \alpha_1 = 3/2, \; \alpha_2 = 1/2, $  it is non-unitary.
$ \Lg( \mid \; \mid^{\alpha} \chi_{\omega_{E/F}} ;  \chi_{\omega_{E/F}} \rtimes \lambda') $ is unitary for $ 0 < \alpha
 \leq 1/2. $ For $ \alpha > 1 $ it is non-unitary.}

\label{figure2}

\end{figure}

\begin{Proof}
1. $ \chi_{\omega_{E/F}} \times \chi_{\omega_{E/F}} \rtimes \lambda' $ is irreducible by \cite[Th. 4.1]{CS} and unitary.
For $ 0 < \alpha_2 \leq \alpha_1 < 1/2 \; \mid \; \mid^{\alpha_1} \chi_{\omega_{E/F}} \times \mid \; \mid^{\alpha_2}
 \chi_{\omega_{E/F}} \rtimes \lambda' $ is 
irreducible by \cite[Th. 4.3]{CS} and equal to its Langlands quotient $ \Lg( \mid \; \mid^{\alpha_1} 
\chi_{\omega_{E/F}} ; \mid \; \mid^{\alpha_2} \chi_{\omega_{E/F}} ; \lambda'). $
Let $ \alpha_1 = 1/2 $ and fix $ 0 < \alpha_2 \leq 1/2. $ For $ t \in [0,1], $ let $  \pi_{(t 1/2, t \alpha_2)} =: \pi_t $ 
denote the continuous one-parameter family of hermitian
representations $ \mid \; \mid^{t 1/2} \chi_{\omega_{E/F}} \times \mid \; \mid^{t \alpha_2} \chi_{\omega_{E/F}} \rtimes
\lambda'. $ For $ t \in [0,1) $ these representations are equal to their own Langlands quotient
$ \Lg( \mid \; \mid^{t 1/2} \chi_{\omega_{E/F}} ; \mid \; \mid^{t \alpha_2} \chi_{\omega_{E/F}} ; \lambda') $ and by 
\cite[Lemma 3.6]{CS} unitary. For $ t = 1 $ the representations $ \pi_t $ reduce 
for the
first time. By \cite{MR0324429} $ \Lg( \mid \; \mid^{1/2} \chi_{\omega_{E/F}} ;
 \mid \; \mid^{\alpha_2} \chi_{\omega_{E/F}} ; \lambda') $ is unitary.

\medskip

2. Let $ 1 < \alpha_1 < 3/2 $ and let $ \alpha_1 - \alpha_2 =1. $ Then 
$ \mid \; \mid^{\alpha_1} \chi_{\omega_{E/F}} \times \mid \; \mid^{\alpha_2} \chi_{\omega_{E/F}} \rtimes \lambda' $ is
reducible by \cite[Th. 4.3]{CS}. The subquotients 
$ \mid \; \mid^{\frac{\alpha_1 +\alpha_2}{2}} \chi_{\omega_{E/F}} 1_{\GL_2} \rtimes \lambda' $ are irreducible by
\cite[Th. 4.7]{CS}. They form a continuous 1-parameter family of irreducible
hermitian representations, that similar as before,
we realize on the same vector space V. Let $ \alpha_1 = 3/2, \; \alpha_2 = 1/2. \;  \mid \; \mid^{\frac{3/2 + 1/2}{2}}
\chi_{\omega_{E/F}} 1_{\GL_2} \rtimes \lambda' =  \mid \; \mid \chi_{\omega_{E/F}} 1_{\GL_2} \rtimes \lambda' $ reduces 
by \cite[Th. 5.4]{CS}.
Let $ \pi_{1, \chi_{\omega_{E/F}}} $ be the unique square-integrable subquotient of $ \mid \; \mid^{1/2}
 \chi_{\omega_{E/F}} \rtimes \lambda' $ \cite{Ky}.
By Theorem 1.1 and Remark 4.7 in \cite{MR2652536} the 
irreducible subquotient $ \Lg(\mid \; \mid^{3/2} \chi_{\omega_{E/F}} ; \pi_{1, \chi_{\omega_{E/F}}}) $ of $ \mid \; \mid
\chi_{\omega_{E/F}} 1_{\GL_2} \rtimes \lambda' $ is non-unitary. By \cite{MR0324429} and \cite[Lemma 3.6]{CS}
the representations
$ \mid \; \mid^{\frac{\alpha_1 +\alpha_2}{2}} \chi_{\omega_{E/F}} 1_{\GL_2} \rtimes \lambda', \; 1 < \alpha_1 <
3/2, \; \alpha_1 - \alpha_2 =
1, $ that are equal to $ \Lg( \mid \; \mid^{\alpha_1} 
\chi_{\omega_{E/F}} ; \mid \; \mid^{\alpha_2} \chi_{\omega_{E/F}} ; \lambda'), \; 1 < \alpha_1 < 3/2, \;
 \alpha_1 - \alpha_2 =
1, $ are non-unitary. 

\medskip

Let $ 1/2 < \alpha_1 < 3/2, \; \alpha_2 = 1/2. $ The representations 
$ \mid \; \mid^{\alpha_1} \chi_{\omega_{E/F}} \times \mid \; \mid^{1/2} \chi_{\omega_{E/F}} \rtimes \lambda' $ are
reducible by \cite[Th. 4.3]{CS}. Let $ \pi_{2, \chi_{\omega_{E/F}}} $ be the unique irreducible non-tempered 
subquotient of $ \mid \; \mid^{1/2} \chi_{\omega_{E/F}} \rtimes
\lambda' $ \cite{Ky}. For $  1/2 < \alpha_1 < 3/2, $ the representations 
$ \mid \; \mid^{\alpha_1} \chi_{\omega_{E/F}} \rtimes \pi_{2, \chi_{\omega_{E/F}}}  $ are irreducible by \cite[Th. 4.11]{CS}
and
equal to the
Langlands quotient $ \Lg(\mid \; \mid^{\alpha_1} \chi_{\omega_{E/F}} ; 
\mid \; \mid^{1/2} \chi_{\omega_{E/F}} ; \lambda'). $
They form a 1 - parameter family of irreducible
hermitian representations, that we realise on the same vector space $ V. $ For $ \alpha_1 = 3/2, \; \mid \; \mid^{3/2}
 \chi_{\omega_{E/F}} \rtimes \pi_{2, \chi_{\omega_{E/F}}} $ reduces by \cite[Th. 5.4]{CS}, and by Theorem 1.1 and 
Remark 4.7 in \cite{MR2652536} its 
irreducible subquotient $ \Lg(\mid \; \mid \chi_{\omega_{E/F}} \St_{\GL_2} ; \lambda') $ is non-unitary.
By \cite{MR0324429} and \cite[Lemma 3.6]{CS} $ \mid \; \mid^{\alpha_1} \chi_{\omega_{E/F}} \rtimes \pi_2 =
\Lg( \mid \; \mid^{\alpha_1} 
\chi_{\omega_{E/F}} ; \mid \; \mid^{1/2} \chi_{\omega_{E/F}} ; \lambda') $ is non-unitary for $ 1/2 < \alpha_1 <
3/2. $

\medskip

Representations $ \mid \; \mid^{\alpha_1} \chi_{\omega_{E/F}} \times \mid \; \mid^{\alpha_2} \chi_{\omega_{E/F}} \rtimes
\lambda' $  in II,III,IV,V  of Figure ~\ref{figure2} on page ~\pageref{figure2} are irreducible by \cite[Th. 4.3]{CS}
and equal to their own Langlands quotient
$ \Lg(\mid \; \mid^{\alpha_1}
 \chi_{\omega_{E/F}} ; \mid \; \mid^{\alpha_2} \chi_{\omega_{E/F}} ; \lambda') . $ The Langlands quotients $
 \Lg(\mid \; \mid^{\alpha_1}
 \chi_{\omega_{E/F}} ; \mid \; \mid^{\alpha_2} \chi_{\omega_{E/F}} ; \lambda') $ in  II are non-unitary by \cite{MR0324429}
and \cite[Lemma 3.6]{CS}.  $  \Lg(\mid \; \mid^{\alpha_1}
 \chi_{\omega_{E/F}} ; \mid \; \mid^{\alpha_2} \chi_{\omega_{E/F}} ; \lambda') $ in III, IV and V are non-unitary by
\cite[Lemmas 3.6 and 3.8]{CS}.

\medskip

3. $ \chi_{\omega_{E/F}} \times \chi_{\omega_{E/F}} \rtimes \lambda' $ is irreducible by \cite[Th. 4.1]{CS} and unitary.
For $ 0 < \alpha < 1/2 \; \mid \; \mid^{\alpha} \chi_{\omega_{E/F}} \times \chi_{\omega_{E/F}} \rtimes \lambda' $ is 
irreducible by Theorem \cite[Th. 4.5]{CS} and equal to its Langlands quotient $ \Lg( \mid \; \mid^{\alpha} \chi_{\omega_{E/F}} ; \chi_{\omega_{E/F}} \rtimes
\lambda'). $ By \cite[Lemma 3.6]{CS} these Langlands quotients are unitary.

For $ \alpha = 1/2, \; \mid \; \mid^{1/2} \chi_{\omega_{E/F}} \times \chi_{\omega_{E/F}} \rtimes
\lambda' $ reduces for the first time (\cite[Th. 4.5]{CS}). By \cite{MR0324429} $ \Lg(\mid \; \mid^{1/2} \chi_{\omega_{E/F}} ; \chi_{\omega_{E/F}} \rtimes
\lambda') $ is unitary.

\medskip

4. For $ \alpha > 1, \; \mid \; \mid^{\alpha} \chi_{\omega_{E/F}} \times \chi_{\omega_{E/F}} \rtimes
\lambda' $ is irreducible by \cite[Th. 4.5]{CS} and equal to its Langlands quotient
$ \Lg(\mid \; \mid^{\alpha} \chi_{\omega_{E/F}} ; \chi_{\omega_{E/F}} ;
\lambda'). $ By \cite[Lemmas 3.6 and 3.8]{CS} these Langlands quotients are non-unitary.
\end{Proof}

\begin{Remark}
We do not have a proof that the representation $ \Lg( \mid \; \mid^{3/2} \chi_{\omega_{E/F}} \times \mid \; \mid^{1/2} 
\chi_{\omega_{E/F}} \rtimes \lambda') $ is unitary. It is the Aubert dual of a square-integrable representation and should
be unitary, see \cite[Th. 5.4]{CS}. See \cite{MR2460906}, where the proof is given for orthogonal and symplectic
groups.
\end{Remark}

\begin{Remark}

\label{remark1chiomega}

 In the Grothendieck group of the category of admissible representations of finite length one has
$ \mid \; \mid \chi_{\omega_{E/F}} \times \chi_{\omega_{E/F}} \rtimes \lambda' = \mid \; \mid^{1/2} \chi_{\omega_{E/F}}
\St_{\GL_2} \rtimes \lambda' + \mid \; \mid^{1/2} \chi_{\omega_{E/F}} 1_{\GL_2} \rtimes \lambda'. $ If we assume that
$ \mid \; \mid^{1/2} \chi_{\omega_{E/F}}
\St_{\GL_2} \rtimes \lambda'$ and $ \mid \; \mid^{1/2} \chi_{\omega_{E/F}} 1_{\GL_2} \rtimes \lambda' $ are irreducible
(see \cite[Remark 5.5]{CS}), we are able to prove that
$  \Lg(\mid \; \mid^{\alpha_1}
 \chi_{\omega_{E/F}} ; \mid \; \mid^{\alpha_2} \chi_{\omega_{E/F}} ; \lambda') $ is non-unitary for 
$ 1/2 < \alpha_1 < 1, \; \alpha_2 \leq 1 - \alpha_1 $ and that $ \Lg(\mid \; \mid^{\alpha} \chi_{\omega_{E/F}} ; 
\chi_{\omega_{E/F}} ;
\lambda') $ is non-unitary for $ 1/2 < \alpha \leq 1: $

\medskip

Let $ 1/2 < \alpha_1 < 1, \; \alpha_2 = 1 - \alpha_1. $ By \cite[Th. 4.3]{CS} representations
$ \mid \; \mid^{\alpha_1} 
\chi_{\omega_{E/F}} \times \mid \; \mid^{\alpha_2} \chi_{\omega_{E/F}} \rtimes \lambda' $ are reducible. By 
\cite[Th. 4.7]{CS}
the subquotient $ \mid \; \mid^{\frac{\alpha_1 - \alpha_2}{2}} \chi_{\omega_{E/F}} 1_{\GL_2} \rtimes 
\lambda' $ is irreducible. It is equal to $  \Lg(\mid \; \mid^{\alpha_1}
 \chi_{\omega_{E/F}} ; \mid \; \mid^{\alpha_2} \chi_{\omega_{E/F}} ; \lambda'). $ By assumption $ \mid \; \mid^{1/2}
 \chi_{\omega_{E/F}} 1_{\GL_2} \rtimes \lambda' $ is irreducible, it is equal to $ \Lg(\mid \; \mid \chi_{\omega_{E/F}} ;
 \chi_{\omega_{E/F}} \rtimes \lambda'). $ Hence we can extend the argument 2 in the proof of \ref{lgalpha12chiomega}:
$ \Lg(\mid \; \mid^{\alpha_1}
 \chi_{\omega_{E/F}} ; \mid \; \mid^{\alpha_2} \chi_{\omega_{E/F}} ; \lambda') $ is non-unitary for $ 1/2 < \alpha_1 <  1$
 and $ \alpha_2 = 1 - \alpha_1 $ and for $ 1 < \alpha_1 < 3/2, \; \alpha_2 = \alpha_1 - 1, $ and 
 $ \Lg(\mid \; \mid \chi_{\omega_{E/F}} ; \chi_{\omega_{E/F}} \rtimes \lambda') $ is non-unitary.

\medskip

Let $ 1/2 < \alpha_1 < 1, \; \alpha_2 < 1 - \alpha_1. $ By \cite[Th. 4.3]{CS} the representations
$  \mid \; \mid^{\alpha_1}
 \chi_{\omega_{E/F}} \times \mid \; \mid^{\alpha_2} \chi_{\omega_{E/F}} \rtimes \lambda' $ are irreducible, they are
equal to their Langlands-quotient $  \Lg(\mid \; \mid^{\alpha_1}
 \chi_{\omega_{E/F}} ; \mid \; \mid^{\alpha_2} \chi_{\omega_{E/F}} ; \lambda') $ (I in Figure ~\ref{figure2}, page
~\pageref{figure2}). By \cite{MR0324429}
and by \cite[Lemma 3.6]{CS} these Langlands-quotients are non-unitary.

\medskip

Let $ 1/2 < \alpha < 1. $ By \cite[Th. 4.5]{CS} representations $  \mid \; \mid^{\alpha}
 \chi_{\omega_{E/F}} \times \chi_{\omega_{E/F}} \rtimes \lambda' $ are irreducible, they are
equal to their Langlands-quotient $  \Lg(\mid \; \mid^{\alpha}
 \chi_{\omega_{E/F}} ; \chi_{\omega_{E/F}} \rtimes \lambda'). \;  \mid \; \mid
 \chi_{\omega_{E/F}} \times \chi_{\omega_{E/F}} \rtimes \lambda' $ is reducible by \cite[Th. 4.5]{CS}, $
 \Lg(\mid \; \mid \chi_{\omega_{E/F}} \times 
\chi_{\omega_{E/F}} ; \lambda') $ is non-unitary by the foregoing argument. By \cite{MR0324429} and by \cite[Lemma 3.6]{CS}
$ \Lg(\mid \; \mid^{\alpha}
 \chi_{\omega_{E/F}} ; \chi_{\omega_{E/F}} \rtimes \lambda') $ is non-unitary for $ 1/2 < \alpha < 1. $
\end{Remark}

\medskip

\bigskip

\begin{center} $ \Lg( \mid \; \mid^{\alpha_1} \chi_{1_{F^*}} ;
\mid \; \mid^{\alpha_2} \chi_{1_{F^*}} ; \lambda'), 0 < \alpha_2 \leq \alpha_1 $ \end{center}

\bigskip

\begin{Theorem}
 Let $ \chi_{1_{F^*}} \in X_{1_{F^*}}. $

1. Let $ 0 < \alpha_2 \leq \alpha_1, \; \alpha_1 + \alpha_2 > 1. $ Then 
$ \Lg( \mid \; \mid^{\alpha_1} \chi_{1_{F^*}} ;
\mid \; \mid^{\alpha_2} \chi_{1_{F^*}} ; \lambda') $ is non-unitary.

2. Let $ 0 < \alpha_2 \leq \alpha_1 , \; \alpha_1 + \alpha_2 = 1. $ Then $ \Lg( \mid \; \mid^{\alpha_1} \chi_{1_{F^*}} ;
\mid \; \mid^{\alpha_2} \chi_{1_{F^*}} ; \lambda') $ is unitary.

\end{Theorem}

\begin{Proof}

1. Let $ \alpha_1 - \alpha_2 = 1. $ Then $ \mid \; \mid^{\alpha_1} \chi_{1_{F^*}} \times \mid \; \mid^{\alpha_2}
\chi_{1_{F^*}} \rtimes \lambda' $ is reducible by \cite[Th. 4.3]{CS}.
Let $ \alpha: = \frac{\alpha_1 + \alpha_2}{2}, $ then $ \alpha > 1/2. $  The subquotient $ \mid \; \mid^{\alpha} 
\chi_{1_{F^*}}
1_{\GL_2} \rtimes \lambda' $ of $ \mid \; \mid^{\alpha_1} \chi_{1_{F^*}} \times \mid \; \mid^{\alpha_2} \chi_{1_{F^*}} 
\rtimes \lambda'$ is irreducible by  \cite[Th. 4.7]{CS} and equal to the Langlands quotient $ \Lg( \mid \; \mid^{\alpha_1} \chi_{1_{F^*}} ; 
\mid \; \mid^{\alpha_2} \chi_{1_{F^*}} ; \lambda'). $ Representations $ \mid \; \mid^{\alpha} \chi_{1_{F^*}}
1_{\GL_2} \rtimes \lambda' $ form a continuous 1 - parameter family
of irreducible hermitian representations that we realise on the same vector space $ V. $ If there existed $ \alpha > 1/2 $
such that $ \mid \; \mid^{\alpha} \chi_{1_{F^*}} 1_{\GL_2} \rtimes
\lambda' $ was unitary, then by \cite[Lemma 3.6]{CS} $ \mid \; \mid^{\alpha} \chi_{1_{F^*}} 1_{\GL_2} \rtimes
\lambda' $ would be unitary $ \forall \alpha > 1/2, $ in contradiction to \cite[Lemma 3.8]{CS}.

\medskip 

The representations  $ \mid \; \mid^{\alpha_1} \chi_{1_{F^*}} \times \mid \; \mid^{\alpha_2}
\chi_{1_{F^*}} \rtimes \lambda' $ in II and III of Figure ~\ref{figure3}, page ~\pageref{figure3}, are irreducible by
 \cite[Th. 4.3]{CS} and equal to their own 
Langlands quotient
$ \Lg(\mid \; \mid^{\alpha_1} \chi_{1_{F^*}} ; \mid \; \mid^{\alpha_2}
\chi_{1_{F^*}} ; \lambda'). $ Similar as in Theorem \ref{Lg11} (2) we show that
$ \Lg(\mid \; \mid^{\alpha_1} \chi_{1_{F^*}} ; \mid \; \mid^{\alpha_2} \chi_{1_{F^*}} ; \lambda') $ is non-unitary if
 $ \alpha_1 + \alpha_2 > 1 $ and $ \alpha_1 - \alpha_2 < 1, $
and if $ \alpha_1 + \alpha_2 > 1 $ and $ \alpha_1 - \alpha_2 > 1 $ (\cite[Lemmas 3.6 and 3.8]{CS}).

\bigskip
 
2. $ \mid \; \mid^{\alpha_1} \chi_{1_{F^*}} \times \mid \; \mid^{\alpha_2} \chi_{1_{F^*}} \rtimes \lambda' $ is reducible by
\cite[Th. 4.3]{CS}.
Let $ \alpha: = \frac{\alpha_1 - \alpha_2}{2}. $ The subquotient 
$ \mid \; \mid^{\alpha} \chi_{1_{F^*}}
1_{\GL_2} \rtimes \lambda' $ of $ \mid \; \mid^{\alpha_1} \chi_{1_{F^*}} \times \mid \; \mid^{\alpha_2} 
\chi_{1_{F^*}} \rtimes \lambda' $ is irreducible by \cite[Th. 6.7]{CS} and equal to the Langlands quotient
$ \Lg( \mid \; \mid^{\alpha_1} \chi_{1_{F^*}} ; 
\mid \; \mid^{\alpha_2} \chi_{1_{F^*}} ; \lambda'). $ These representations form a continuous 1 - parameter family
of irreducible hermitian representations that we realise on the same vector space $ V. \; \chi_{1_{F^*}} 1_{\GL_2} \rtimes
\lambda' $ is irreducible by \cite[Prop. 4.8]{CS}. It is unitary. For $ \alpha = 1/2, \; \mid \; \mid^{1/2} \chi_{1_{F^*}}
1_{\GL_2} \rtimes \lambda' $ reduces for the first time (\cite[Th. 4.5]{CS}). By \cite[3.6]{CS} the representations 
$ \mid \; \mid^{\alpha} \chi_{1_{F^*}} 1_{\GL_2} \rtimes \lambda' = \Lg( \mid \; \mid^{\alpha_1} \chi_{1_{F^*}} ;
\mid \; \mid^{\alpha_2} \chi_{1_{F^*}} ; 
\lambda') $ are unitary for $ 0 \leq \alpha < 1/2, $ i.e. for $ 0 < \alpha_2 \leq \alpha_1, \; \alpha_1 + \alpha_2 = 1. $
\end{Proof}

\begin{Remark}
Let $ 0 < \alpha_2 \leq \alpha_1, \; \alpha_1 + \alpha_2 < 1. $ We do not have a proof that
$ \Lg( \mid \; \mid^{\alpha_1} \chi_{1_{F^*}} ;
\mid \; \mid^{\alpha_2} \chi_{1_{F^*}} ; 
\lambda') $ is non-unitary.
\end{Remark}

\bigskip

\begin{figure}

\begin{center}
\begin{tikzpicture}
 \draw[->] (0,0) -- (3,0);
\draw (3,0) node[below right] {$\alpha_1$};
\draw[->] (0,0) -- (0,3);
\draw (0,3) node[above left] {$ \alpha_2 $};

\draw [blue] (0,0) -- (1,0) -- (0,1) --cycle;

\draw [red] (0,1) -- (2,3);
\draw [red] (1,0) -- (3,2);

\draw (1,0) node[below] {$ 1 $};
\draw (0,1) node[left] {$ 1 $};
\draw (2,0) node[below] {$ 2 $};
\draw (0,2) node[left] {$ 2 $};
\draw (1,0) node[blue] {$ \bullet $};

\draw (1,0) node[blue, above] {$ \cdot \cdot \cdot \cdot $};
\draw (0.4, -0.1) node[blue, above] {$ \cdot \cdot $};
\draw (0.55, 0.4) node[blue, above] { $ \cdot \cdot$};
\draw (2, 0.6) node[above] { $ \cdot \cdot$};
\draw (2, -0.1) node[above] { $ \cdot \cdot $};

\draw (0.3,0.3) node {I};
\draw (2.2,0.5) node {II};
\draw (1.1, 0.8) node {III};

\draw (1.5, -0.5) node[below] { Figure 3 };
\draw (3, 2.2) node[red, right] { $ - $};
\draw (3.4, 2.2) node[right] {reducible; irreducible subquotients non-unitary};
\draw (3, 1.8) node[blue, right] { $ - $};
\draw (3.4, 1.8) node[right] {reducible; irreducible subquotients unitary};
\draw (3.1,1.4) node[right] { $ \cdot \; $ irreducible non-unitary subquotient };
\draw (3.1, 1) node[blue, right] { $ \cdot $};
\draw (3.4, 1) node[right] { irreducible unitary subquotient}; 
\end{tikzpicture}
\end{center}

\caption[$ \Lg(\mid \; \mid^{\alpha_1} \chi_{1_{F^*}} ; \mid \; \mid^{\alpha_2} \chi_{1_{F^*}} ; \lambda'), \;
\alpha_1, \; \alpha_2 \geq 0, \; \chi_{1_{F^*}} \in X_{1_{F^*}} $]{Let $ \alpha_1, \alpha_2 \geq 0, $ let $ \chi_{1_{F^*}}
\in X_{1_{F^*}}. $ Figure 3 shows lines and points of reducibility of
 the representation 
$ \mid \; \mid^{\alpha_1} \chi_{1_{F^*}} \times \mid \; \mid^{\alpha_2} \chi_{1_{F^*}} \rtimes \lambda' $ and the unitary
dual. Let $ 0 < \alpha_2 \leq \alpha_1. \; \Lg(\mid \; \mid^{\alpha_1} \chi_{1_{F^*}} ; \mid \; \mid^{\alpha_2}
 \chi_{1_{F^*}} ; \lambda') $ is unitary for $ \alpha_1 + \alpha_2 = 1. $ It is non-unitary for $ \alpha_1 + \alpha_2 > 1. $}

\label{figure3}

\end{figure}

\bigskip

In the following Theorems \ref{Lg1chiomega}, \ref{Lg1chi1F*} and \ref{Lgchiomegachi1F*}, when speaking of the
Langlands-quotient, 
we will exceptionally allow that $ \alpha_1 < \alpha_2 $ for ease of notation.

\bigskip

\begin{center} $ \Lg(\mid \; \mid^{\alpha_1} 1 ; \mid \; \mid^{\alpha_2} \chi_{\omega_{E/F}} ;
 \lambda'), \; \alpha_1, \alpha_2 > 0, \;
 \Lg(\mid \; \mid^{\alpha} 1 ; \chi_{\omega_{E/F}} \rtimes
 \lambda'), \; \Lg(\mid \; \mid^{\alpha} \chi_{\omega_{E/F}} ; 1 \rtimes \lambda'), \; \alpha > 0, \; \chi_{\omega_{E/F}} 
\in X_{\omega_{E/F}}
$ \end{center}

\bigskip

Let $ \chi_{\omega_{E/F}} \in X_{\omega_{E/F}}. $ Let $ \pi_{1, \chi_{\omega_{E/F}}} $ denote the unique square-integrable
 subquotient and let $ \pi_{2, \chi_{\omega_{E/F}}} $
denote the unique non-tempered
subquotient of $ \mid \; \mid^{1/2} \chi_{\omega_{E/F}} \rtimes \lambda'. $

\medskip

\begin{Theorem}

\label{Lg1chiomega}

Let $ \chi_{\omega_{E/F}} \in X_{\omega_{E/F}}. $ 

 1. Let $ 0 < \alpha_1 \leq 1, \; 0 < \alpha_2  \leq 1/2. \; \Lg( \mid \; \mid^{\alpha_1} 1 ; \mid \; 
\mid^{\alpha_2} \chi_{\omega_{E/F}} ; \lambda') $ is unitary.

2. Let $ \alpha_1 > 1, \; \alpha_2 > 0, $ or let $ 0 < \alpha_1 \leq 1, \; \alpha_2 > 1/2. \;
\Lg( \mid \; \mid^{\alpha_1} 1 ; \mid \; \mid^{\alpha_2} \chi_{\omega_{E/F}} ; \lambda') $ is non-unitary.

3. Let $ 0 < \alpha \leq 1. \; \Lg( \mid \; \mid^{\alpha} 1 ; \chi_{\omega_{E/F}} \rtimes \lambda') $
is unitary.

4. Let $ \alpha > 1. \;  \Lg( \mid \; \mid^{\alpha} 1 ; \chi_{\omega_{E/F}} \rtimes \lambda') $ is non-unitary.

5. Let $ \alpha \leq 1/2. \; \Lg(\mid \; \mid^{\alpha} \chi_{\omega_{E/F}}; 1 \rtimes \lambda') $ is unitary.

6. Let $ \alpha > 1/2. \; \Lg(\mid \; \mid^{\alpha} \chi_{\omega_{E/F}}; 1 \rtimes \lambda') $ is non-unitary.
\end{Theorem}

\begin{Proof}
1. $ 1 \times \chi_{\omega_{E/F}} \rtimes \lambda' $ is irreducible by \cite[Th. 4.1]{CS} and unitary. 
$ \mid \; \mid^{\alpha_1} 1
 \times \mid \; \mid^{\alpha_2} \chi_{\omega_{E/F}} \rtimes \lambda' $ is irreducible by
\cite[Th. 4.3]{CS} and
equal to its own Langlands-quotient $ \Lg( \mid \; \mid^{\alpha_1} 1 \times \mid \; \mid^{\alpha_2} \chi_{\omega_{E/F}} 
\rtimes \lambda'). $ Like in \ref{Lg11} we construct continuous one-parameter families of irreducible hermitian representations
and obtain that $ \Lg( \mid \; \mid^{\alpha_1} 1 \times \mid \; \mid^{\alpha_2} \chi_{\omega_{E/F}} 
\rtimes \lambda') $ is unitary for $ 0 < \alpha_1 < 1, \; 0 < \alpha_2 < 1/2 $ (\cite[Lemma 3.6]{CS}).

\smallskip

For $ \alpha_1 = 1, \; 0 < \alpha_2 \leq 1/2, \; \mid \; \mid 1 \times \mid \; \mid^{\alpha_2} 
\chi_{\omega_{E/F}} \rtimes \lambda' $ reduces for the first time (\cite[Th. 4.3]{CS}). For 
$ 0 < \alpha_1 < 1, \; \alpha_2 = 1/2, \;
\mid \; \mid^{\alpha_ 1} 1 \times \mid \; \mid^{1/2} \chi_{\omega_{E/F}} \rtimes \lambda' $ reduces for the first time
(\cite[Th. 4.3]{CS}).
By \cite{MR0324429} $ \Lg( \mid \; \mid^{\alpha_1} 1 \times \mid \; \mid^{\alpha_2} 
\chi_{\omega_{E/F}} \rtimes \lambda') $ is unitary for $ \alpha_1 = 1,\;  0 < \alpha_2 \leq 1/2, $ and for
$ 0 < \alpha_1 < 1, \; \alpha_2 = 1/2. $

\bigskip

2. Let $ \alpha_2 = 1/2. $ Let $ \alpha > 1. $ The representation $ \mid \; \mid^{\alpha} 1 \times \mid \; \mid^{1/2}
\chi_{\omega_{E/F}} \rtimes \lambda' $ is reducible by \cite[Th. 4.3]{CS}. Its subquotient 
$ \mid \; \mid^{\alpha} 1 \rtimes 
\pi_{2, \chi_{\omega_{E/F}}} $ is
irreducible by \cite[Th. 4.11]{CS} and equal to 
$ \Lg(\mid \; \mid^{\alpha} 1; \mid \; \mid^{1/2} \chi_{\omega_{E/F}} ; \lambda'). $
The argument
is similar to \ref{Lg11};
if there existed $ \alpha > 1, $ such that $ \Lg(\mid \; \mid^{\alpha} 1; \mid \; \mid^{1/2} 
\chi_{\omega_{E/F}} ; \lambda') $  was unitary, then by \cite[Lemma 3.6]{CS} $ \Lg(\mid \; \mid^{\alpha} 1; \mid \; \mid^{1/2}
 \chi_{\omega_{E/F}} ; \lambda') $ would be unitary $ \forall \alpha > 1, $ in contradiction to \cite[Lemma 6.8]{CS}.

\smallskip

Let $ \alpha_1 = 1. $ Let $ \alpha > 1/2. $ The representation $ \mid \; \mid 1 \times \mid \; \mid^{\alpha} 
\chi_{\omega_{E/F}} \rtimes \lambda' $ is reducible by \cite[Th. 4.3]{CS}. The subquotient $ \mid \; \mid^{\alpha} 
\chi_{\omega_{E/F}}
\rtimes \lambda'(\det) 1_{U(3)} $ is irreducible by \cite[Th. 4.9]{CS} and equal to 
$ \Lg( \mid \; \mid 1 ; \mid \; \mid^{\alpha} 
\chi_{\omega_{E/F}} ; \lambda'). $
By \cite[Lemmas 3.6 and 3.8]{CS} $ \Lg( \mid \; \mid 1 ; \mid \; \mid^{\alpha} 
\chi_{\omega_{E/F}} ; \lambda') $ is non-unitary for $ \alpha > 1/2 .$

\smallskip

Let $ \alpha_1 > 1, \; 1/2 \neq \alpha_2 > 0, $ or let $ 0 < \alpha_1 < 1, \; \alpha_2 > 1/2. $ Representations $ \mid \;
\mid^{\alpha_1} 1 \times \mid \; \mid^{\alpha_2} \chi_{\omega_{E/F}} \rtimes \lambda' $ are irreducible by 
\cite[Th. 4.3]{CS} and equal to
their own Langlands-quotient $ \Lg( \mid \;
\mid^{\alpha_1} 1 ; \mid \; \mid^{\alpha_2} \chi_{\omega_{E/F}} ; \lambda'). $ Like in \ref{Lg11} we obtain that these 
Langlands-quotients are non-unitary (II, III and IV in Figure ~\ref{figure4}, page ~\pageref{figure4}).

\bigskip

3 +4. $  1 \times \chi_{\omega_{E/F}} \rtimes \lambda' $ is irreducible by \cite[Th. 4.1]{CS} and unitary. 
$ \mid \; \mid^{\alpha} 1 \times
\chi_{\omega_{E/F}} \rtimes \lambda' $ is irreducible for $ 0 < \alpha < 1 $ by \cite[Th. 4.5]{CS} and equal to its own 
Langlands -quotient
$ \Lg(\mid \; \mid^{\alpha} 1 ; \chi_{\omega_{E/F}} \rtimes \lambda'). $ By \cite[Lemma 3.6]{CS} these subquotients are 
unitary. The representation $ \mid \; \mid 1 \times \chi_{\omega_{E/F}} \rtimes \lambda' $ is reducible by 
\cite[Th. 4.5]{CS}.
The subquotient $ \chi_{\omega_{E/F}} \rtimes \lambda'(\det) 1_{U(3)} $ is irreducible by \cite[Prop. 4.10]{CS} and equal
 to
$ \Lg(\mid \; \mid 1 ; \chi_{\omega_{E/F}} \rtimes \lambda'). $ It is unitary.
By \cite[Th. 4.5]{CS} representations $ \mid \; \mid^{\alpha} 1 \times \chi_{\omega_{E/F}} \rtimes \lambda' $
 are irreducible for $ \alpha > 1 $
and equal
to their own Langlands quotient $ \Lg(\mid \; \mid^{\alpha} 1 ; \chi_{\omega_{E/F}} \rtimes \lambda'). $ By
\cite[Lemmas 3.6 and 3.8]{CS} these Langlands-quotients are non-unitary.

\smallskip

5+6. $  1 \times \chi_{\omega_{E/F}} \rtimes \lambda' $ is irreducible by \cite[Th. 4.1]{CS} and unitary.
By \cite[Th. 4.5]{CS} $ 1 \times \mid \; \mid^{\alpha} \chi_{\omega_{E/F}} \rtimes \lambda' $ is irreducible for
 $ 0 < \alpha < 1/2 $ 
and equal to its own Langlands -quotient
$ \Lg(\mid \; \mid^{\alpha} \chi_{\omega_{E/F}}; 1 \rtimes \lambda'). $ By \cite[Lemma 3.6]{CS} these subquotients are 
unitary. By \cite[Th. 4.5]{CS} the representation $ 1 \times \mid \; \mid^{1/2} \chi_{\omega_{E/F}} \rtimes \lambda' $
is reducible. Its
subquotient $ 1 \rtimes \pi_{2, \chi_{\omega_{E/F}}} $ is irreducible by \cite[Prop. 4.12]{CS} and equal to the
 Langlands-quotient
 $ \Lg( \mid \; \mid^{1/2} 
\chi_{\omega_{E/F}},  1 \rtimes \lambda'). $ It is unitary.
Representations $ 1 \times \mid \; \mid^{\alpha} \chi_{\omega_{E/F}} \rtimes \lambda' $ are irreducible for $ \alpha > 
1/2 $ by \cite[Th. 4.5]{CS}
and equal
to their own Langlands quotient $ \Lg(\mid \; \mid^{\alpha} \chi_{\omega_{E/F}} ; 1 \rtimes \lambda'). $ 
By \cite[Lemmas 3.6 and 3.8]{CS} these Langlands-quotients are non-unitary.
\end{Proof}

\medskip

\begin{figure}[!h]

\begin{center}
\begin{tikzpicture}
 \draw[->] (0,0) -- (4,0);
\draw (4,0) node[below right] {$\alpha_1$};
\draw[->] (0,0) -- (0,4);
\draw (0,4) node[above left] {$ \alpha_2 $};

\draw [fill=blue!50] (0,0) -- (2,0) -- (2,1) -- (0,1) -- cycle;

\draw [red] (2 ,1) -- (4,1);
\draw [red] (2,1) -- (2,4);
\draw [blue] (2, 0) -- (2, 1);
\draw [blue] (0, 1) -- (2, 1);

\draw (1,0) node[below] {$ 1/2 $};
\draw (0, 1) node[left] {$ 1/2 $};
\draw (2,0) node[below] {$ 1 $};
\draw (0,2) node[left] {$ 1 $};
\draw (2, 1) node[blue] {$ \bullet $};
\draw (2,0.95) node[blue, above] { $ \cdot \cdot \cdot \cdot $};

\draw (1, 0.9) node[blue, above] { $ \cdot \cdot $};
\draw (2.15, 0.25) node[blue, above] { $ \cdot \cdot $};
\draw (3, 0.9) node[above]  {$ \cdot \cdot $};

\draw (2, -0.5) node[below] { Figure 4 };

\draw (4.5, 2.7) node[red, right] { $ - $};
\draw (4.9, 2.7) node[right] {reducible; irreducible subquotients non-unitary};
\draw (4.5, 2.3) node[blue, right] { $ - $};
\draw (4.9, 2.3) node[right] {reducible; irreducible subquotients unitary};
\draw (4.6,1.9) node[right] { $ \cdot \; $ irreducible non-unitary subquotient };
\draw (4.6, 1.5) node[blue, right] { $ \cdot $};
\draw (4.9, 1.5) node[right] { irreducible unitary subquotient}; 

\draw (1, 0.5) node {I};
\draw (3, 0.5) node {II};
\draw ( 3, 1.5) node {III};
\draw (1, 1.5) node {IV};

\end{tikzpicture}
\end{center}

\caption[$ \Lg(\mid \; \mid^{\alpha_1} 1 ; \mid \; \mid^{\alpha_2} \chi_{\omega_{E/F}} ;
 \lambda'), \; \alpha_1, \alpha_2 > 0, \;
 \Lg(\mid \; \mid^{\alpha} 1 ; \chi_{\omega_{E/F}} \rtimes
 \lambda'), \; \Lg(\mid \; \mid^{\alpha} \chi_{\omega_{E/F}} ; 1 \rtimes \lambda'), \; \alpha > 0, \; \chi_{\omega_{E/F}} 
\in X_{\omega_{E/F}}
$]{Let $ \alpha_1, \alpha_2 \geq 0, $ let $ \chi_{\omega_{E/F}} 
\in X_{\omega_{E/F}}. $ Figure 4 shows lines and points of reducibility of the representation 
$ \mid \; \mid^{\alpha_1} 1 \times \mid \; \mid^{\alpha_2} \chi_{\omega_{E/F}} \rtimes \lambda' $ and the unitary dual.
$ \Lg(\mid \; \mid^{\alpha_1} 1 ; \mid \; \mid^{\alpha_2} \chi_{\omega_{E/F}} ;
 \lambda') $ is unitary for $ 0 < \alpha_1 \leq 1, \; 0 < \alpha_2 \leq 1/2. $ Otherwise it is non-unitary.
$ \Lg(\mid \; \mid^{\alpha} 1 ; \; \chi_{\omega_{E/F}} \rtimes
 \lambda') $ is unitary for $ 0 < \alpha \leq 1, \; \Lg(\mid \; \mid^{\alpha} \chi_{\omega_{E/F}} ; 1 \rtimes \lambda') $
is unitary for $ 0 < \alpha \leq 1/2. $ Otherwise these Langlands-quotients are non-unitary.} 

\label{figure4}

\end{figure}

\bigskip

\begin{center} $ \Lg( \mid \; \mid^{\alpha_1} 1 ; \mid \; \mid^{\alpha_2} \chi_{1_{F^*}} ; \lambda'),
\;  \alpha_1, \alpha_2 > 0, \;  \Lg( \mid \; \mid^{\alpha} \chi_{1_{F^*}} ; 1 \rtimes \lambda'), \; \alpha > 0, \;
\chi_{1_{F^*}} \in X_{1_{F^*}} $ \end{center}

\medskip

Let $ \chi_{1_{F^*}} \in X_{1_{F^*}}. $

\begin{Theorem}

\label{Lg1chi1F*}

1. Let $ \alpha_1, \alpha_2 > 0. \; \Lg(\mid \; \mid^{\alpha_1} 1 ; \mid \; \mid^{\alpha_2} \chi_{1_{F^*}} ; \lambda') $
is non-unitary.

2. Let $ \alpha > 0. \; \Lg(\mid \; \mid^{\alpha} \chi_{1_{F^*}} ; 1 \rtimes \lambda') $ is non-unitary.
\end{Theorem}

\begin{figure}

\begin{center}
\begin{tikzpicture}
 \draw[->] (0,0) -- (3,0);
\draw (3,0) node[below right] {$\alpha_1$};
\draw[->] (0,0) -- (0,3);
\draw (0,3) node[above left] {$ \alpha_2 $};

\draw [blue] (0,0) -- (1,0);

\draw [red] (1,0) -- (3,0);
\draw [red] (1,0) -- (1,3);

\draw (1,0) node[below] {$ 1 $};
\draw (0,1) node[left] {$ 1 $};
\draw (2,0) node[below] {$ 2 $};
\draw (0,2) node[left] {$ 2 $};
\draw (1,0) node[blue] {$ \bullet $};

\draw (1,0) node[blue, above] {$ \cdot \cdot \cdot \cdot $};
\draw (0.4, -0.1) node[blue, above] {$ \cdot \cdot $};
\draw (2, -0.1) node[above] { $ \cdot \cdot $};
\draw (1.2, 1.5) node { $ \cdot \cdot $};

\draw (1.5, -0.5) node[below] { Figure 5 };
\draw (3, 2.2) node[red, right] { $ - $};
\draw (3.4, 2.2) node[right] {reducible; irreducible subquotients non-unitary};
\draw (3, 1.8) node[blue, right] { $ - $};
\draw (3.4, 1.8) node[right] {reducible; irreducible subquotients unitary};
\draw (3.1,1.4) node[right] { $ \cdot \; $ irreducible non-unitary subquotient };
\draw (3.1, 1) node[blue, right] { $ \cdot $};
\draw (3.4, 1) node[right] { irreducible unitary subquotient}; 
\end{tikzpicture}
\end{center}

\caption[$ \Lg( \mid \; \mid^{\alpha_1} 1 ; \mid \; \mid^{\alpha_2} \chi_{1_{F^*}} ; \lambda'),
\;  \alpha_1, \alpha_2 > 0, \;  \Lg( \mid \; \mid^{\alpha} \chi_{1_{F^*}} ; 1 \rtimes \lambda'), \; \alpha > 0, \;
\chi_{1_{F^*}} \in X_{1_{F^*}} $]{Let $ \alpha_1, \alpha_2 \geq 0, $ let $ \chi_{1_{F^*}} \in X_{1_{F^*}}. $ 
Figure 5 shows lines and points of reducibility of
 the representation 
$ \mid \; \mid^{\alpha_1} 1 \times \mid \; \mid^{\alpha_2} \chi_{1_{F^*}} \rtimes \lambda' $ and the unitary dual.
$ \Lg( \mid \; \mid^{\alpha_1} 1 ; \mid \; \mid^{\alpha_2} \chi_{1_{F^*}} ; \lambda') $ is non-unitary $ \forall 
\alpha_1, \alpha_2 > 0.
\; \Lg( \mid \; \mid^{\alpha} \chi_{1_{F^*}} ; 1 \rtimes \lambda') $ is non-unitary $ \forall \alpha > 0. $ }

\label{figure5}

\end{figure}
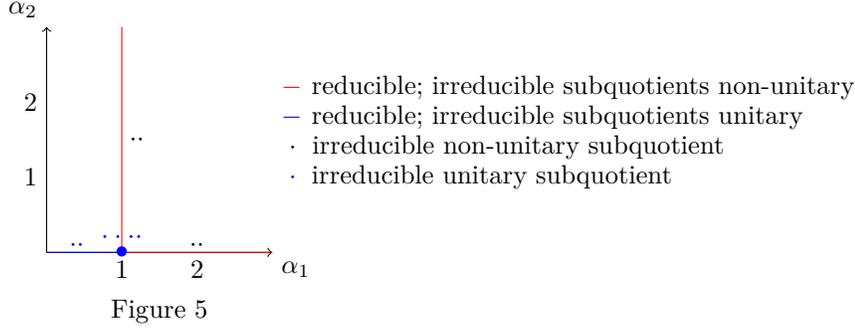

\begin{Proof} 
1. Let $ \alpha_1 = 1, \alpha_2 > 0. $ The representation $ \mid \; \mid 1 \times \mid \; \mid^{\alpha_2} \chi_{1_{F^*}} 
\rtimes \lambda' $ is reducible by \cite[Th. 4.3]{CS}. Its subquotient 
$ \mid \; \mid^{\alpha_2} \chi_{1_{F^*}} \rtimes 
\lambda'(\det) 1_{U(3)} $
is irreducible by \cite[Th. 4.9]{CS} and equal to $ \Lg(\mid \; \mid 1 ; \mid \; \mid^{\alpha_2} \chi_{1_{F^*}}
; \lambda'). $ By \cite[Lemmas 3.6 and 3.8]{CS} these Langlands-quotients are non-unitary.

\medskip

Let $ \alpha_1, \alpha_2 > 0, \; \alpha_1 \neq 1. $ Representations $ \mid \; \mid^{\alpha_1} 1 \times 
\mid \; \mid^{\alpha_2} \chi_{1_{F^*}} 
\rtimes \lambda' $ are irreducible by \cite[Th. 4.3]{CS} and equal to their own Langlands-quotient
$ \Lg(\mid \; \mid^{\alpha_1} 1 ;
 \mid \; \mid^{\alpha_2} \chi_{1_{F^*}} ; \lambda'). $ Similar as in \ref{Lg11} (2), by \cite[Lemmas 3.6 and 3.8]{CS},
we find that these Langlands-quotients are non-unitary.

\bigskip

2. Let $ \alpha > 0. $ Representations $ \mid \; \mid^{\alpha} \chi_{1_{F^*}} \times 1 \rtimes \lambda' $ are irreducible
by \cite[Th. 4.5]{CS} and equal to their own Langlands-quotient 
$  \Lg(\mid \; \mid^{\alpha} \chi_{1_{F^*}} ; 1 \rtimes \lambda'). $ By \cite[Lemmas 3.6 and 3.8]{CS}
these Langlands-quotients are non-unitary. (Figure ~\ref{figure5}, 
page ~\pageref{figure5})

\end{Proof}

\bigskip

\bigskip

\begin{center} $ \Lg( \mid \; \mid^{\alpha_1} \chi_{\omega_{E/F}} ; \mid \; \mid^{\alpha_2} \chi_{1_{F^*}} ; \lambda'), \; 
\alpha_1, \; \alpha_2 > 0, \; \Lg( \mid \; \mid^{\alpha} \chi_{1_{F^*}} ;  \chi_{\omega_{E/F}} \rtimes \lambda'), \;
\alpha > 0, \; \chi_{\omega_{E/F}} \in X_{\omega_{E/F}}, \; \chi_{1_{F^*}} \in X_{1_{F^*}} $ \end{center}

\bigskip

Let $ \chi_{\omega_{E/F}} \in X_{\omega_{E/F}}, $ let $ \chi_{1_{F^*}} \in X_{1_{F^*}}. $
\medskip

\begin{Theorem}

\label{Lgchiomegachi1F*}

 1. Let $ \alpha_1, \alpha_2 > 0. \; \Lg( \mid \; \mid^{\alpha_1} \chi_{\omega_{E/F}} ; \mid \; \mid^{\alpha_2} \chi_{1_{F^*}}
 ; \lambda') $ is non-unitary.

2. Let $ \alpha > 0. \; \Lg(\mid \; \mid^{\alpha} \chi_{1_{F^*}} ; \chi_{\omega_{E/F}} \rtimes \lambda') $ is non-unitary.
\end{Theorem}

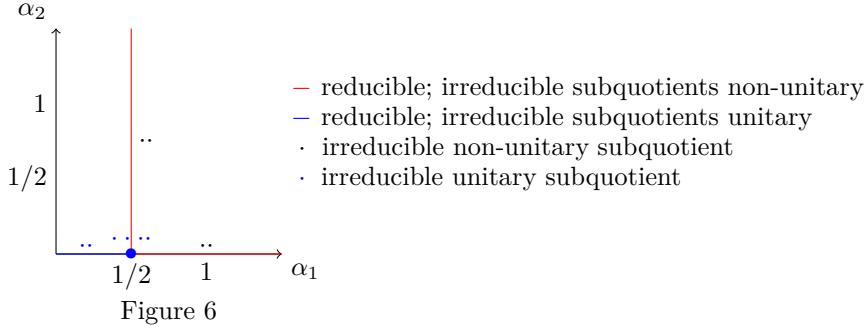
\begin{figure}[!h] 

\begin{center}
\begin{tikzpicture}
 \draw[->] (0,0) -- (3,0);
\draw (3,0) node[below right] {$\alpha_1$};
\draw[->] (0,0) -- (0,3);
\draw (0,3) node[above left] {$ \alpha_2 $};

\draw [blue] (0,0) -- (1,0);

\draw [red] (1,0) -- (3,0);
\draw [red] (1,0) -- (1,3);

\draw (1,0) node[below] {$ 1/2 $};
\draw (0,1) node[left] {$ 1/2 $};
\draw (2,0) node[below] {$ 1 $};
\draw (0,2) node[left] {$ 1 $};
\draw (1,0) node[blue] {$ \bullet $};

\draw (1,0) node[blue, above] {$ \cdot \cdot \cdot \cdot $};
\draw (0.4, -0.1) node[blue, above] {$ \cdot \cdot $};
\draw (2, -0.1) node[above] { $ \cdot \cdot $};
\draw (1.2, 1.5) node { $ \cdot \cdot $};

\draw (1.5, -0.5) node[below] { Figure 6 };
\draw (3, 2.2) node[red, right] { $ - $};
\draw (3.4, 2.2) node[right] {reducible; irreducible subquotients non-unitary};
\draw (3, 1.8) node[blue, right] { $ - $};
\draw (3.4, 1.8) node[right] {reducible; irreducible subquotients unitary};
\draw (3.1,1.4) node[right] { $ \cdot \; $ irreducible non-unitary subquotient };
\draw (3.1, 1) node[blue, right] { $ \cdot $};
\draw (3.4, 1) node[right] { irreducible unitary subquotient}; 
\end{tikzpicture}
\end{center}

\caption[$ \Lg( \mid \; \mid^{\alpha_1} \chi_{\omega_{E/F}} ; \mid \; \mid^{\alpha_2} \chi_{1_{F^*}} ; \lambda'), \; 
\alpha_1, \; \alpha_2 > 0, \; \Lg( \mid \; \mid^{\alpha} \chi_{1_{F^*}} ;  \chi_{\omega_{E/F}} \rtimes \lambda'), \;
\alpha > 0, \; \chi_{\omega_{E/F}} \in X_{\omega_{E/F}}, \; \chi_{1_{F^*}} \in X_{1_{F^*}} $]{Let 
$ \alpha_1, \alpha_2 \geq 0, $ let $ \chi_{\omega_{E/F}} \in X_{\omega_{E/F}}, \chi_{1_{F^*}} \in X_{1_{F^*}}. $ Figure 6
 shows lines and points of reducibility of the representation 
$ \mid \; \mid^{\alpha_1} \chi_{\omega_{E/F}} \times \mid \; \mid^{\alpha_2} \chi_{1_{F^*}} \rtimes \lambda' $ and the 
unitary dual. $ \Lg( \mid \; \mid^{\alpha_1} \chi_{\omega_{E/F}} ; \mid \; \mid^{\alpha_2} \chi_{1_{F^*}} ; \lambda') $
is non-unitary $ \forall \alpha_1, \alpha_2 > 0. \; \Lg( \mid \; \mid^{\alpha} \chi_{1_{F^*}} ; \chi_{\omega_{E/F}} 
\rtimes \lambda') $ is non-unitary $ \forall \alpha > 0. $ } 

\label{figure6}

\end{figure}

\begin{Proof} 
1. Let $ \alpha_1 = 1/2, \alpha_2 > 0. $ The representation $ \mid \; \mid^{1/2} \chi_{\omega_{E/F}} \times
 \mid \; \mid^{\alpha_2} \chi_{1_{F^*}} 
\rtimes \lambda' $ is reducible by \cite[Th. 4.3]{CS}. Let $ \pi_{2, \chi_{\omega_{E/F}}} $ be the unique
irreducible non-tempered subquotient of $ \mid \; \mid^{1/2} \chi_{\omega_{E/F}} \rtimes \lambda'. $
The subquotient $ \mid \; \mid^{\alpha_2} \chi_{1_{F^*}} \rtimes \pi_{2, \chi_{\omega_{E/F}}} $ of $ 
\mid \; \mid^{1/2} \chi_{\omega_{E/F}}
\times \mid \; \mid^{\alpha_2} \chi_{1_{F^*}} 
\rtimes \lambda' $
is irreducible by \cite[Th. 4.11]{CS} and equal to $ \Lg(\mid \; \mid^{1/2} \chi_{\omega_{E/F}} ; 
\mid \; \mid^{\alpha_2} \chi_{1_{F^*}} ; 
\lambda'). $ By \cite[Lemmas 3.6 and 3.8]{CS} these Langlands-quotients are non-unitary.

\medskip

Let $ \alpha_1, \alpha_2 > 0, \; \alpha_1 \neq 1/2. $ Representations $ \mid \; \mid^{\alpha_1} \chi_{\omega_{E/F}}
 \times 
\mid \; \mid^{\alpha_2} \chi_{1_{F^*}} 
\rtimes \lambda' $ are irreducible by \cite[Th. 4.3]{CS} and equal to their own Langlands-quotient 
$ \Lg(\mid \; \mid^{\alpha_1} \chi_{\omega_{E/F}} ;
 \mid \; \mid^{\alpha_2} \chi_{1_{F^*}} ; \lambda'). $ Similar as in \ref{Lg11} (2), by \cite[Lemmas 3.6 and 3.8]{CS},
 we find that these Langlands-quotients are non-unitary.

\bigskip

2. Let $ \alpha > 0. $ Representations $ \mid \; \mid^{\alpha} \chi_{1_{F^*}} \times \chi_{\omega_{E/F}} \rtimes \lambda'
 $ are irreducible by \cite[Th. 4.5]{CS}
and equal to their own Langlands-quotient $  \Lg(\mid \; \mid^{\alpha} \chi_{1_{F^*}} ; \chi_{\omega_{E/F}} \rtimes 
\lambda'). $ By \cite[Lemmas 3.6 and 3.8]{CS} these Langlands-quotients are non-unitary.
(Figure ~\ref{figure6}, page ~\pageref{figure6})
\end{Proof}

\bigskip

\bigskip

\begin{center} $ \Lg(\mid \; \mid^{\alpha_1} \chi_{\omega_{E/F},1} ; \mid \; \mid^{\alpha_2} \chi_{\omega_{E/F},2} ;
 \lambda'), \; 0 < \alpha_2 \leq \alpha_1, \; \Lg(\mid \; \mid^{\alpha} \chi_{\omega_{E/F},1} ; \chi_{\omega_{E/F},2}
 \rtimes
 \lambda'), \alpha > 0, \; \chi_{\omega_{E/F},1}, \chi_{\omega_{E/F},2} \in X_{\omega_{E/F}}, \;
\chi_{\omega_{E/F},1} \ncong \chi_{\omega_{E/F},2}  $ \end{center}

\bigskip

\begin{Theorem}

\label{Lgchiomega12}

Let $ \chi_{\omega_{E/F},1}, \chi_{\omega_{E/F},2} \in X_{\omega_{E/F}}, $ such that $ \chi_{\omega_{E/F},1} \ncong
\chi_{\omega_{E/F},2}. $

 1. Let $ 0 < \alpha_2 \leq \alpha_1 \leq 1/2. \; \Lg( \mid \; \mid^{\alpha_1} \chi_{\omega_{E/F},1} ; \mid \; 
\mid^{\alpha_2} \chi_{\omega_{E/F},2} ; \lambda') $ is unitary.

2. Let $ \alpha_1 > 1/2, \; 0 < \alpha_2 \leq \alpha_1. \; \Lg( \mid \; \mid^{\alpha_1} \chi_{\omega_{E/F},1} ; \mid \; 
\mid^{\alpha_2} \chi_{\omega_{E/F},2} \rtimes \lambda') $ is non-unitary.

3. Let $ 0 < \alpha \leq 1/2. \; \Lg( \mid \; \mid^{\alpha} \chi_{\omega_{E/F},1} ; \chi_{\omega_{E/F},2} \rtimes \lambda') $
is unitary.

4. Let $ \alpha > 1/2. \;  \Lg( \mid \; \mid^{\alpha} \chi_{\omega_{E/F},1}  ; \chi_{\omega_{E/F},2} \rtimes \lambda') $
 is non-unitary.
\end{Theorem}

\begin{figure}[!h]
 
\begin{center}
\begin{tikzpicture}
 \draw[->] (0,0) -- (4,0);
\draw (4,0) node[below right] {$\alpha_1$};
\draw[->] (0,0) -- (0,4);
\draw (0,4) node[above left] {$ \alpha_2 $};

\draw [fill=blue!50] (0,0) -- (1,0) -- (1,1) -- (0,1) -- cycle;

\draw [red] (1 ,1) -- (4,1);
\draw [red] (1,1) -- (1,4);
\draw [blue] (1, 0) -- (1, 1);
\draw [blue] (0, 1) -- (1, 1);

\draw (1,0) node[below] {$ 1/2 $};
\draw (0, 1) node[left] {$ 1/2 $};
\draw (2,0) node[below] {$ 1 $};
\draw (0,2) node[left] {$ 1 $};
\draw (1, 1) node[blue] {$ \bullet $};
\draw (1,0.95) node[blue, above] { $ \cdot \cdot \cdot \cdot $};

\draw (1.15, 0.25) node[blue, above] { $ \cdot \cdot $};
\draw (2.5, 0.9) node[above]  {$ \cdot \cdot $};

\draw (2, -0.5) node[below] { Figure 7 };
\draw (4.5, 2.7) node[red, right] { $ - $};
\draw (4.9, 2.7) node[right] {reducible; irreducible subquotients non-unitary};
\draw (4.5, 2.3) node[blue, right] { $ - $};
\draw (4.9, 2.3) node[right] {reducible; irreducible subquotients unitary};
\draw (4.6,1.9) node[right] { $ \cdot \; $ irreducible non-unitary subquotient };
\draw (4.6, 1.5) node[blue, right] { $ \cdot $};
\draw (4.9, 1.5) node[right] { irreducible unitary subquotient}; 

\draw (0.5, 0.5) node {I};
\draw (2.5, 0.5) node {II};
\draw ( 2.5, 1.5) node {III};

\end{tikzpicture}
\end{center}

\caption[$ \Lg(\mid \; \mid^{\alpha_1} \chi_{\omega_{E/F},1} ; \mid \; \mid^{\alpha_2} \chi_{\omega_{E/F},2} ;
 \lambda'), \; 0 < \alpha_2 \leq \alpha_1, \; \Lg(\mid \; \mid^{\alpha} \chi_{\omega_{E/F},1} ; \chi_{\omega_{E/F},2}
 \rtimes
 \lambda'), \alpha > 0, \; \chi_{\omega_{E/F},1}, \chi_{\omega_{E/F},2} \in X_{\omega_{E/F}}, \;
\chi_{\omega_{E/F},1} \ncong \chi_{\omega_{E/F},2}  $]{Let $ \alpha_1, \alpha_2 \geq 0, $ let $ \chi_{\omega_{E/F},1}, 
\chi_{\omega_{E/F},2} \in X_{\omega_{E/F}}, $ such that
$ \chi_{\omega_{E/F},1} \ncong \chi_{\omega_{E/F},2}. $ Figure 7 shows 
lines and points of reducibility of the representation 
$ \mid \; \mid^{\alpha_1} \chi_{\omega_{E/F},1} \times \mid \; \mid^{\alpha_2} \chi_{\omega_{E/F},2} \rtimes \lambda' $ and 
the unitary dual. $ \Lg(\mid \; \mid^{\alpha_1} \chi_{\omega_{E/F},1} ; \mid \; \mid^{\alpha_2} \chi_{\omega_{E/F},2} ;
 \lambda') $ is unitary for $ 0 < \alpha_2 \leq \alpha_1 \leq 1/2. $ Otherwise it is non-unitary. $
\Lg(\mid \; \mid^{\alpha} \chi_{\omega_{E/F},1} ; \chi_{\omega_{E/F},2}
 \rtimes
 \lambda') $ is unitary for $ 0 < \alpha \leq 1/2. $ Otherwise it is non-unitary.} 

\label{figure7}

\end{figure}
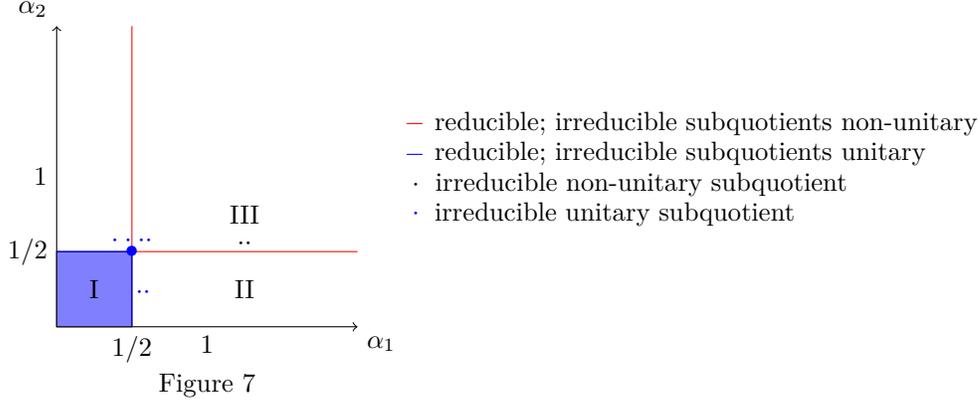

\begin{Proof} 
1. $ \chi_{\omega_{E/F},1} \times \chi_{\omega_{E/F},2} \rtimes \lambda' $ is irreducible by \cite[Th. 4.1]{CS}
 and unitary. 
Let $  \alpha_1 = 1/2,$ fix $ 0 < \alpha_2 \leq 1/2 $ and let  $ t \in [0,1]. $ Representations 
$ \mid \; \mid^{t 1/2} \chi_{\omega_{E/F},1}
\times \mid \; \mid^{t\alpha_2} \chi_{\omega_{E/F},2} \rtimes \lambda' $ give a continuous 1- parameter 
family of hermitian representations that we realise on the same vector space V, similar as in \ref{Lg11}.
For $ t \in [0,1), \; \mid \; \mid^{t 1/2} \chi_{\omega_{E/F},1}
\times \mid \; \mid^{t\alpha_2} \chi_{\omega_{E/F},2} \rtimes \lambda' $ is irreducible by \cite[Th. 4.3]{CS} and
unitary by \cite[Lemma 3.6]{CS}. For $ t = 1, \; \mid \; \mid^{1/2} \chi_{\omega_{E/F},1}
\times \mid \; \mid^{\alpha_2} \chi_{\omega_{E/F},2} \rtimes \lambda' $
reduces for the first time. By \cite{MR0324429}
$ \Lg(\mid \; \mid^{1/2} \chi_{\omega_{E/F},1}; \mid \; \mid^{\alpha_2} \chi_{\omega_{E/F},2} ; \lambda')
$ is unitary. Hence $ \Lg(\mid \; \mid^{\alpha_1} \chi_{\omega_{E/F},1} ; \mid \; \mid^{\alpha_2} \chi_{\omega_{E/F},2} ; 
\lambda') $
is unitary for $ 0 < \alpha_2 \leq \alpha_1 \leq 1/2. $

\medskip

2. Let $ \alpha_1 > 1/2, \; \alpha_2 = 1/2. \; \mid \; \mid^{\alpha_1} \chi_{\omega_{E/F},1} \times \mid \; \mid^{1/2}
\chi_{\omega_{E/F},2} \rtimes \lambda' $ is reducible by \cite[Th. 4.3]{CS}. Let
$ \pi_{2, \chi_{\omega_{E/F}}'} $ be the unique irreducible non-tempered subquotient of $ \mid \; \mid^{1/2} 
\chi_{\omega_{E/F},2} 
\rtimes \lambda'. $ By \cite[Th. 4.11]{CS} $ \mid \; \mid^{\alpha_1}
\chi_{\omega_{E/F},1} \rtimes \pi_{2, \chi_{\omega_{E/F}}'} = \Lg(\mid \; \mid^{\alpha_1} \chi_{\omega_{E/F},1} ;
\mid \; \mid^{1/2}
\chi_{\omega_{E/F},2} ; \lambda') $ is irreducible. 
For $ \alpha_1 > 1/2, \; \Lg(\mid \; \mid^{\alpha_1} \chi_{\omega_{E/F},1} ; \mid \; \mid^{1/2}
\chi_{\omega_{E/F},2} ; \lambda') $ gives a continuous 1-parameter-family of irreducible hermitian representations that we realize on the same
vector space V, similar as in \ref{Lg11}. If there existed $ \alpha_1 > 1/2 $ such that 
$ \Lg(\mid \; \mid^{\alpha_1} \chi_{\omega_{E/F},1} ; \mid \; \mid^{1/2}
\chi_{\omega_{E/F},2} ; \lambda') $ was unitary, then by \cite[Lemma 3.6]{CS} 
$ \Lg(\mid \; \mid^{\alpha_1} \chi_{\omega_{E/F},1} ;
\mid \; \mid^{1/2}
\chi_{\omega_{E/F},2} ; \lambda') $ would be unitary for all $ \alpha_1 > 1/2, $ in contradiction to
\cite[Lemma 3.8]{CS}.

\medskip

We show that representations in II of Figure ~\ref{figure7}, page ~\pageref{figure7}, are non-unitary.
Let $ \alpha_1 > 1/2 $ and $ 0 < \alpha_2 < 1/2. $ Representations $ \mid \; \mid^{\alpha_1} \chi_{\omega_{E/F},1}
\times \mid \; \mid^{\alpha_2} \chi_{\omega_{E/F},2} \rtimes \lambda' $ are irreducible by \cite[Th. 4.3]{CS}
and equal their own
Langlands quotient $ \Lg( \mid \; \mid^{\alpha_1} \chi_{\omega_{E/F},1} ; \mid \; \mid^{\alpha_2} \chi_{\omega_{E/F},2} ;
 \lambda') $. Fix $ 0 < \alpha_2 < 1/2. \; \Lg( \mid \; \mid^{\alpha_1} \chi_{\omega_{E/F},1} ; \mid \; \mid^{\alpha_2}
\chi_{\omega_{E/F},2} ;
 \lambda') $ give a continuous 1-parameter family of irreducible hermitian representations that we realize on the same vector space V.
 If there existed $ \alpha_1 > 1/2 $
such that $ \Lg( \mid \; \mid^{\alpha_1} \chi_{\omega_{E/F},1} ; \mid \; \mid^{\alpha_2} \chi_{\omega_{E/F},2} ;
 \lambda') $ was unitary, then by \cite[Lemma 3.6]{CS} representations
$ \Lg( \mid \; \mid^{\alpha_1} \chi_{\omega_{E/F},1} ; \mid \; \mid^{\alpha_2} \chi_{\omega_{E/F},2} ;
 \lambda') $ would be unitary $ \forall \alpha_1 > 1/2, $ in contradiction to \cite[Lemma 3.8]{CS}.

\medskip

We show that representations in III of Figure ~\ref{figure7}, page ~\pageref{figure7}, are non-unitary.
Let $ \alpha_1, \alpha_2 > 1/2, \; \alpha_2 \leq \alpha_1. $ The representations 
$ \mid \; \mid^{\alpha_1} \chi_{\omega_{E/F},1}
\times \mid \; \mid^{\alpha_2} \chi_{\omega_{E/F},2} \rtimes \lambda' $ are irreducible by \cite[Th. 4.3]{CS}
and equal to their own
Langlands quotient $ \Lg( \mid \; \mid^{\alpha_1} \chi_{\omega_{E/F},1} ; \mid \; \mid^{\alpha_2} \chi_{\omega_{E/F},2} ;
 \lambda') $. Fix $ \alpha_1 = \alpha_2 =: \alpha. $ Let $ t \in [1, \infty[. \; \Lg( \mid \; \mid^{t\alpha} 
\chi_{\omega_{E/F},1} ; \mid \; \mid^{
\alpha} \chi_{\omega_{E/F},2} ;
 \lambda') $ gives a continuous 1-parameter family of irreducible hermitian representations that we realize on the same vector space V,
similar as
in \ref{Lg11}. If there existed $ t \in [ 1, \infty [ $ such that  $ \Lg( \mid \; \mid^{t\alpha} 
\chi_{\omega_{E/F},1} ; \mid \; \mid^{
\alpha} \chi_{\omega_{E/F},2} ;
 \lambda') $ was unitary, then by \cite[Lemma 3.6]{CS}  $ \Lg( \mid \; \mid^{t\alpha} \chi_{\omega_{E/F},1} ; \mid \; \mid^{
\alpha} \chi_{\omega_{E/F},2} ; \lambda') $ would be unitary $ \forall t \in [ 1, \infty [, $ in contradiction to Lemma 
\cite[Lemma 3.8]{CS}.
Hence $ \Lg( \mid \; \mid^{\alpha_1} \chi_{\omega_{E/F},1} ; \mid \; \mid^{\alpha_2} \chi_{\omega_{E/F},2} ;
 \lambda') $ is non-unitary for $ \alpha_1, \alpha_2 > 1/2, \alpha_2 \leq \alpha_1. $

\medskip

3. $ \chi_{\omega_{E/F},1} \times \chi_{\omega_{E/F},2} \rtimes \lambda' $ is irreducible by \cite[Th. 4.1]{CS} and unitary.
$ \mid \; \mid^{\alpha}
\chi_{\omega_{E/F},1} \times \chi_{\omega_{E/F},2} \rtimes \lambda' $ is irreducible for $ 0 < \alpha < 1/2 $ by
\cite[Th. 4.5]{CS}
and equal to its 
Langlands-quotient
$ \Lg(\mid \; \mid^{\alpha}
\chi_{\omega_{E/F},1} ; \chi_{\omega_{E/F},2} \rtimes \lambda'). $  By \cite[Lemma 3.6]{CS} these 
Langlands-quotients are unitary. For $ \alpha = 1/2, \; \mid \; \mid^{1/2}
\chi_{\omega_{E/F},1} \times \chi_{\omega_{E/F},2} \rtimes \lambda' $ reduces for the first time. By \cite{MR0324429}
 $ \Lg(\mid \; \mid^{1/2}
\chi_{\omega_{E/F},1} ; \chi_{\omega_{E/F},2} \rtimes \lambda') $ is unitary. 

\medskip

4. For $ \alpha > 1/2, \; \mid \; \mid^{\alpha}
\chi_{\omega_{E/F},1} \times \chi_{\omega_{E/F},2} \rtimes \lambda' $ is irreducible by \cite[Th. 4.5]{CS} and equal
to its 
Langlands-quotient
$ \Lg(\mid \; \mid^{\alpha}
\chi_{\omega_{E/F},1} ; \chi_{\omega_{E/F},2} \rtimes \lambda'). $ By \cite[Lemmas 3.6 and 3.8]{CS}
these Langlands-quotients
are non-unitary. (Figure ~\ref{figure7}, page ~\pageref{figure7})
\end{Proof} 

\bigskip

\bigskip

\begin{center} $ \Lg( \mid \; \mid^{\alpha_1} \chi_{1_{F^*},1} ; \mid \; \mid^{\alpha_2} \chi_{1_{F^*},2} ; \lambda'), \;
0 < \alpha_2 \leq \alpha_1 $ \end{center}

\medskip

Let $ \chi_{1_{F^*},1}, \chi_{1_{F^*},2} \in X_{1_{F^*}}, $ such that $ \chi_{1_{F^*}} \ncong \chi_{2_{F^*}}. $

\begin{Theorem}

\label{Lgchi1F*12}

Let $ \alpha_1, \alpha_2 > 0. \; \Lg(\mid \; \mid^{\alpha_1} \chi_{1_{F^*},1} ; \mid \; \mid^{\alpha_2}
 \chi_{1_{F^*},2} ; \lambda') $ is non-unitary.
\end{Theorem}

\begin{Proof} Representations $ \mid \; \mid^{\alpha_1} \chi_{1_{F^*},1} \times \mid \; \mid^{\alpha_2}
 \chi_{1_{F^*},2} \rtimes \lambda' $ are irreducible by \cite[Th. 4.3]{CS} and equal to their own Langlands-quotient
$ \Lg(\mid \; \mid^{\alpha_1} \chi_{1_{F^*},1} ; \mid \; \mid^{\alpha_2}
 \chi_{1_{F^*},2} ; \lambda'). $ As in \ref{Lg11} (2), by \cite[Lemmas 3.6 and 3.8]{CS}, we obtain that 
$ \Lg(\mid \; \mid^{\alpha_1} \chi_{1_{F^*},1} ; \mid \; \mid^{\alpha_2}
 \chi_{1_{F^*},2} ; \lambda') $ is non-unitary. (Figure ~\ref{figure8}, page ~\pageref{figure8})
\end{Proof}

\begin{figure}

\begin{center}
\begin{tikzpicture}
 \draw[->] (0,0) -- (3,0);
\draw (3,0) node[below right] {$\alpha_1$};
\draw[->] (0,0) -- (0,3);
\draw (0,3) node[above left] {$ \alpha_2 $};

\draw [red] (0,0) -- (3,0);
\draw [red] (0,0) -- (0,3);

\draw (0,0) node[blue] {$ \bullet $};

\draw (0,0) node[blue, above] {$ \cdot \cdot \cdot \cdot $};
\draw (1.5, -0.1) node[above] { $ \cdot \cdot $};
\draw (0.2, 1.5) node { $ \cdot \cdot $};

\draw (1.5, -0.5) node[below] { Figure 8 };
\draw (3, 2.2) node[red, right] { $ - $};
\draw (3.4, 2.2) node[right] {reducible; irreducible subquotients non-unitary};
\draw (3.1,1.8) node[right] { $ \cdot \; $ irreducible non-unitary subquotient };
\draw (3.1, 1.4) node[blue, right] { $ \cdot $};
\draw (3.4, 1.4) node[right] { irreducible unitary subquotient}; 
\end{tikzpicture}
\end{center}

\caption[$ \Lg( \mid \; \mid^{\alpha_1} \chi_{1_{F^*},1} ; \mid \; \mid^{\alpha_2} \chi_{1_{F^*},2} ; \lambda'), \;
0 < \alpha_2 \leq \alpha_1 $]{Let $ \alpha_1, \alpha_2 \geq 0, $ let $ \chi_{1_{F^*},1}, \chi_{1_{F^*},2} \in X_{1_{F^*}} $ such that $ \chi_{1_{F^*},1}
\ncong \chi_{1_{F^*},2}. $ Figure 8
 shows lines and points of reducibility of
 the representation 
$ \mid \; \mid^{\alpha_1} \chi_{1_{F^*},1} \times \mid \; \mid^{\alpha_2} \chi_{1_{F^*},2} \rtimes \lambda' $ and the 
unitary dual. $ \Lg( \mid \; \mid^{\alpha_1} \chi_{1_{F^*},1} ; \mid \; \mid^{\alpha_2} \chi_{1_{F^*},2} ; \lambda') $
is non-unitary for $ 0 < \alpha_2 \leq \alpha_1. $} 

\label{figure8}

\end{figure}

\smallskip

\subsection{Representations induced from $ M_1, $ with cuspidal support in $ M_0, $ not fully-induced}

\bigskip
\smallskip

\begin{center} Irreducible subquotients of $ \chi \St_{\GL_2} \rtimes \lambda' $ \end{center}

\medskip

Let $ \chi $ be a unitary character of $ E^*. $

\medskip

 $ \chi \St_{\GL_2} \rtimes \lambda' $ is tempered, unitary. Hence all irreducible subquotients of $
\chi \St_{\GL_2} \rtimes \lambda' $ are tempered and unitary.

\medskip

\begin{Remark}
 By \cite[Prop. 4.8]{CS} $ \chi \St_{\GL_2} \rtimes \lambda' $ is reducible if and only if $ \chi =: 
\chi_{\omega_{E/F}} \in
 X_{\omega_{E/F}}. \; \chi_{\omega_{E/F}} \St_{\GL_2} \rtimes \lambda' $ has two tempered
subquotients.
\end{Remark}

\bigskip

\begin{center} $ \Lg( \mid \; \mid^{\alpha} \chi \St_{\GL_2} ; \lambda'), \; \alpha > 0 $ \end{center}

\medskip

\begin{Theorem}
    Let $ \chi $ be a unitary character of $ E^* $ such that $ \chi \notin X_{N_{E/F}(E^*)}. \; \Lg( \mid \; \mid^{\alpha}
 \chi \St_{\GL_2} ; \lambda') $ is non-unitary $ \forall
\alpha > 0. $
\end{Theorem}

\begin{Proof}
If $ \chi \notin X_{N_{E/F}(E^*}), $ the representations $ \Lg( \mid \; \mid^{\alpha} \chi \St_{\GL_2} ; \lambda') $ are
not
hermitian.
\end{Proof}

\smallskip

\begin{Theorem}
 
1. Let $ 0 < \alpha \leq 1/2. \; \Lg( \mid \; \mid^{\alpha} \St_{\GL_2} ; \lambda') $ is unitary.

2. Let $ \alpha > 1/2. \; \Lg( \mid \; \mid^{\alpha} \St_{\GL_2} ; \lambda') $ is non-unitary.
\end{Theorem}

\begin{figure}
 
\begin{center}
\begin{tikzpicture}
 \draw[->] (0,0) -- (3,0);
\draw (3,0) node[below right] {$\alpha_1$};
\draw[->] (0,0) -- (0,3);
\draw (0,3) node[above left] {$ \alpha_2 $};

\draw  (0,0) -- (1,0) -- (0,1) --cycle;
\draw [blue] (1,0) -- (0,1);

\draw  (0,1) -- (3,1);
\draw [red] (0,1) -- (2,3);
\draw  (1,0) -- (1,3);
\draw [red] (1,0) -- (3,2);

\draw (1,0) node[below] {$ 1 $};
\draw (0,1) node[left] {$ 1 $};
\draw (2,0) node[below] {$ 2 $};
\draw (0,2) node[left] {$ 2 $};
\draw (1,0) node[blue] {$ \bullet $};
\draw (2,1) node {$ \bullet $};
\draw (1,0) node[above] {$ \cdot \cdot \cdot \cdot $};
\draw (1.9,1) node[blue, above] { $ \cdot \cdot $};
\draw (2.1, 1) node[above] { $ \cdot \cdot $};
\draw (0.55, 0.4) node[above] { $ \cdot \cdot$};
\draw (1.55, 0.4) node[above] { $ \cdot \cdot$};


\draw (1.4, 0.4) node[right] { \begin{tiny} $ \Lg( \mid \; \mid^{\alpha} \St_{\GL_2} ; \lambda') $ \end{tiny}};

\draw (1.5, -0.5) node[below] { Figure 9 };
\draw (3, 2.2) node[red, right] { $ - $};
\draw (3.4, 2.2) node[right] {reducible; irreducible subquotients non-unitary};
\draw (3, 1.8) node[blue, right] { $ - $};
\draw (3.4, 1.8) node[right] {reducible; irreducible subquotients unitary};
\draw (3.1,1.4) node[right] { $ \cdot \; $ irreducible non-unitary subquotient };
\draw (3.1, 1) node[blue, right] { $ \cdot $};
\draw (3.4, 1) node[right] { irreducible unitary subquotient}; 
\end{tikzpicture}
\end{center}

\caption[$ \Lg(\mid \; \mid^{\alpha} \St_{\GL_2} \rtimes \lambda'), \; \alpha > 0 $]{
Let $ \alpha_1, \alpha_2 \geq 0. $ Figure 9 shows lines and points of reducibility of the representation 
$ \mid \; \mid^{\alpha_1} 1 \times \mid \; \mid^{\alpha_2} 1 \rtimes \lambda'. $
For $ 0 < \alpha \leq 1/2, \; \Lg(\mid \; \mid^{\alpha} \St_{\GL_2} \rtimes \lambda') $ is unitary, for $ \alpha > 1/2 $
it is non-unitary.}

\label{figure9}

\end{figure}
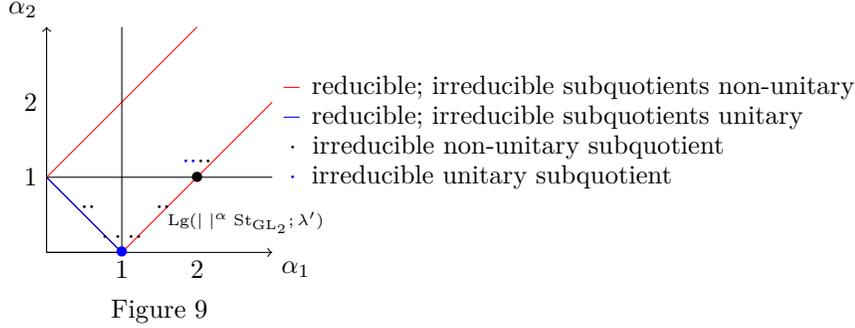

\begin{Proof}
1. Let $ 0 < \alpha < 1/2. $ The representations $ \mid \; \mid^{\alpha} \St_{\GL_2} \rtimes \lambda' $ are irreducible by
\cite[Th. 4.7]{CS}. They form a continuous
one-parameter family of irreducible hermitian representations, that, similar as in \ref{Lg11}, we realize on the
same vector space $ V. \; \St_{\GL_2} \rtimes \lambda' $ is irreducible by \cite[Prop. 4.8]{CS} and tempered, hence 
unitary. By \cite[Lemma 3.6]{CS} the representations 
$ \mid \; \mid^{\alpha} \St_{\GL_2} \rtimes \lambda' = \Lg(\mid \; \mid^{\alpha} \St_{\GL_2} ; \lambda') $ are unitary for $
0 < \alpha < 1/2. $ By \cite[Th. 4.7]{CS} $ \mid \; \mid^{1/2} \St_{\GL_2} \rtimes \lambda' $ is reducible.
 By \cite{MR0324429}
 $ \Lg(\mid \; \mid^{1/2} \St_{\GL_2} 
; \lambda') $ is unitary.

\medskip

2. Let $ 1/2 < \alpha < 3/2. $ The representations $  \mid \; \mid^{\alpha} \St_{\GL_2} \rtimes \lambda' $ are irreducible
by \cite[Th. 4.7]{CS} and equal to their Langlands quotient $ 
\Lg(\mid \; \mid^{\alpha} \St_{\GL_2} ; \lambda'). $ They form a continuous 1 - parameter family of irreducible hermitian
representations, that we realize on the same vector space $ V. $
For $ \alpha = 3/2, \; \Lg(\mid \; \mid^{3/2} \St_{\GL_2} ; \lambda') $ is a subquotient of the
representation $  \mid \; \mid^2 1 \times \mid \; \mid^1 1 \rtimes \lambda'  $ (\cite[Th. 5.2]{CS}). By results of Casselmann
\cite{MR656064}, page 915, it is non-unitary.
By \cite{MR0324429} and \cite[Lemma 6.3]{CS} $ \Lg( \mid \; \mid^{\alpha} \St_{\GL_2} ; \lambda') $ is
not unitary for $ 1/2 < \alpha < 3/2. $

The representations $ \Lg(\mid \; \mid^{\alpha} \St_{\GL_2} ; \lambda'), \; \alpha > 3/2, $ form a continuous 1-paramater
family of irreducible hermitian representations. If there existed $ \alpha > 3/2 $ such that $ \Lg(\mid \; \mid^{\alpha}
 \St_{\GL_2} ; \lambda') $ was unitary, then by \cite[Lemma 6.3]{CS} $ \Lg(\mid \; \mid^{\alpha} \St_{\GL_2} ; \lambda') $ would 
be unitary for all $ \alpha > 3/2, $ in contradiction to \cite[Lemma 6.8]{CS}. (Figure ~\ref{figure9},
page ~\pageref{figure9})
\end{Proof}
 
\smallskip

\begin{Theorem}
Let $ \chi_{\omega_{E/F}} \in X_{\omega_{E/F}}. $

Let $ \alpha > 1/2. \; \Lg(\mid \; \mid^{\alpha} \chi_{\omega_{E/F}} \St_{\GL_2} ; \lambda') $ is non-unitary.
\end{Theorem}

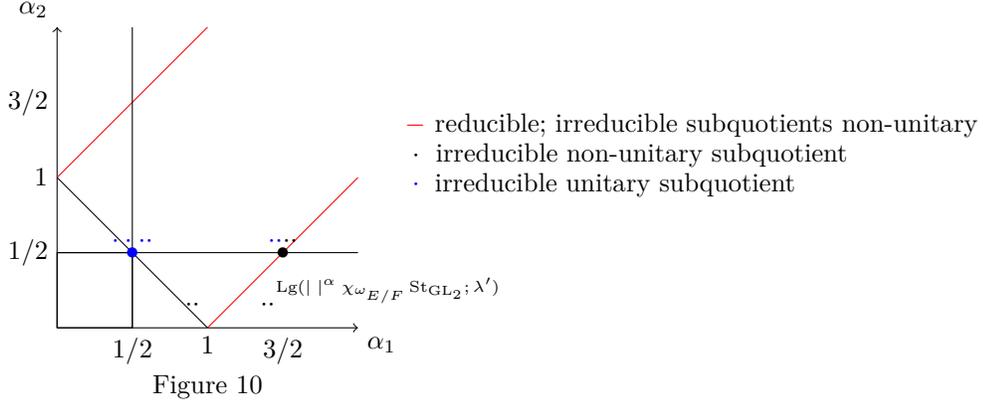
\begin{figure}
 
\begin{center}
\begin{tikzpicture}
 \draw[->] (0,0) -- (4,0);
\draw (4,0) node[below right] {$\alpha_1$};
\draw[->] (0,0) -- (0,4);
\draw (0,4) node[above left] {$ \alpha_2 $};

\draw  (0,0) -- (1,0) -- (1,1) -- (0,1) -- cycle;

\draw  (1 ,1) -- (4,1);
\draw [red] (0,2) -- (2,4);
\draw  (1,1) -- (1,4);
\draw  [red] (2,0) -- (4,2);
\draw  (2,0) -- (0,2);
\draw  (1, 0) -- (1, 1);
\draw (0, 1) -- (1, 1);

\draw (1,0) node[below] {$ 1/2 $};
\draw (0, 1) node[left] {$ 1/2 $};
\draw (2,0) node[below] {$ 1 $};
\draw (0,2) node[left] {$ 1 $};
\draw (3,0) node[below] { $ 3/2 $};
\draw (0, 3) node[left] { $ 3/2 $};
\draw (1, 1) node[blue] {$ \bullet $};
\draw (3, 1) node {$ \bullet $};
\draw (2.9,0.95) node[blue, above] {$ \cdot \cdot$};
\draw (3.1,0.95) node[above]{$ \cdot \cdot$};
\draw (1,0.95) node[blue, above] { $ \cdot \cdot \cdot \cdot $};
\draw (1.8, 0.1) node[above] { $ \cdot \cdot $};
\draw (2.8, 0.1) node[above] { $ \cdot \cdot $};

\draw (2.7, 0.5) node[right] { \begin{tiny} $ \Lg(\mid \; \mid^{\alpha} \chi_{\omega_{E/F}} \St_{\GL_2} ; \lambda') 
$ \end{tiny} };

\draw (2, -0.5) node[below] { Figure 10 };
\draw (4.5, 2.7) node[red, right] { $ - $};
\draw (4.9, 2.7) node[right] {reducible; irreducible subquotients non-unitary};
\draw (4.6,2.3) node[right] { $ \cdot \; $ irreducible non-unitary subquotient };
\draw (4.6, 1.9) node[blue, right] { $ \cdot $};
\draw (4.9, 1.9) node[right] { irreducible unitary subquotient}; 


\end{tikzpicture}
\end{center}

\caption[$ \Lg(\mid \; \mid^{\alpha} \chi_{\omega_{E/F}} \St_{\GL_2} ; \lambda'), \; \alpha > 0, \chi_{\omega_{E/F}} \in
X_{\omega_{E/F}} $]{Let $ \alpha_1, \alpha_2 \geq 0, $ let $ \chi_{\omega_{E/F}} \in X_{\omega_{E/F}}. $
Figure 10 shows lines and points of reducibility of the representation 
$ \mid \; \mid^{\alpha_1} \chi_{\omega_{E/F}} \times \mid \; \mid^{\alpha_2} \chi_{\omega_{E/F}} \rtimes \lambda'. $
Let $ \alpha > 1/2. \; \Lg(\mid \; \mid^{\alpha} \chi_{\omega_{E/F}} \St_{\GL_2} ; \lambda') $ is non-unitary.} 

\label{figure10}

\end{figure}

\begin{Proof}
Let $ 1/2 < \alpha < 1. \; \mid \; \mid^{\alpha} \chi_{\omega_{E/F}} \St_{\GL_2} \rtimes \lambda' $ is irreducible by 
\cite[Th. 4.7]{CS} and
equal to its own Langlands quotient $ \Lg(\mid \; \mid^{\alpha} \chi_{\omega_{E/F}} \St_{\GL_2} ; \lambda').$
For $ \alpha = 1, $ by \cite[Th. 5.4]{CS} $ \mid \; \mid \chi_{\omega_{E/F}} \St_{\GL_2} \rtimes \lambda' $ is
reducible, by Remark 4.7
in \cite{MR2652536}, 
$ \Lg(\mid \; \mid \chi_{\omega_{E/F}} \St_{\GL_2} ; \lambda') $ is non-unitary. By \cite{MR0324429} and 
\cite[Lemma 6.3]{CS}
$ \Lg(\mid \; \mid^{\alpha} \chi_{\omega_{E/F}} \St_{\GL_2} ; \lambda') $ is not unitary for $  1/2 < \alpha < 1. $

\medskip

Let $ \alpha > 1. \;  \mid \; \mid^{\alpha} \chi_{\omega_{E/F}} \St_{\GL_2} \rtimes \lambda' $ is irreducible by
\cite[Th. 4.7]{CS}
and equal to its own Langlands quotient $ \Lg(\mid \; \mid^{\alpha} \chi_{\omega_{E/F}} \St_{\GL_2} ; \lambda') $.
If there existed $ \alpha > 1 $ such that $ \Lg(\mid \; \mid^{\alpha} \chi_{\omega_{E/F}} \St_{\GL_2} ; \lambda') $ was
unitary, then by \cite[Lemma 6.3]{CS} $ \Lg(\mid \; \mid^{\alpha} \chi_{\omega_{E/F}} \St_{\GL_2} ; \lambda') $ would
 be unitary
$ \forall \alpha > 1, $ in contradiction to \cite[Lemma 6.8]{CS}. (Figure ~\ref{figure10}, page ~\pageref{figure10})
\end{Proof}

\begin{Remark}

\label{remarkalphachiomegast}

 Let $ 0 < \alpha < 1/2. $ Then $ \mid \; \mid^{\alpha} \chi_{\omega_{E/F}} \St_{\GL_2} \rtimes \lambda' $ is 
irreducible by \cite[Th. 4.7]{CS}
and
equal to its own Langlands quotient $ \Lg(\mid \; \mid^{\alpha} \chi_{\omega_{E/F}} \St_{\GL_2} ; \lambda').$
If we assume that $ \mid \; \mid^{1/2} \chi_{\omega_{E/F}} \St_{\GL_2} \rtimes \lambda' $ is irreducible and equal to
$ \Lg(\mid \; \mid^{1/2} \chi_{\omega_{E/F}} \St_{\GL_2} ; \lambda'), $ see \cite[Remark 5.5]{CS}, then we can extend the argument
that $ \Lg( \mid \; \mid^{\alpha}
\chi_{\omega_{E/F}} \St_{\GL_2} ; \lambda') $ is non-unitary to $ \alpha > 0. $
\end{Remark}

\bigskip

\medskip

\begin{Theorem}
Let $ \chi_{1_{F^*}} \in X_{1_{F^*}}. $

Let $ 0 < \alpha \leq 1/2. \; \Lg(\mid \; \mid^{\alpha} \chi_{1_{F^*}} \St_{\GL_2} ; \lambda') $ is unitary.

Let $ \alpha > 1/2. \; \Lg(\mid \; \mid^{\alpha} \chi_{1_{F^*}} \St_{\GL_2} ; \lambda') $ is non-unitary.
\end{Theorem}

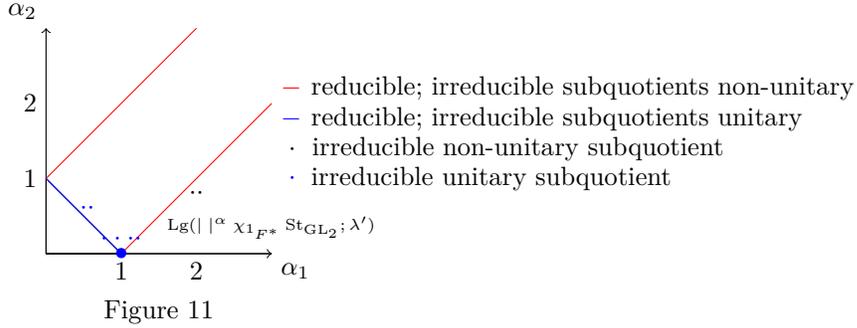
\begin{figure}

\begin{center}
\begin{tikzpicture}
 \draw[->] (0,0) -- (3,0);
\draw (3,0) node[below right] {$\alpha_1$};
\draw[->] (0,0) -- (0,3);
\draw (0,3) node[above left] {$ \alpha_2 $};

\draw  (0,0) -- (1,0) -- (0,1) --cycle;

\draw [blue] (1,0) -- (0,1);
\draw [red] (0,1) -- (2,3);
\draw  [red] (1,0) -- (3,2);
\draw  (0,0) -- (1,0);
\draw  (1,0) -- (3,0);

\draw (1,0) node[below] {$ 1 $};
\draw (0,1) node[left] {$ 1 $};
\draw (2,0) node[below] {$ 2 $};
\draw (0,2) node[left] {$ 2 $};
\draw (1,0) node[blue] {$ \bullet $};

\draw (1,0) node[blue, above] {$ \cdot \cdot \cdot \cdot $};
\draw (0.55, 0.4) node[blue, above] { $ \cdot \cdot$};
\draw (2, 0.6) node[above] { $ \cdot \cdot$};


\draw (1.4, 0.35) node[right] { \begin{tiny} $
  \Lg( \mid \; \mid^{\alpha} \chi_{1_{F^*}} \St_{\GL_2} ; \lambda')                                 
                                 $ \end{tiny} };

\draw (1.5, -0.5) node[below] { Figure 11 };
\draw (3, 2.2) node[red, right] { $ - $};
\draw (3.4, 2.2) node[right] {reducible; irreducible subquotients non-unitary};
\draw (3, 1.8) node[blue, right] { $ - $};
\draw (3.4, 1.8) node[right] {reducible; irreducible subquotients unitary};
\draw (3.1,1.4) node[right] { $ \cdot \; $ irreducible non-unitary subquotient };
\draw (3.1, 1) node[blue, right] { $ \cdot $};
\draw (3.4, 1) node[right] { irreducible unitary subquotient}; 
\end{tikzpicture}
\end{center}

\caption[$ \Lg( \mid \; \mid^{\alpha} \chi_{1_{F^*}} \St_{\GL_2} ; \lambda'), \; \alpha > 0, \; \chi_{1_{F^*}} \in X_{1_{F^*}}
 $]{Let $ \alpha_1, \alpha_2 \geq 0, $ let $ \chi_{1_{F^*}} \in X_{1_{F^*}}. $ Figure 11 shows lines and points of reducibility of
 the representation 
$ \mid \; \mid^{\alpha_1} \chi_{1_{F^*}} \times \mid \; \mid^{\alpha_2} \chi_{1_{F^*}} \rtimes \lambda'. $
 For $ 0 < \alpha \leq 1/2, $
$ \Lg( \mid \; \mid^{\alpha} \chi_{1_{F^*}} \St_{\GL_2} ; \lambda') $ is unitary, for $ \alpha > 1/2 $ it is
non-unitary.}

\label{figure11}

\end{figure}

\begin{Proof}
By Proposition \cite[Prop. 4.8]{CS} $ \chi_{1_{F^*}} \St_{\GL_2} \rtimes \lambda' $ is irreducible. It is unitary. Let $ 0 < \alpha < 1/2. $
The representations  $ \mid \; \mid^{\alpha} \chi_{1_{F^*}} \St_{\GL_2} \rtimes \lambda' $ are irreducible by 
\cite[Th. 4.7]{CS} and equal
to their own Langlands quotient $ \Lg(\mid \; \mid^{\alpha} \chi_{1_{F^*}} \St_{\GL_2} ; \lambda'). $ They form a continuous
1-parameter family of irreducible hermitian representations. Similar as in \ref{Lg11} we realize them on the 
same vector space $ V. $ By Lemma \cite[Lemma 3.6]{CS} $ \Lg(\mid \; \mid^{\alpha} \chi_{1_{F^*}} \St_{\GL_2} ;
 \lambda') $ is unitary for $ 0 < \alpha < 1/2. $  For $ \alpha = 1/2, \; \mid \; \mid^{1/2} \chi_{1_{F^*}} \St_{\GL_2}
\rtimes \lambda' $ reduces
for the first time (\cite[Th. 4.7]{CS}). By \cite{MR0324429} $ \Lg(\mid \; \mid^{1/2} \chi_{1_{F^*}} \St_{\GL_2} ; 
\lambda') $
is unitary.

\medskip

Let $ \alpha > 1/2. $ Representations $ \mid \; \mid^{\alpha} \chi_{1_{F^*}} \St_{\GL_2} \rtimes \lambda' $ are 
irreducible by \cite[Th. 4.7]{CS} and equal to their own Langlands quotient 
$ \Lg(\mid \; \mid^{\alpha} \chi_{1_{F^*}} \St_{\GL_2} ; \lambda'). $
If there existed $ \alpha > 1/2 $ such that $ \Lg(\mid \; \mid^{\alpha} \chi_{1_{F^*}} \St_{\GL_2} ; \lambda') $ is
unitary, then by \cite[Lemma 3.6]{CS} all representations $ \Lg(\mid \; \mid^{\alpha} \chi_{1_{F^*}} \St_{\GL_2} ; \lambda') $
would be unitary for $ \alpha > 1/2, $ in contradiction to \cite[Lemma 3.8]{CS}.
(Figure ~\ref{figure11}, page ~\pageref{figure11})
\end{Proof}

\smallskip

\subsection{Representations induced from $ M_2, $ with cuspidal support in $ M_0, $ not fully-induced}

\bigskip
\smallskip

\begin{center} Irreducible subquotients of $ \chi \rtimes \tau, \; \tau $ tempered non-cuspidal of $ U(3), $
 not fully-induced \end{center}

\medskip

Recall that $ \lambda'(\det) \St_{U(3)} $ is the unique square-integrable subquotient of $ \mid \; \mid 1 
\rtimes \lambda' $ \cite{Ky}. Let $ \chi_{\omega_{E/F}} \in X_{\omega_{E/F}}.  $ Let $ \pi_{1, \chi_{\omega_{E/F}}} $ denote the unique square-integrable irreducible 
subquotient of $ \mid \; \mid^{1/2} \chi_{\omega_{E/F}} \rtimes \lambda'. $ Let $ \chi_{1_{F^*}} \in X_{1_{F^*}}. $
We have $ \chi_{1_{F^*}} \rtimes \lambda' = \sigma_{1, \chi_{1_{F^*}}} \oplus \sigma_{2, \chi_{1_{F^*}}}, $ where
$ \sigma_{1, \chi_{1_{F^*}}} $ and $
\sigma_{2, \chi_{1_{F^*}}} $ are irreducible tempered \cite{Ky}.

\medskip

$ \lambda'(\det) \St_{U(3)}, \pi_{1, \chi_{\omega_{E/F}}}, \sigma_{1, \chi_{1_{F^*}}} $ and $ \sigma_{2, \chi_{1_{F^*}}} $
are all non-cuspidal tempered representations of $ U(3) $ that
are not fully induced \cite{Ky}.

\bigskip

Let $ \chi $ be a unitary character of $ E^*. $

\medskip

The representations $ \chi \rtimes \lambda'(\det) \St_{U(3)}, \;  \chi \rtimes \pi_{1, \chi_{\omega_{E/F}}}, \;
\chi \rtimes \sigma_{1, \chi_{1_{F^*}}} $ and $ \chi \rtimes \sigma_{2, \chi_{1_{F^*}}} $ are tempered, hence unitary.
Hence all their irreducible subquotients are tempered, hence unitary.

\begin{Remark}
 By \cite[Prop. 4.10]{CS} $ \chi \rtimes \lambda'(\det) \St_{U(3)} $ is reducible if and only if $ \chi = 1 $ or
 $ \chi \in
X_{1_{F^*}}. $ By \cite[Prop. 4.12]{CS} $ \chi \rtimes \pi_{1, \chi_{\omega_{E/F}}} $ is reducible if and only if
$ \chi \in X_{1_{F^*}}. $ By \cite[Th. 4.1]{CS} $ \chi \rtimes \sigma_{1, \chi_{1_{F^*}}} $ and 
$ \chi \rtimes \sigma_{2, \chi_{1_{F^*}}} $ are reducible if and only if $ \chi \in X_{1_{F^*}} $ but $ \chi \ncong
\chi_{1_{F^*}}. $

Let $ \chi_{1_{F^*}} \in X_{1_{F^*}}. \; 1 \rtimes \lambda'(\det) \St_{U(3)}, \; \chi_{1_{F^*}} \rtimes \lambda'(\det) \St_{U(3)}, \;
\chi_{1_{F^*}} \rtimes \pi_{1,{\chi_{\omega_{E/F}}}}, \; \chi \rtimes \sigma_{1, \chi_{1_{F^*}}} $ and $ \chi \rtimes
\sigma_{2, \chi_{1_{F^*}}}, $ where $ \chi \in X_{1_{F^*}} $ but $ \chi \ncong \chi_{1_{F^*}} $ have two tempered
subquotients (\cite[Propositions 4.10 and 4.12, Th. 4.1]{CS}).
\end{Remark}
 
\bigskip

\begin{center} $ \Lg(\mid \; \mid^{\alpha} \chi ; \tau), \; \alpha > 0, \; \tau $ tempered non-cuspidal of $ U(3), $ not
fully-induced \end{center}

Recall that $ \lambda'(\det) \St_{U(3)} $ is the unique square-integrable subquotient of $ \mid \; \mid 1 
\rtimes \lambda' $ \cite{Ky}. Let $ \chi_{\omega_{E/F}} \in X_{\omega_{E/F}}.  $ Let $ \pi_{1, \chi_{\omega_{E/F}}} $ denote the unique square-integrable irreducible 
subquotient of $ \mid \; \mid^{1/2} \chi_{\omega_{E/F}} \rtimes \lambda'. $ Let $ \chi_{1_{F^*}} \in X_{1_{F^*}}. $
We have $ \chi_{1_{F^*}} \rtimes \lambda' = \sigma_{1, \chi_{1_{F^*}}} + \sigma_{2, \chi_{1_{F^*}}}, $ where
$ \sigma_{1, \chi_{1_{F^*}}} $ and $
\sigma_{2, \chi_{1_{F^*}}} $ are irreducible tempered \cite{Ky}.

\medskip

$ \lambda'(\det) \St_{U(3)}, \pi_{1, \chi_{\omega_{E/F}}}, \sigma_{1, \chi_{1_{F^*}}} $ and $ \sigma_{2, \chi_{1_{F^*}}} $
are all non-cuspidal tempered representations of $ U(3) $ that
are not fully induced \cite{Ky}.

\bigskip

\begin{Theorem}
Let $  \chi \notin X_{N_{E/F}}(E^*). $

Let $ \alpha > 0. \;  \Lg(\mid \; \mid^{\alpha} \chi ; \lambda'(\det) \St_{U(3)}), \Lg( \mid \; \mid^{\alpha} \chi, 
\pi_{1, \chi_{\omega_{E/F}}}),
\Lg( \mid \; \mid^{\alpha} \chi; \sigma_{1, \chi_{1_{F^*}}}) $ and $ \Lg( \mid \; \mid^{\alpha} \chi ;
\sigma_{2, \chi_{1_{F^*}}}) $ are non-unitary.

\end{Theorem}

\begin{Proof} The representations $ \Lg(\mid \; \mid^{\alpha} \chi ; \lambda'(\det) \St_{U(3)}), 
\Lg( \mid \; \mid^{\alpha} \chi, \pi_{1, \chi_{1_{F^*}}}),
\Lg( \mid \; \mid^{\alpha} \chi; \sigma_{1, \chi_{1_{F^*}}}) $ and $ \Lg( \mid \; \mid^{\alpha} \chi ; 
\sigma_{2, \chi_{1_{F^*}}}) $ are not hermitian, hence not
unitary.
\end{Proof}

\bigskip 

Let $ \chi \in X_{N_{E/F}}(E^*) = 1 \cup X_{\omega_{E/F}} \cup X_{1_{F^*}}. $

Let $ \chi_{\omega_{E/F}}, \; \chi_{\omega_{E/F}}' \in X_{\omega_{E/F}}, $ such that $ \chi_{\omega_{E/F}}' \neq 
\chi_{\omega_{E/F}}, $ and let
$ \chi_{1_{F^*}}, \; \chi_{1_{F^*}}' \in X_{1_{F^*}}, $ such that $ \chi_{1_{F^*}}' \neq \chi_{1_{F^*}}. $ 

\bigskip

\begin{Theorem}

\label{alphachitau}

1. Let $ \alpha > 1. \; \Lg(\mid \; \mid^{\alpha} 1 ; \lambda'(\det) \St_{U(3)} ) $ is non-unitary.

2. Let $ 0 < \alpha \leq 1. \Lg(\mid \; \mid^{\alpha} 1 ; \pi_{1, \chi_{\omega_{E/F}}}), \; \Lg(\mid \; \mid^{\alpha} 1 ; 
\sigma_{1, \chi_{1_{F^*}}} ), \;
 \Lg(\mid \; \mid^{\alpha} 1 ; \sigma_{2, \chi_{1_{F^*}}} ), \Lg(\mid \; \mid^{\alpha} \chi_{1_{F^*}} ; 
\sigma_{1, \chi_{1_{F^*}}} ) $
 and $
 \Lg(\mid \; \mid^{\alpha} \chi_{1_{F^*}} ; \sigma_{2, \chi_{1_{F^*}}} )
$ are unitary.

3. Let $ \alpha > 1.  \Lg(\mid \; \mid^{\alpha} 1 ; \pi_{1,\chi_{\omega_{E/F}}}), \; \Lg(\mid \; \mid^{\alpha} 1 ; 
\sigma_{1, \chi_{1_{F^*}}} ), $
$ \Lg(\mid \; \mid^{\alpha} 1 ; \sigma_{2, \chi_{1_{F^*}}} ), \Lg(\mid \; \mid^{\alpha} \chi_{1_{F^*}} ; 
\sigma_{1, \chi_{1_{F^*}}} ) $
 and $
 \Lg(\mid \; \mid^{\alpha} \chi_{1_{F^*}} ; \sigma_{2, \chi_{1_{F^*}}} )
$ are non-unitary.

4. Let $ 0 < \alpha \leq 1/2. \; \Lg(\mid \; \mid^{\alpha} \chi_{\omega_{E/F}} ; \lambda'(\det) \St_{U(3)}),
\Lg(\mid \; \mid^{\alpha} \chi_{\omega_{E/F}} ; \pi_{1, \chi_{\omega_{E/F}}}), \; \Lg(\mid \; \mid^{\alpha}
 \chi_{\omega_{E/F}}' ; \pi_{1, \chi_{\omega_{E/F}}}),
 \; \Lg(\mid \; \mid^{\alpha} \chi_{\omega_{E/F}} ; \sigma_{1, \chi_{1_{F^*}}} ) $
 and $
 \Lg(\mid \; \mid^{\alpha} \chi_{\omega_{E/F}} ; \sigma_{2, \chi_{1_{F^*}}} ) $ are unitary.

5. Let $ \alpha > 1/2. \; \Lg(\mid \; \mid^{\alpha} \chi_{\omega_{E/F}} ; \lambda'(\det) \St_{U(3)}),
\Lg(\mid \; \mid^{\alpha} \chi_{\omega_{E/F}} ; \pi_{1, \chi_{\omega_{E/F}}}), \; 
\Lg(\mid \; \mid^{\alpha} \chi_{\omega_{E/F}}' ; \pi_{1, \chi_{\omega_{E/F}}}),
    \;
 \Lg(\mid \; \mid^{\alpha} \chi_{\omega_{E/F}} ; \sigma_{1, \chi_{1_{F^*}}} ) $
 and $
 \Lg(\mid \; \mid^{\alpha} \chi_{\omega_{E/F}} ; \sigma_{2, \chi_{1_{F^*}}} ) $ are non-unitary.

6. Let $ \alpha > 0. \; \Lg(\mid \; \mid^{\alpha} \chi_{1_{F^*}} ; \lambda'(\det) \St_{U(3)}), \;
\Lg(\mid \; \mid^{\alpha} \chi_{1_{F^*}} ; \pi_{1, \chi_{\omega_{E/F}}}), \; \Lg(\mid \; \mid^{\alpha} \chi_{1_{F^*}}' ;
 \sigma_{1, \chi_{1_{F^*}}}), $ and 
$  \Lg(\mid \; \mid^{\alpha} \chi_{1_{F^*}}' ; \sigma_{2, \chi_{1_{F^*}}}), $ are non-unitary.
\end{Theorem}

\begin{Proof}
1. Let $ 1 < \alpha < 2. $ By \cite[Th. 4.9]{CS} $ \mid \; \mid^{\alpha} 1 \rtimes \lambda'(\det) \St_{U(3)} $ is 
irreducible
 and equal to its
Langlands-quotient $ \Lg(\mid \; \mid^{\alpha} 1 ; \lambda'(\det) \St_{U(3)} ). $ By \cite{MR656064} the
representation $ \mid \; \mid^2 1 \rtimes \lambda'(\det) \St_{U(3)} $ is reducible. By the
same author $ \Lg(\mid \; \mid^2 1; \lambda'(\det) \St_{U(3)}) $ is non-unitary.
 By \cite{MR0324429} and \cite[Lemma 3.6]{CS}
$ \Lg(\mid \; \mid^{\alpha} 1 ; \lambda'(\det) \St_{U(3)}) $ is non-unitary for $ 1 < \alpha < 2. $

\smallskip

Let $ \alpha > 2. \;  \mid \; \mid^{\alpha} 1 \rtimes \lambda'(\det) \St_{U(3)} $ is irreducible by
\cite[Th. 4.9]{CS}
and equal
 to its Langlands-quotient $ \Lg(\mid \; \mid^{\alpha} 1 ; \lambda'(\det) \St_{U(3)} ). $ By \cite[Lemmas 3.6 and 3.8]{CS} and 
$ \Lg(\mid \; \mid^{\alpha} 1 ; \lambda'(\det) \St_{U(3)} ) $ is non-unitary for $ \alpha > 2. $

\smallskip

Let $ 0 < \alpha \leq 1. $ We have no proof that $ \Lg(\mid \; \mid^{\alpha} 1 ; \lambda'(\det) \St_{U(3)} ) $ is non-unitary.
(Figure ~\ref{figure12}, page ~\pageref{figure12}).
\medskip

2. The representations $ 1 \rtimes \pi_{1, \chi_{\omega_{E/F}}} $ are irreducible by \cite[Prop. 4.12]{CS},
representations $
 1 \rtimes \sigma_{1, \chi_{1_{F^*}}}, 1 \rtimes 
\sigma_{2, \chi_{1_{F^*}}}, \chi_{1_{F^*}} \rtimes \sigma_{1, \chi_{1_{F^*}}} $ and
$ \chi_{1_{F^*}} \rtimes \sigma_{2, \chi_{1_{F^*}}} $ are irreducible by \cite[Th. 4.1]{CS}. All representations are unitary.
For $ 0 < \alpha < 1, $ representations $ \mid \; \mid^{\alpha} 1
 \rtimes \pi_{1, \chi_{\omega_{E/F}}} $ are irreducible by \cite[Th. 4.11]{CS}, representations $
 \mid \; \mid^{\alpha} 1 \rtimes \sigma_{1, \chi_{1_{F^*}}}, \; 
\mid \; \mid^{\alpha} 1 \rtimes \sigma_{2, \chi_{1_{F^*}}},
 \mid \; \mid^{\alpha} \chi_{1_{F^*}} \rtimes \sigma_{1, \chi_{1_{F^*}}}  $ and $ \mid \; \mid^{\alpha} 
\chi_{1_{F^*}} \rtimes \sigma_{2, \chi_{1_{F^*}}} $ are
 irreducible by \cite[Th. 4.13]{CS}. The representations are equal to their own Langlands-quotients 
$ \Lg(\mid \; \mid^{\alpha} 1 ; \pi_{1, \chi_{\omega_{E/F}}}), \; 
\Lg(\mid \; \mid^{\alpha} 1 ; \sigma_{1, \chi_{1_{F^*}}} ), \; \Lg(\mid \; \mid^{\alpha} 1 ; \sigma_{2, \chi_{1_{F^*}}} ), \;
 \Lg(\mid \; \mid^{\alpha} \chi_{1_{F^*}} ; \sigma_{1, \chi_{1_{F^*}}} ) $ and $ 
\Lg(\mid \; \mid^{\alpha} \chi_{1_{F^*}} ; \sigma_{2, \chi_{1_{F^*}}} ), $
 respectively.
By Lemma \cite[Lemma 3.6]{CS} these Langlands-quotients are unitary. For $ \alpha = 1, \; \mid \; \mid 1 \rtimes
 \pi_{1, \chi_{\omega_{E/F}}}, \; 
\mid \; \mid 1 \rtimes \sigma_{1, \chi_{1_{F^*}}}, \; \mid \; \mid 1 \rtimes \sigma_{2, \chi_{1_{F^*}}}, 
\mid \; \mid \chi_{1_{F^*}} \rtimes \sigma_{1, \chi_{1_{F^*}}}  $
and $ \mid \; \mid \chi_{1_{F^*}} \rtimes \sigma_{2, \chi_{1_{F^*}}}  $
reduce for
the first time (\cite[Th. 4.11 and 4.13]{CS}. By \cite{MR0324429} 
$ \Lg( \mid \; \mid 1; \pi_{1, \chi_{\omega_{E/F}}}), \; 
\Lg(\mid \; \mid 1 ; \sigma_{1, \chi_{1_{F^*}}} ), \; \Lg(\mid \; \mid 1 ; \sigma_{2, \chi_{1_{F^*}}} ),
\;
\Lg(\mid \; \mid^{\alpha} \chi_{1_{F^*}} ; \sigma_{1, \chi_{1_{F^*}}} ) $ and $ \Lg(\mid \; \mid^{\alpha} \chi_{1_{F^*}} ; 
\sigma_{2, \chi_{1_{F^*}}} ) $ 
are unitary.

\medskip

3. For $ \alpha > 1, $ representations $  \mid \; \mid^{\alpha} 1 \rtimes \pi_{1, \chi_{\omega_{E/F}}},  \; 
\mid \; \mid^{\alpha} 1 \rtimes
 \sigma_{1, \chi_{1_{F^*}}}, \; \mid \; \mid^{\alpha} 1 \rtimes \sigma_{2, \chi_{1_{F^*}}}, \mid \; \mid^{\alpha} 
\chi_{1_{F^*}} \rtimes \sigma_{1, \chi_{1_{F^*}}}  $
 and $ \mid \; \mid^{\alpha} \chi_{1_{F^*}} \rtimes \sigma_{2, \chi_{1_{F^*}}} $ are irreducible
(\cite[Th. 4.11 and 4.13]{CS} and equal to their Langlands-quotients
$ \Lg( \mid \; \mid^{\alpha} 1 ; \pi_{1, \chi_{\omega_{E/F}}}), \; \Lg(\mid \; \mid^{\alpha} 1 ;
 \sigma_{1, \chi_{1_{F^*}}} ), \Lg(\mid \; \mid^{\alpha} 1 ; 
\sigma_{2, \chi_{1_{F^*}}} ), \Lg( \mid \; \mid^{\alpha} \chi_{1_{F^*}} ; \sigma_{1, \chi_{1_{F^*}}}) $ and 
$ \Lg(\mid \; \mid^{\alpha} \chi_{1_{F^*}}, \sigma_{2, \chi_{1_{F^*}}}), $  respectively.
By \cite[Lemmas 3.6 and 3.8]{CS} these Langlands-quotients are non-unitary.

\medskip

4. Representations $ \chi_{\omega_{E/F}} \rtimes \lambda'(\det) \St_{U(3)}, \; \chi_{\omega_{E/F}} \rtimes
\pi_{1, \chi_{\omega_{E/F}}}, \;
\chi_{\omega_{E/F}}' \rtimes \pi_{1, \chi_{\omega_{E/F}}}, \;
\chi_{\omega_{E/F}} \rtimes \sigma_{1, \chi_{1_{F^*}}} $ and $ \chi_{\omega_{E/F}} \rtimes \sigma_{2, \chi_{1_{F^*}}} $
are irreducible (\cite[Prop. 4.10, 4.12 and Th. 4.1]{CS}) and unitary.
For $ 0 < \alpha < 1/2, $ representations $ \mid \; \mid^{\alpha} \chi_{\omega_{E/F}} \rtimes \lambda'(\det) \St_{U(3)},
\mid \; \mid^{\alpha} \chi_{\omega_{E/F}} \rtimes \pi_{1, \chi_{\omega_{E/F}}}, \; \mid \; \mid^{\alpha} 
\chi_{\omega_{E/F}}' \rtimes \pi_{1, \chi_{\omega_{E/F}}},\;
\mid \; \mid^{\alpha} \chi_{\omega_{E/F}} \rtimes \sigma_{1, \chi_{1_{F^*}}} $
 and $
 \mid \; \mid^{\alpha} \chi_{\omega_{E/F}} \rtimes \sigma_{2, \chi_{1_{F^*}}} $ are irreducible 
(\cite[Th. 4.9, 4.11 and 4.13]{CS}) and equal 
to their Langlands-quotients
 $ \Lg(\mid \; \mid^{\alpha} \chi_{\omega_{E/F}} ; \lambda'(\det) \St_{U(3)}),
\Lg(\mid \; \mid^{\alpha} \chi_{\omega_{E/F}} ; \pi_{1, \chi_{\omega_{E/F}}}), \; \Lg(\mid \; \mid^{\alpha} 
\chi_{\omega_{E/F}}' ; \pi_{1, \chi_{\omega_{E/F}}}),
 \;\Lg(\mid \; \mid^{\alpha} \chi_{\omega_{E/F}} ; \sigma_{1, \chi_{1_{F^*}}} ) $
 and $
 \Lg(\mid \; \mid^{\alpha} \chi_{\omega_{E/F}} ; \sigma_{2, \chi_{1_{F^*}}} ), $ respectively. By \cite[Lemma 3.6]{CS} these
Langlands-quotients are
unitary.  For $ \alpha = 1/2, \; \mid \; \mid^{1/2} \chi_{\omega_{E/F}} \rtimes \lambda'(\det) \St_{U(3)},
\mid \; \mid^{1/2} \chi_{\omega_{E/F}} \rtimes \pi_{1, \chi_{\omega_{E/F}}}, \; \mid \; \mid^{1/2} \chi_{\omega_{E/F}}'
 \rtimes \pi_{1, \chi_{\omega_{E/F}}}, \; 
\mid \; \mid^{1/2} \chi_{\omega_{E/F}} \rtimes \sigma_{1, \chi_{1_{F^*}}} $
 and $
 \mid \; \mid^{1/2} \chi_{\omega_{E/F}} \rtimes \sigma_{2, \chi_{1_{F^*}}} $ reduce for the first time (\cite[Th. 4.9, 4.11 and
4.13]{CS}. By \cite{MR0324429} $
 \Lg(\mid \; \mid^{1/2} \chi_{\omega_{E/F}} ; \lambda'(\det) \St_{U(3)}),
\Lg(\mid \; \mid^{1/2} \chi_{\omega_{E/F}} ; \pi_{1, \chi_{\omega_{E/F}}}), \; \Lg(\mid \; \mid^{1/2} 
\chi_{\omega_{E/F}}' ; \pi_{1, \chi_{\omega_{E/F}}}), 
\;\Lg(\mid \; \mid^{1/2} \chi_{\omega_{E/F}} ; \sigma_{1, \chi_{1_{F^*}}} ) $
 and $
 \Lg(\mid \; \mid^{1/2} \chi_{\omega_{E/F}} ; \sigma_{2, \chi_{1_{F^*}}} ) $ are unitary. 

\medskip

5. For $ \alpha > 1/2, $ representations $ \mid \; \mid^{\alpha} \chi_{\omega_{E/F}} \rtimes \lambda'(\det) \St_{U(3)},
\mid \; \mid^{\alpha} \chi_{\omega_{E/F}} \rtimes \pi_{1, \chi_{\omega_{E/F}}}, \; \mid \; \mid^{\alpha} 
\chi_{\omega_{E/F}}' \rtimes \pi_{1, \chi_{\omega_{E/F}}}, 
\; \mid \; \mid^{\alpha} \chi_{\omega_{E/F}} \rtimes \sigma_{1, \chi_{1_{F^*}}} $
 and $
 \mid \; \mid^{\alpha} \chi_{\omega_{E/F}} \rtimes \sigma_{2, \chi_{1_{F^*}}} $ are irreducible (
\cite[Th. 4.9, 4.11 and 4.13]{CS})
and equal to their Langlands-quotients
 $ \Lg(\mid \; \mid^{\alpha} \chi_{\omega_{E/F}} ; \lambda'(\det) \St_{U(3)}),
\Lg(\mid \; \mid^{\alpha} \chi_{\omega_{E/F}} ; \pi_{1, \chi_{\omega_{E/F}}}), \; \Lg(\mid \; \mid^{\alpha} 
\chi_{\omega_{E/F}}' ; \pi_{1, \chi_{\omega_{E/F}}}),
 \;\Lg(\mid \; \mid^{\alpha} \chi_{\omega_{E/F}} ; \sigma_{1, \chi_{1_{F^*}}} ) $
 and $
 \Lg(\mid \; \mid^{\alpha} \chi_{\omega_{E/F}} ; \sigma_{2, \chi_{1_{F^*}}} ), $ respectively. By 
\cite[Lemmas 3.6 and 3.8]{CS} these Langlands-quotients are non-unitary.

\medskip

6. Let $ \alpha > 0. $ The representations $ \mid \; \mid^{\alpha} \chi_{1_{F^*}} \rtimes \lambda'(\det) \St_{U(3)}, \;
 \mid \; \mid^{\alpha} \chi_{1_{F^*}} \rtimes \pi_{1, \chi_{\omega_{E/F}}}, \; \mid \; \mid^{\alpha} \chi_{1_{F^*}}' 
\rtimes \sigma_{1, \chi_{1_{F^*}}} $ and
$ \mid \; \mid^{\alpha} \chi_{1_{F^*}}' \rtimes \sigma_{2, \chi_{1_{F^*}}} $ are irreducible
(\cite[Th. 4.9, 4.11 and 4.13]{CS}) and equal to their Langlands-quotients
$ \Lg(\mid \; \mid^{\alpha} \chi_{1_{F^*}} ; \lambda'(\det) \St_{U(3)}), \; 
\Lg(\mid \; \mid^{\alpha} \chi_{1_{F^*}} ; \pi_{1, \chi_{\omega_{E/F}}}), \; \Lg(\mid \; \mid^{\alpha} \chi_{1_{F^*}}' ;
 \sigma_{1, \chi_{1_{F^*}}}) $ and
 $ \Lg(\mid \; \mid^{\alpha} \chi_{1_{F^*}}' ; \sigma_{2, \chi_{1_{F^*}}}). $ By \cite[Lemmas 3.6 and 3.8]{CS} these 
Langlands-quotients are non-unitary.
\end{Proof}

\smallskip

\begin{figure}
 
\begin{center}
\begin{tikzpicture}
 \draw[->] (0,0) -- (3,0);
\draw (3,0) node[below right] {$\alpha_1$};
\draw[->] (0,0) -- (0,3);
\draw (0,3) node[above left] {$ \alpha_2 $};

\draw  (0,0) -- (1,0) -- (0,1) --cycle;
\draw  (1,0) -- (0,1);

\draw  (0,1) -- (3,1);
\draw  (0,1) -- (2,3);
\draw  (1,0) -- (1,1);
\draw [red] (1,1) -- (1,3);
\draw  (1,0) -- (3,2);

\draw (1,0) node[below] {$ 1 $};
\draw (0,1) node[left] {$ 1 $};
\draw (2,0) node[below] {$ 2 $};
\draw (0,2) node[left] {$ 2 $};
\draw (1,1) node {$ \bullet $};



\draw (1, 1.35) node[right] { \begin{tiny} $ {\Lg(\mid \; \mid^{\alpha} 1 ; \lambda'(\det) \St_{U(3)})}  $ \end{tiny}};

\draw (1.5, -0.5) node[below] { Figure 12 };
\draw (3, 2.2) node[red, right] { $ - $};
\draw (3.4, 2.2) node[right] {reducible; irreducible subquotients non-unitary};
\draw (3.1,1.8) node[right] { $ \cdot \; $ irreducible non-unitary subquotient };
\end{tikzpicture}
\end{center}

\caption[$ \Lg(\mid \; \mid^{\alpha} 1 \rtimes \lambda'(\det) \St_{U(3)}), \; \alpha > 0 $]{Let $ \alpha_1, \alpha_2 \geq 0.
 $ Figure 12 shows lines and points of reducibility of the representation 
$ \mid \; \mid^{\alpha_1} 1 \times \mid \; \mid^{\alpha_2} 1 \rtimes \lambda'. $
For $ \alpha > 1, \; \Lg(\mid \; \mid^{\alpha} 1 \rtimes \lambda'(\det) \St_{U(3)}) $ is non-unitary.}

\label{figure12}

\end{figure}

\medskip

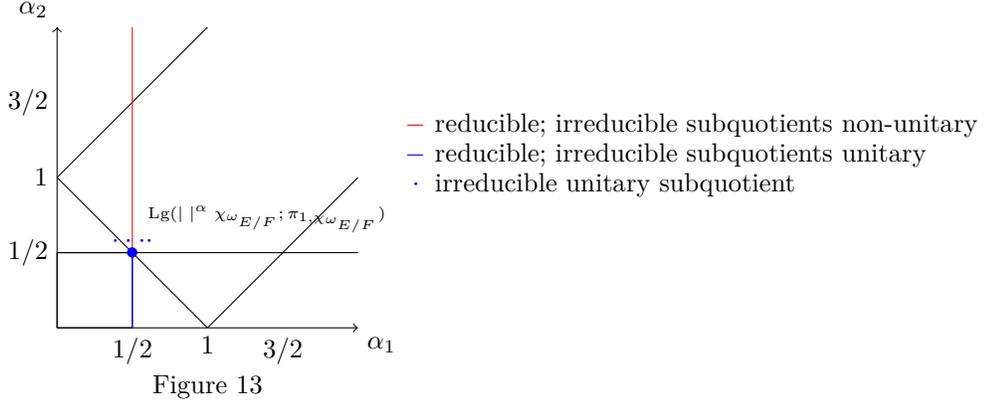
\begin{figure}
 
\begin{center}
\begin{tikzpicture}
 \draw[->] (0,0) -- (4,0);
\draw (4,0) node[below right] {$\alpha_1$};
\draw[->] (0,0) -- (0,4);
\draw (0,4) node[above left] {$ \alpha_2 $};

\draw  (0,0) -- (1,0) -- (1,1) -- (0,1) -- cycle;

\draw  (1 ,1) -- (4,1);
\draw  (0,2) -- (2,4);
\draw [red] (1,1) -- (1,4);
\draw  (2,0) -- (4,2);
\draw  (2,0) -- (0,2);
\draw [blue] (1, 0) -- (1, 1);
\draw (0, 1) -- (1, 1);

\draw (1,0) node[below] {$ 1/2 $};
\draw (0, 1) node[left] {$ 1/2 $};
\draw (2,0) node[below] {$ 1 $};
\draw (0,2) node[left] {$ 1 $};
\draw (3,0) node[below] { $ 3/2 $};
\draw (0, 3) node[left] { $ 3/2 $};
\draw (1, 1) node[blue] {$ \bullet $};
\draw (1,0.95) node[blue, above] { $ \cdot \cdot \cdot \cdot $};

\draw (1, 1.45) node[right] { \begin{tiny} $ \Lg(\mid \; \mid^{\alpha} \chi_{\omega_{E/F}} ; 
\pi_{1, \chi_{\omega_{E/F}}}) $ \end{tiny} };

\draw (2, -0.5) node[below] { Figure 13 };
\draw (4.5, 2.7) node[red, right] { $ - $};
\draw (4.9, 2.7) node[right] {reducible; irreducible subquotients non-unitary};
\draw (4.5, 2.3) node[blue, right] { $ - $};
\draw (4.9, 2.3) node[right] {reducible; irreducible subquotients unitary};
\draw (4.6, 1.9) node[blue, right] { $ \cdot $};
\draw (4.9, 1.9) node[right] { irreducible unitary subquotient}; 


\end{tikzpicture}
\end{center}

\caption[$ \Lg(\mid \; \mid^{\alpha} \chi_{\omega_{E/F}} ; \pi_{1, \chi_{\omega_{E/F}}}), \; \alpha > 0, \;
\chi_{\omega_{E/F}} \in X_{\omega_{E/F}} $]{Let $ \alpha_1, \alpha_2 \geq 0, $ let $ \chi_{\omega_{E/F}} 
\in X_{\omega_{E/F}}. $
Figure 13 shows lines and points of reducibility of the representation 
$ \mid \; \mid^{\alpha_1} \chi_{\omega_{E/F}} \times \mid \; \mid^{\alpha_2} \chi_{\omega_{E/F}} \rtimes \lambda'. $
Let $ \pi_{1, \chi_{\omega_{E/F}}} $ be the
unique square-integrable subquotient of $ \mid \; \mid^{1/2} \chi_{\omega_{E/F}} \rtimes \lambda'. $
$ \Lg(\mid \; \mid^{\alpha} \chi_{\omega_{E/F}} ; \pi_{1, \chi_{\omega_{E/F}}})
 $ is unitary for $ 0 < \alpha \leq 1/2. $ It is non-unitary for $ \alpha > 1/2. $} 

\label{figure13}

\end{figure}

\medskip

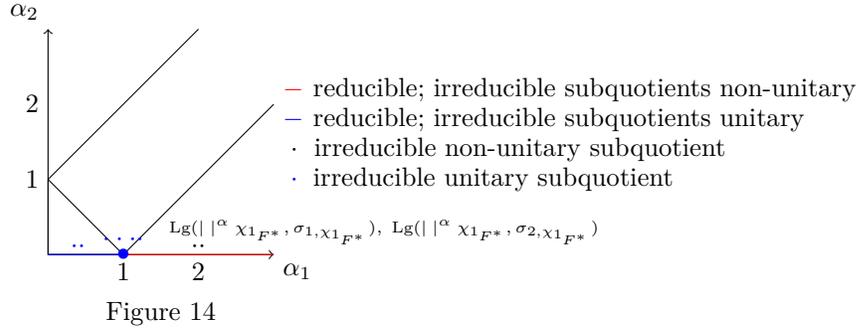
\begin{figure}
 
\begin{center}
\begin{tikzpicture}
 \draw[->] (0,0) -- (3,0);
\draw (3,0) node[below right] {$\alpha_1$};
\draw[->] (0,0) -- (0,3);
\draw (0,3) node[above left] {$ \alpha_2 $};

\draw  (0,0) -- (1,0) -- (0,1) --cycle;

\draw  (0,1) -- (2,3);
\draw  (1,0) -- (3,2);
\draw [blue] (0,0) -- (1,0);
\draw [red] (1,0) -- (3,0);

\draw (1,0) node[below] {$ 1 $};
\draw (0,1) node[left] {$ 1 $};
\draw (2,0) node[below] {$ 2 $};
\draw (0,2) node[left] {$ 2 $};
\draw (1,0) node[blue] {$ \bullet $};

\draw (1,0) node[blue, above] {$ \cdot \cdot \cdot \cdot $};
\draw (0.4, -0.1) node[blue, above] {$ \cdot \cdot $};
\draw (2, -0.1) node[above] { $ \cdot \cdot $};


\draw (1.4, 0.3) node[right] { \begin{tiny} $
  \Lg( \mid \; \mid^{\alpha} \chi_{1_{F^*}}, \sigma_{1, \chi_{1_{F^*}}}), \; \Lg(\mid \; \mid^{\alpha} \chi_{1_{F^*}}, 
\sigma_{2, \chi_{1_{F^*}}})                                 
                                 $ \end{tiny} };

\draw (1.5, -0.5) node[below] { Figure 14 };
\draw (3, 2.2) node[red, right] { $ - $};
\draw (3.4, 2.2) node[right] {reducible; irreducible subquotients non-unitary};
\draw (3, 1.8) node[blue, right] { $ - $};
\draw (3.4, 1.8) node[right] {reducible; irreducible subquotients unitary};
\draw (3.1,1.4) node[right] { $ \cdot \; $ irreducible non-unitary subquotient };
\draw (3.1, 1) node[blue, right] { $ \cdot $};
\draw (3.4, 1) node[right] { irreducible unitary subquotient}; 
\end{tikzpicture}
\end{center}

\caption[$ \Lg( \mid \; \mid^{\alpha} \chi_{1_{F^*}}, \sigma_{1, \chi_{1_{F^*}}}) $ and $ \Lg(\mid \; \mid^{\alpha}
 \chi_{1_{F^*}}, \sigma_{2, \chi_{1_{F^*}}}), \; \alpha > 0 $]{Let $ \alpha_1, \alpha_2 \geq 0, $ let $ \chi_{1_{F^*}}
\in X_{1_{F^*}}. $ Figure 14 shows lines and points of reducibility of
 the representation $ \mid \; \mid^{\alpha_1} \chi_{1_{F^*}} \times \mid \; \mid^{\alpha_2} \chi_{1_{F^*}} \rtimes \lambda'. $ Let $ \alpha > 0. $
$ \Lg( \mid \; \mid^{\alpha} \chi_{1_{F^*}}, \sigma_{1, \chi_{1_{F^*}}}) $ and $ \Lg(\mid \; \mid^{\alpha}
 \chi_{1_{F^*}}, \sigma_{2, \chi_{1_{F^*}}}) $
are unitary for $ 0 < \alpha \leq 1 $ and non-unitary for $ \alpha > 1. $}

\label{figure14}

\end{figure}

\begin{figure}
 
\begin{center}
\begin{tikzpicture}
 \draw[->] (0,0) -- (4,0);
\draw (4,0) node[below right] {$\alpha_1$};
\draw[->] (0,0) -- (0,4);
\draw (0,4) node[above left] {$ \alpha_2 $};

\draw (0,0) -- (2,0) -- (2,1) -- (0,1) -- cycle;

\draw [red] (2 ,1) -- (4,1);
\draw [red] (2,1) -- (2,4);
\draw [blue] (2, 0) -- (2, 1);
\draw [blue] (0, 1) -- (2, 1);

\draw (1,0) node[below] {$ 1/2 $};
\draw (0, 1) node[left] {$ 1/2 $};
\draw (2,0) node[below] {$ 1 $};
\draw (0,2) node[left] {$ 1 $};
\draw (2, 1) node[blue] {$ \bullet $};
\draw (2,0.95) node[blue, above] { $ \cdot \cdot \cdot \cdot $};

\draw (1, 0.9) node[blue, above] { $ \cdot \cdot $};
\draw (2.15, 0.25) node[blue, above] { $ \cdot \cdot $};
\draw (3, 0.9) node[above]  {$ \cdot \cdot $};
\draw (2.15, 2.2) node[above] { $ \cdot \cdot $};

\draw (2, -0.5) node[below] { Figure 15 };

\draw (2, 3) node[right] {\begin{tiny} $ \Lg( \mid \; \mid^{\alpha} \chi_{\omega_{E/F}} ; \lambda'(\det) \St_{U(3)}) $ 
\end{tiny}};
\draw (2.4, 0.7) node[right] { \begin{tiny} $ \Lg( \mid \; \mid^{\alpha} 1 ; \pi_{1, \chi_{\omega_{E/F}}}) $ \end{tiny}};

\draw (4.5, 2.2) node[red, right] { $ - $};
\draw (4.9, 2.2) node[right] {reducible; irreducible subquotients non-unitary};
\draw (4.5, 1.8) node[blue, right] { $ - $};
\draw (4.9, 1.8) node[right] {reducible; irreducible subquotients unitary};
\draw (4.6,1.4) node[right] { $ \cdot \; $ irreducible non-unitary subquotient };
\draw (4.6, 1) node[blue, right] { $ \cdot $};
\draw (4.9, 1) node[right] { irreducible unitary subquotient}; 


\end{tikzpicture}
\end{center}

\caption[$ \Lg( \mid \; \mid^{\alpha} 1 ; \pi_{1, \chi_{\omega_{E/F}}}) $ and $ \Lg( \mid \; \mid^{\alpha} \chi_{\omega_{E/F}} ; 
\lambda'(\det) \St_{U(3)}), \; \alpha > 0 $]{Let $ \alpha_1, \alpha_2 \geq 0, $ let $ \chi_{\omega_{E/F}} \in X_{\omega_{E/F}}.
 $ Figure 15 shows lines and points 
of reducibility of the representation $ \mid \; \mid^{\alpha_1} 1 \times \mid \; \mid^{\alpha_2} \chi_{\omega_{E/F}} 
\rtimes 
\lambda'. \;  \Lg( \mid \; \mid^{\alpha} 1 ; \pi_{1, \chi_{\omega_{E/F}}}) $ is unitary for $ 0 < \alpha \leq 1. $ It is
non-unitary for $ \alpha > 1. \; \Lg( \mid \; \mid^{\alpha} \chi_{\omega_{E/F}} ; \lambda'(\det) \St_{U(3)}) $ is unitary for
$ 0 < \alpha \leq 1/2. $ It is non-unitary for $ \alpha > 1/2. $}

\label{figure15}

\end{figure}
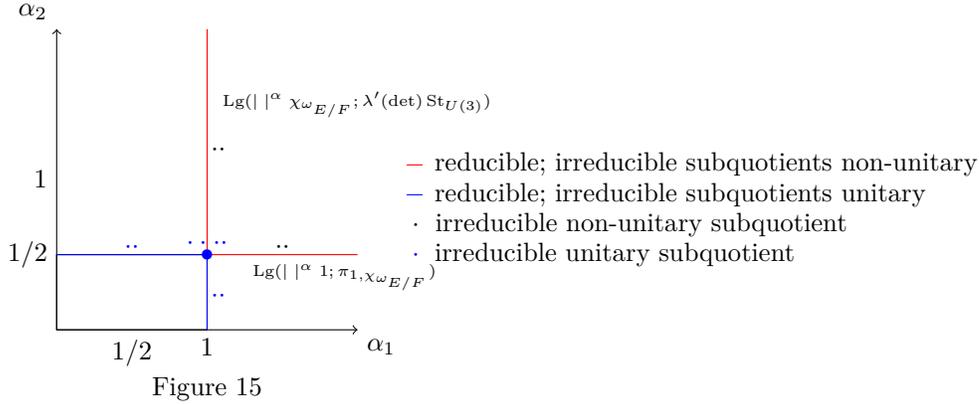

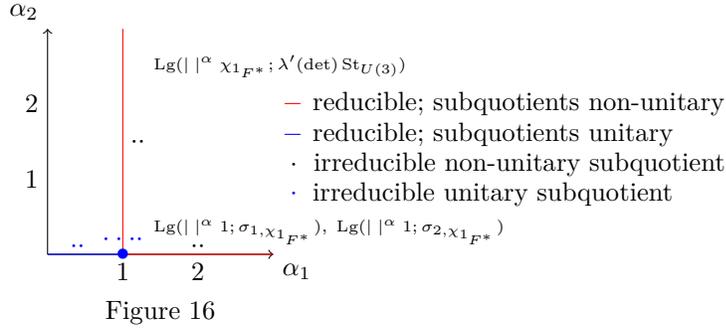
\begin{figure}

\begin{center}
\begin{tikzpicture}
 \draw[->] (0,0) -- (3,0);
\draw (3,0) node[below right] {$\alpha_1$};
\draw[->] (0,0) -- (0,3);
\draw (0,3) node[above left] {$ \alpha_2 $};

\draw [blue] (0,0) -- (1,0);

\draw [red] (1,0) -- (3,0);
\draw [red] (1,0) -- (1,3);

\draw (1,0) node[below] {$ 1 $};
\draw (0,1) node[left] {$ 1 $};
\draw (2,0) node[below] {$ 2 $};
\draw (0,2) node[left] {$ 2 $};
\draw (1,0) node[blue] {$ \bullet $};

\draw (1,0) node[blue, above] {$ \cdot \cdot \cdot \cdot $};
\draw (0.4, -0.1) node[blue, above] {$ \cdot \cdot $};
\draw (2, -0.1) node[above] { $ \cdot \cdot $};
\draw (1.2, 1.5) node { $ \cdot \cdot $};

\draw(1.2, 2.5) node[right] { \begin{tiny} $ \Lg(\mid \; \mid^{\alpha} \chi_{1_{F^*}}; \lambda'(\det) \St_{U(3)})
                               $ \end{tiny}};

\draw (1.2,0.3) node[right] { \begin{tiny} $ \Lg(\mid \; \mid^{\alpha} 1 ; \sigma_{1, \chi_{1_{F^*}}}), \; 
\Lg( \mid \; \mid^{\alpha} 1 ;
\sigma_{2, \chi_{1_{F^*}}}) $ \end{tiny}};

\draw (1.5, -0.5) node[below] { Figure 16 };
\draw (3, 2.0) node[red, right] { $ - $};
\draw (3.4, 2.0) node[right] {reducible; subquotients non-unitary};
\draw (3, 1.6) node[blue, right] { $ - $};
\draw (3.4, 1.6) node[right] {reducible; subquotients unitary};
\draw (3.1,1.2) node[right] { $ \cdot \; $ irreducible non-unitary subquotient };
\draw (3.1, 0.8) node[blue, right] { $ \cdot $};
\draw (3.4, 0.8) node[right] { irreducible unitary subquotient}; 
\end{tikzpicture}
\end{center}

\caption[$ \Lg(\mid \; \mid^{\alpha} 1 ; \sigma_{1, \chi_{1_{F^*}}}), \; 
\Lg( \mid \; \mid^{\alpha} 1 ;
\sigma_{2, \chi_{1_{F^*}}}) $ and $ \Lg(\mid \; \mid^{\alpha} \chi_{1_{F^*}}; \lambda'(\det) \St_{U(3)}), \; \alpha > 0 $]{
Let $ \alpha_1, \alpha_2 \geq 0, $ let $ \chi_{1_{F^*}} \in X_{1_{F^*}}. $ Figure 16 shows lines and points of reducibility
of the representation $ \mid \; \mid^{\alpha_1} 1 \times \mid \; \mid^{\alpha_2} \chi_{1_{F^*}} \rtimes \lambda'. \;
\Lg(\mid \; \mid^{\alpha} 1 ; \sigma_{1, \chi_{1_{F^*}}}) $ and $ \Lg( \mid \; \mid^{\alpha} 1 ;
\sigma_{2, \chi_{1_{F^*}}}) $ are unitary for $ 0 < \alpha \leq 1 $ and non-unitary for $ \alpha > 1. \;
\Lg(\mid \; \mid^{\alpha} \chi_{1_{F^*}}; \lambda'(\det) \St_{U(3)}) $ is non-unitary $ \forall \alpha > 0. $}

\label{figure16}

\end{figure}

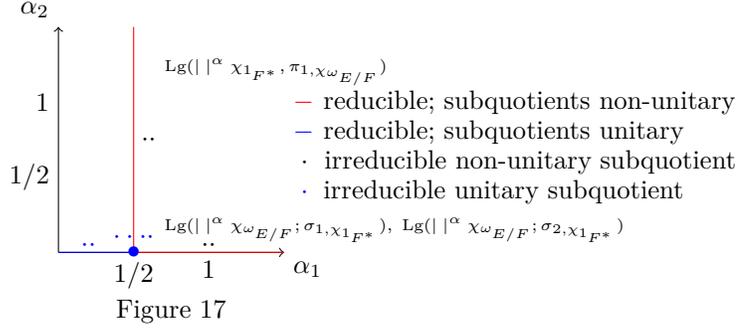
\begin{figure}
 
\begin{center}
\begin{tikzpicture}
 \draw[->] (0,0) -- (3,0);
\draw (3,0) node[below right] {$\alpha_1$};
\draw[->] (0,0) -- (0,3);
\draw (0,3) node[above left] {$ \alpha_2 $};

\draw [blue] (0,0) -- (1,0);

\draw [red] (1,0) -- (3,0);
\draw [red] (1,0) -- (1,3);

\draw (1,0) node[below] {$ 1/2 $};
\draw (0,1) node[left] {$ 1/2 $};
\draw (2,0) node[below] {$ 1 $};
\draw (0,2) node[left] {$ 1 $};
\draw (1,0) node[blue] {$ \bullet $};

\draw (1,0) node[blue, above] {$ \cdot \cdot \cdot \cdot $};
\draw (0.4, -0.1) node[blue, above] {$ \cdot \cdot $};
\draw (2, -0.1) node[above] { $ \cdot \cdot $};
\draw (1.2, 1.5) node { $ \cdot \cdot $};

\draw  (1.2, 2.4) node[right] { \begin{tiny} $ \Lg( \mid \; \mid^{\alpha} \chi_{1_{F^*}}, \pi_{1, \chi_{\omega_{E/F}}}) $ \end{tiny} };
\draw (1.2, 0.3) node[right] {\begin{tiny} $ \Lg( \mid \; \mid^{\alpha} \chi_{\omega_{E/F}} ; \sigma_{1, \chi_{1_{F^*}}}), \;
                        \Lg( \mid \; \mid^{\alpha} \chi_{\omega_{E/F}}; \sigma_{2, \chi_{1_{F^*}}}) $ \end{tiny} };

\draw (1.5, -0.5) node[below] { Figure 17 };
\draw (3, 2.0) node[red, right] { $ - $};
\draw (3.4, 2.0) node[right] {reducible; subquotients non-unitary};
\draw (3, 1.6) node[blue, right] { $ - $};
\draw (3.4, 1.6) node[right] {reducible; subquotients unitary};
\draw (3.1,1.2) node[right] { $ \cdot \; $ irreducible non-unitary subquotient };
\draw (3.1, 0.8) node[blue, right] { $ \cdot $};
\draw (3.4, 0.8) node[right] { irreducible unitary subquotient}; 
\end{tikzpicture}
\end{center}

\caption[$ \Lg( \mid \; \mid^{\alpha} \chi_{\omega_{E/F}} ; \sigma_{1, \chi_{1_{F^*}}}), \; \Lg( \mid \; \mid^{\alpha}
\chi_{\omega_{E/F}}; \sigma_{2, \chi_{1_{F^*}}}) $ and $ \Lg( \mid \; \mid^{\alpha} \chi_{1_{F^*}},
\pi_{1, \chi_{\omega_{E/F}}}), \; \alpha > 0 $]{Let $ \alpha_1, \alpha_2 \geq 0, $ let $ \chi_{\omega_{E/F}} 
\in X_{\omega_{E/F}} $ and $ \chi_{1_{F^*}} \in X_{1_{F^*}}. $  Figure 17
 shows lines and points of reducibility of
 the representation $ \mid \; \mid^{\alpha_1} \chi_{\omega_{E/F}} \times \mid \; \mid^{\alpha_2} \chi_{1_{F^*}} \rtimes
 \lambda'. \; \Lg( \mid \; \mid^{\alpha} \chi_{\omega_{E/F}} ; \sigma_{1, \chi_{1_{F^*}}}) $ and $ 
\Lg( \mid \; \mid^{\alpha} 
\chi_{\omega_{E/F}}; \sigma_{2, \chi_{1_{F^*}}}) $ are unitary for $ 0 < \alpha \leq 1/2 $ and non-unitary for
$ \alpha > 1/2. \; \Lg( \mid \; \mid^{\alpha} \chi_{1_{F^*}},
 \pi_{1, \chi_{\omega_{E/F}}}) $ is non-unitary $ \forall \alpha > 0.$}

\label{figure17}

\end{figure}

\smallskip

\subsection{Representations with cuspidal support in $ M_1 $}

\medskip

Recall that $ M_1 \cong \GL(2,E) \times E^1. $

\bigskip

\begin{center} $ \Lg(\mid \; \mid^{\alpha} \pi ; \lambda'), \alpha > 0, \; \pi $ a cuspidal unitary representation of
 $ \GL_2(E) $ \end{center}

\bigskip

Let $ \pi $ be a cuspidal unitary representation of $ \GL_2(E). $ Let $ \alpha > 0. $
Assume $ \exists \; g \in \GL(2,E) $ s. t. $ \pi(g) \neq \pi((\overline{g}^t)^{-1}). $ Then the induced representations $ \pi \rtimes \lambda' $ and 
$ \mid \; \mid^{\alpha} \pi \rtimes \lambda' $ are irreducible. $ \pi \rtimes \lambda' $ is unitary. $
\mid \; \mid^{\alpha} \pi \rtimes \lambda' $ is not hermitian $ \forall \alpha > 0, $ hence not unitary.

\bigskip

Assume $ \pi(g) = \pi((\overline{g}^t)^{-1}) \; \forall g \in \GL_2(E). $ Then $ \pi $ is obtained by base change lift from
$ U(2) $ to $ \GL(2,E), $ that is by 
endoscopic 
liftings from endoscopic data of $ U(2) $ to data of $ \GL_2(E) $ \cite{MR1081540}.

\medskip

Let $ \G:= U(2) $  and $ \overset{\sim}{\G} = \Res_{E/F}\G = \GL(2,E). $

Let $ \chi_{\omega_{E/F}} \in X_{\omega_{E/F}}. $ Let $ \Langdual{\G} $ be the Langlands dual group of $ \G. $ 
Recall that $ \sigma $ is defined to be the non-trivial element of $ \Gal(E,F), $ i.e. $ \sigma(x) = \overline{x} \; \forall
x \in E. $ Let $ \overset{\sim}{\sigma} $
denote the $ F - $ automorphism of $ \overset{\sim}{\G} $ associated to $ \sigma $ by the $ F - $ structure of
$ \overset{\sim}{\G}. $
Let $ \Gamma $ denote the absolute Galois group of $ E, $ let $ W_F $ and $ W_E $ denote the Weil groups of $ F $ and $ E, $
respectively. Let $ \rho_{\G} $ denote an L-action of $ \Gamma $ on $ \G $ and let $ \rho_{\overset{\sim}{\G}} $ denote an
L-action of $ \Gamma $ on $ \overset{\sim}{\G}. $
One fixes $ \omega_{\sigma} \in W_F \backslash W_E. $

\begin{Lemma}[\cite{MR1081540}, 4.7]
Up to isomorphisms, the base change problems for $ U(2) $ are the
endoscopic liftings from endoscopic data $ (\G, \Langdual{\G}, 1, \xi) $ and $ (\G, \Langdual{\G}, 1, \xi_{\chi_{\omega_{E/F}}}) $ for
$ (\overset{\sim}{\G}, \overset{\sim}{\sigma},1) $ to $ \overset{\sim}{\G}. $ Here

\medskip

$ \xi: \; \; \Langdual{\G} \backepsilon g \rtimes_{\rho_{\G}} \omega \mapsto (g,g) \rtimes_{\rho_{\overset{\sim}{\G}}} \omega \in
\Langdual{\overset{\sim}{\G}} $

\bigskip

$ \xi_{\chi_{\omega_{E/F}}}: \; \; \Langdual{\G} \backepsilon g \rtimes_{\rho_{\G}} \omega \mapsto \left[^{(g \chi_{\omega_{E/F}}(
\omega), g \chi_{\omega_{E/F}}(\omega)) \rtimes_{\rho_{\overset{\sim}{\G}}} \omega \in \Langdual{\overset{\sim}{\G}} \; \;
 \text{if} \; \omega \in W_E}_{(g, -g) \rtimes_{\rho_{\overset{\sim}{\G}}} \omega_{\sigma} \in \Langdual{\G} \; \; \text{if} \; 
\omega = \omega_{\sigma}} \right. $
\end{Lemma}

\bigskip

$ \xi $ is called standard base change and $ \xi_{\chi_{\omega_{E/F}}} $ is called twisted base change (\cite{Ko}).

\bigskip
Let $ \Pi_{temp}(\G) $ be the set of equivalence classes of irreducible admissible tempered 
representations of $ \G. $
Let $ \Pi_{temp}(\overset{\sim}{\G}) $ be the set of equivalence classes of irreducible admissible tempered 
representations of $ \overset{\sim}{\G}. $
Let $ \Pi $ be a tempered L-packet of $ \G, $ then $ \xi(\Pi), \; \xi_{\chi_{\omega_{E/F}}}(\Pi) \in 
\Pi_{temp}(\overset{\sim}{\G}). $

As before let $ \pi $ be a cuspidal unitary representation of $ \GL(2,E). $ If $ \pi(g) = \pi((\overline{g}^t)^{-1}) \; \forall g \in \GL(2,E), $ then $ \pi = \xi_{\chi_{\omega_{E/F}}}(\Pi) $ or
$ \pi = \xi(\Pi) $ (\cite[4.2]{MR1081540}).

\bigskip

Let  $ \pi $ be a cuspidal unitary representation of $ \GL_2(E). $

\medskip

(1) \; If $ \pi = \xi_{\chi_{\omega_{E/F}}}(\Pi), $ then $ \pi \rtimes \lambda' $ is reducible  
(\cite[4.2]{Ko}; \cite[6.2]{MR1266747}).
$  \pi \rtimes \lambda' = \tau_1(\pi) + \tau_2(\pi), $ where $ \tau_1(\pi) $ and $ \tau_2(\pi) $
are irreducible tempered.

$  \mid \;
\mid^{\alpha} \pi \rtimes \lambda' $ is irreducible and never unitarisable for $ \alpha > 0 $
(\cite[6.3]{MR1266747}).

\bigskip

(2) \; If $ \pi = \xi(\Pi), $ then $ \pi \rtimes \lambda' $ is irreducible ( \cite[4.2]{Ko},
; \cite[6.2]{MR1266747}).

\medskip

By results of Goldberg (\cite[6.3]{MR1266747}) one has:

(a) $ \mid \; \mid^{\alpha} \pi \rtimes \lambda' $ is irreducible and unitarisable
for $ 0 < \alpha < 1/2. $

(b) $ \mid \; \mid^{1/2} \pi \rtimes \lambda' $ is reducible. One has
 $ \mid \; \mid^{1/2} \pi \rtimes \lambda' = \sigma + \Lg(\mid \; \mid^{1/2} \pi ; \lambda'),$ where
$ \sigma $ is
a generic, non-supercuspidal
and square-integrable subrepresentation, and $ \Lg(\mid \; \mid^{1/2} \pi ; \lambda') $ is unitary.

(c) $ \mid \; \mid^{\alpha} \pi \rtimes \lambda' $ is irreducible and never unitarisable for
$ \alpha > 1/2. $

\bigskip

We obtain the following

\begin{Theorem}
Let $ \pi = \xi_{\chi_{\omega_{E/F}}}(\Pi), $ let $ \alpha > 0. $ Then

\smallskip
$ \Lg( \mid \; \mid^{\alpha} \pi; \lambda') $ is non-unitary.

\smallskip

$  \pi \rtimes \lambda' = \tau_1(\pi) + \tau_2(\pi), $ where $ \tau_1(\pi) $ and $ \tau_2(\pi) $ are irreducible tempered.

\bigskip

Let $ \pi = \xi(\Pi). $

Let $ 0 < \alpha < 1/2. \; \Lg( \mid \; \mid^{\alpha} \pi, \lambda') $ is unitary.

Let $ \alpha = 1/2. \; \Lg( \mid \; \mid^{1/2} \pi, \lambda') $ is unitary.

Let $ \alpha > 1/2. \;  \Lg( \mid \; \mid^{\alpha} \pi, \lambda') $ is non-unitary.

\smallskip

$ \pi \rtimes \lambda' $ is irreducible and unitary.

\end{Theorem}

We have no classification in terms of the Langlands-quotients for the induced representations of $ U(5) $ with cuspidal
support in $ M_2 \cong E^* \times U(3). $ One reason is that we do not have any parametrizatien for the cuspidal 
representations of $ U(3). $

%% file: amsart2.bbl
\providecommand{\bysame}{\leavevmode\hbox to3em{\hrulefill}\thinspace}
\providecommand{\MR}{\relax\ifhmode\unskip\space\fi MR }
\providecommand{\MRhref}[2]{%
  \href{http://www.ams.org/mathscinet-getitem?mr=#1}{#2}
}
\providecommand{\href}[2]{#2}